\documentclass[10pt]{article}
\usepackage[leqno]{amsmath}
\usepackage{amsfonts}
\usepackage{graphicx}

\usepackage{amsmath}
\usepackage{amssymb}
\usepackage{latexsym}
\usepackage{amsmath, amsfonts,amssymb, amsthm, euscript,makeidx,color,mathrsfs}

\def\sqr#1#2{{\vcenter{\vbox{\hrule height.#2pt
              \hbox{\vrule width.#2pt height#1pt \kern#1pt \vrule width.#2pt}
              \hrule height.#2pt}}}}
\def\signed #1{{\unskip\nobreak\hfil\penalty50
              \hskip2em\hbox{}\nobreak\hfil#1
              \parfillskip=0pt \finalhyphendemerits=0 \par}}
\def\endpf{\signed {$\sqr69$}}

\def \Re{\mbox{\rm Re}}
\def\3n{\negthinspace \negthinspace \negthinspace }
\def\2n{\negthinspace \negthinspace }
\def\1n{\negthinspace }

\def\dbA{\mathbb{A}}
\def\dbB{\mathbb{B}}

\def\dbH{\mathbb{H}}

\def\dbK{\mathbb{K}}

\def\dbM{\mathbb{M}}

\def\dbP{\mathbb{P}}
\def\dbQ{\mathbb{Q}}
\def\dbR{\mathbb{R}}
\def\dbS{\mathbb{S}}

\def\dbX{\mathbb{X}}

\def\={\buildrel \triangle \over =}

\def\ds{\displaystyle}

\def\ns{\noalign{\ss}}
\def\a{\alpha}
\def\b{\beta}
\def\g{\gamma}
\def\d{\delta}
\def\e{\varepsilon}
\def\z{\zeta}
\def\k{\kappa}
\def\l{\lambda}
\def\m{\mu}

\def\si{\sigma}
\def\t{\tau}
\def\f{\varphi}
\def\th{\theta}
\def\o{\omega}

\def\G{\Gamma}
\def\D{\Delta}
\def\Th{\Theta}
\def\L{\Lambda}

\def\F{\Phi}
\def\O{\Omega}

\def\cD{{\cal D}}

\def\cG{{\cal G}}

\def\cL{{\cal L}}

\def\cU{{\cal U}}

\def\BA{{\bf A}}

\def\no{\noindent}

\def\ss{\smallskip}
\def\ms{\medskip}
\def\bs{\bigskip}
\def\q{\quad}
\def\qq{\qquad}
\def\hb{\hbox}

\def\lan{\big\langle}
\def\ran{\big\rangle}

\def\pa{\partial}
\def\h{\widehat}
\def\wt{\widetilde}

\def\cd{\cdot}
\def\cds{\cdots}

\def\ae{\hbox{\rm a.e.{ }}}

\def\span{\hbox{\rm span$\,$}}

\def\coh{\mathop{\overline{\rm co}}}
\def\deq{\triangleq}
\def\les{\leqslant}
\def\ges{\geqslant}
\def\Re{{\mathop{\rm Re}\,}}

\def\({\Big (}
\def\){\Big )}
\def\[{\Big[}
\def\]{\Big]}
\def\bde{\begin{definition}}
\def\ede{\end{definition}}
\def\be{\begin{equation}}
\def\bel{\begin{equation}\label}
\def\ee{\end{equation}}
\def\bt{\begin{theorem}}
\def\et{\end{theorem}}
\def\bc{\begin{corollary}}
\def\ec{\end{corollary}}
\def\bl{\begin{lemma}}
\def\el{\end{lemma}}
\def\bp{\begin{proposition}}
\def\ep{\end{proposition}}
\def\bas{\begin{assumption}}
\def\eas{\end{assumption}}
\def\br{\begin{remark}}
\def\er{\end{remark}}
\def\ba{\begin{array}}
\def\ea{\end{array}}
\def\ed{\end{document}}

\def\square#1{\vbox{\hrule\hbox{\vrule height#1%
     \kern#1\vrule}\hrule}}
\def\rectangle#1#2{\vbox{\hrule\hbox{\vrule height#1%
     \kern#2\vrule}\hrule}}

\font\tenbb=msbm10 \font\sevenbb=msbm7 \font\fivebb=msbm5

\newfam\bbfam
\scriptscriptfont\bbfam=\fivebb \textfont\bbfam=\tenbb
\scriptfont\bbfam=\sevenbb

\newtheorem{lemma}{\hskip 1.4em Lemma}[section]
\newtheorem{remark}{\hskip 1.4em Remark}[section]

\newtheorem{theorem}{\hskip 1.4em Theorem}[section]
\newtheorem{corollary}{\hskip 1.4em Corollary}[section]

\newtheorem{definition}{\hskip 1.4em Definition}[section]
\newtheorem{proposition}{\hskip 1.4em Proposition}[section]
\newtheorem{assumption}{\hskip 1.4em Assumption}[section]

\makeatletter
   
   \@addtoreset{equation}{section}
\makeatother

\begin{document}

\title{\bf Forward-Backward Evolution Equations\\ and Applications\thanks{This work is supported in part by NSF Grant
DMS-1406776.}}
\author{Jiongmin Yong\thanks{Department of Mathematics, University of Central
Florida, Orlando, FL 32816, USA (jiongmin.yong@ucf.edu).}
 }

\maketitle

\no{\bf Abstract:} Well-posedness is studied for a special system of
two-point boundary value problem for evolution equations which is
called a {\it forward-backward evolution equation} (FBEE, for
short). Two approaches are introduced: A decoupling method with some
brief discussions, and a method of continuation with some
substantial discussions. For the latter, we have introduced Lyapunov
operators for FBEEs, whose existence leads to some uniform {\it a
priori} estimates for the mild solutions of FBEEs, which will be
sufficient for the well-posedness. For some special cases, Lyapunov
operators are constructed. Also, from some given Lyapunov operators,
the corresponding solvable FBEEs are identified.

\ms

\no\bf AMS Mathematics Subject Classification. \rm 34D10, 34G20,
35K90, 35L90, 47D06, 47J35, 49K27.

\ms

\no{\bf Keywords}. Forward-backward evolution equations, decoupling
field, Lyapunov operator, method of continuation.

\maketitle

\section{Introduction}

In this paper, we consider the following system of evolution
equations:
\bel{FBEE1}\left\{\2n\ba{ll}
\ns\ds\3n\ba{ll}
\ns\ds\dot y(t)=Ay(t)+b(t,y(t),\psi(t)),\\
\ns\ds\dot\psi(t)=-A^*\psi(t)-g(t,y(t),\psi(t)),\ea\qq t\in[0,T],\\
\ns\ds y(0)=x,\qq\psi(T)=h(y(T)),\ea\right.\ee
where $A:\cD(A)\subseteq X\to X$ generates a $C_0$-semigroup
$e^{At}$ on a real Hilbert space $X$ (identified with its dual
$X^*$), with
$$\big(e^{At}\big)^*=e^{A^{\1n*}t},\qq t\ge0,$$
being the adjoint semigroup generated by $A^*$ (the adjoint operator
of $A$), and $b$, $g$, and $h$ being suitable maps. The above could
be called a {\it two-point boundary value problem}, mimicking a
similar notion for ordinary differential equations. We see that the
equation for $y(\cd)$ is an initial value problem which should be
solved forwardly, and the equation for $\psi(\cd)$ is a terminal
value problem which should be solved backwardly. Therefore, inspired
by the so-called {\it forward-backward stochastic differential
equations} (FBSDEs, for short, see \cite{Ma-Yong 1999, Yong 1997,
Yong 2010} for details), we prefer to call (\ref{FBEE1}) a {\it
forward-backward evolution equation} (FBEE, for short). In common
occasions, two-point boundary value problem is related to certain
eigenvalue problems, for which the well-posedness of the problem
might not be the goal, instead, one might be more interested in the
existence of solutions, not necessarily the uniqueness. See, for
examples, \cite{Bailey-Shampine-Waltman 1968,Zettl
2005,DeCoster-Habets 2006}, and see also \cite{Dower-McEneaney 2015}
for some other considerations. Whereas, in this paper, we are
interested in the well-posedness of (\ref{FBEE1}). On the other
hand, our system has a special structure, involving one {\it
forward} evolution equation and one {\it backward} evolution
equation. Hence, we use the name FBEE to distinguish the current
situation from other situations in the literature.

\ms

A pair of functions $(y(\cd),\psi(\cd))$ is called a {\it strong
solution} of (\ref{FBEE1}) if these functions are differentiable
almost everywhere, with the property
$$y(t)\in\cD(A),\q\psi(t)\in\cD(A^*),\qq\ae t\in[0,T],$$
and the equations are satisfied almost everywhere. A pair
$(y(\cd),\psi(\cd))$ is called a {\it mild solution} (or a {\it weak
solution}) to FBEE (\ref{FBEE1}) if the following system of integral
equations are satisfied:
\bel{mild}\left\{\2n\ba{ll}
\ns\ds y(t)=e^{At}x+\int_0^te^{A(t-s)}b(s,y(s),\psi(s))ds,\\
\ns\ds\psi(t)=e^{A^*(T-t)}h(y(T))+\int_t^Te^{A^*(s-t)}g(s,y(s),\psi(s))ds,
\ea\right.\q t\in[0,T].\ee
Note that in the case $A$ is bounded, (\ref{FBEE1}) and (\ref{mild})
are actually equivalent, and thus, a mild solution
$(y(\cd),\psi(\cd))$ is actually a strong solution.

\ms

Our study of the above system is mainly motivated by the study of
optimal control theory. It is known that for a standard optimal
control problem of an evolution equation with, say, a Bolza type
cost functional, by applying the Pontryagin maximum/minimum
principle, one will obtain an {\it optimality system} of the above
form whose solution will give a candidate for the optimal trajectory
and its adjoint (\cite{Li-Yong 1995}). Therefore, solvability of the
above type system is important, at least for optimal control theory
of evolution equations.

\ms

Roughly speaking, when $T$ is small enough, or the Lipschitz
constants of the involved functions are small enough, one can show
that FBEE (\ref{FBEE1}) will have a unique mild solution, by means
of contraction mapping theorem. On the other hand, if (\ref{FBEE1})
is the optimality system (obtained via Pontryagin maximum/minimum
principle) of a corresponding optimal control problem which admits
an optimal control, then this FBEE admits a mild solution, which
might not be unique. Further, if the corresponding optimal control
has an optimal control and the optimality system admits a unique
mild solution, then this solution can be used to construct the
optimal control(s). Hence, under proper conditions, FBEE
(\ref{FBEE1}) could admit a (unique) mild solution, without
restriction on the length of the time horizon $T$, and/or the size
of the Lipschitz constants of the involved functions. This is
actually the case if the FBEE is the optimality system of a
linear-quadratic (LQ, for sort) optimal control problem satisfying
proper conditions (\cite{Li-Yong 1995}).

\ms

In this paper, we will study the (unique) solvability of FBEE
(\ref{FBEE1}) under some general conditions. Two approaches will be
introduced: decoupling method and method of continuation. The former
is inspired by the so-called {\it invariant embedding} which can be
traced back to \cite{Ambarzumyan 1943, Chandrasekhar 1950,
Bellman-Wing 1975}. Such a method was used in the study of FBSDEs
(see \cite{Ma-Protter-Yong 1994, Ma-Yong 1999}, for details). The
latter is inspired by the method of continuity for elliptic partial
differential equations (see, e.g. \cite{Gilbarg-Trudinger}), and
FBSDEs (\cite{Hu-Peng 1995, Yong 1997, Peng-Wu 1999, Yong 2010}).
Due to the nature of FBEE (\ref{FBEE1}), some technical difficulties
exist in applying either of these methods. We will briefly present
some main idea of the decoupling method and will relatively more
carefully present the method of continuation.

\ms

The rest of this paper is organized as follows. In Section 2, we
will present some preliminary results, including a main motivation
from optimal control theory. Linear FBEEs are carefully discussed in
Section 3. In Section 4, a brief description on the decoupling
method will be given. In Section 5, we will introduce the so-called
{\it Lyapunov operator} which is adopted from \cite{Yong 2010} (for
FBSDEs). The existence of Lyapunov operators lead to some uniform
{\it a priori} estimates for the mild solutions of our FBEE.
Well-posedness of FBEEs will be established in Section 6. In Section
7, we will construct some Lyapunov operators through which some
well-posed FBEEs will be identified. In Section 8, we briefly
discussed some extensions of our main results. In Section 9, several
illustrative examples will be presented. Finally, some concluding
remarks will be made in Section 10.

\section{Preliminaries}

Throughout of this paper, we let $X$ be a separable real Hilbert
space, with the norm $\|\cd\|$ and the inner product
$\lan\cd\,,\cd\ran$. We identify the dual $X^*$ with $X$. The set of
all bounded linear operators from $X$ to itself is denoted by
$\cL(X)$. The set of all self-adjoint operators on $X$ is denoted by
$\dbS(X)$ and the set of all positive semi-definite operators on $X$
is denoted by $\dbS^+(X)$. For the notational simplicity, when there
is no confusion, we will not distinguish between $\l$ and $\l I$
(for any $\l\in\dbR$). For example, we use $\l-A$ to denote $\l
I-A$. Also, if $F$ is in $\dbS^+(X)$, we denote it by $F\ges0$; if
$F-cI\ges0$, we simply denote it by $F\ges c$, and $F\les c$ means
$-F\ges-c$. Next, we denote
$$C([0,T];X)=\Big\{y:[0,T]\to X\bigm|y(\cd)\hb{ is
continuous}\Big\},$$
and
$$\|y(\cd)\|_\infty=\sup_{t\in[0,T]}\|y(t)\|,\qq\forall y(\cd)\in
C([0,T];X).$$

\ms

For convenience and definiteness of our presentation, we introduce
the following standing assumptions:

\ms

{\bf(H0)$'$} $A:\cD(A)\subseteq X\to X$ generates a $C_0$-semigroup
$e^{At}$ on $X$.

\ms

{\bf(H0)} In addition to (H0)$'$, either
\bel{case1}A^*=A,\ee
with the spectrum $\si(A)\subseteq\dbR$ of $A$ satisfying
\bel{si>0}\sup\si(A)\equiv\sup\Re\si(A)=-\si_0<0,\ee
or
\bel{case2}A^*=-A,\ee
for which it holds: $\si(A)\subseteq i\dbR$ and thus
\bel{si=0}\sup\Re\si(A)=\si_0=0.\ee

\ms

Case (\ref{case1}) corresponds to the heat equation (or second order
parabolic equations) with proper lower order terms and proper
boundary conditions. Case (\ref{case2}) corresponds to the wave
equation (or second order hyperbolic equations) with proper boundary
conditions, without damping. Some extensions of the results
presented in this paper are possible. But for the moment, we prefer
not to get into the most general situations, for the simplicity of
our presentation. We should keep in mind that for the case $A^*=A$,
one has $\si_0>0$ and for the case $A^*=-A$, one has $\si_0=0$.

\ms

Let us now look at our main motivation of studying our FBEEs.
Consider the following controlled system:
\bel{state2.5}\left\{\2n\ba{ll}
\ns\ds
\dot y(t)=Ay(t)+f(t,y(t),u(t)),\qq t\in[0,T],\\
\ns\ds y(0)=x,\ea\right.\ee
with cost functional of Bolza type:
\bel{2.2}J(x;u(\cd))=\int_0^Tf^0(s,y(s),u(s))ds+f^1(y(T)).\ee
In the above, $f:[0,T]\times X\times U\to X$, $f^0:[0,T]\times
X\times U\to\dbR$, $f^1:X\to\dbR$ are suitable maps, with $U$ being
a separable metric space. We call $x\in X$ an {\it initial state},
$u(\cd)$ a {\it control}, and $y(\cd)$ a {\it state trajectory},
respectively. Denote
$$\cU=\big\{u:[0,T]\to U\bigm|u(\cd)\hb{ is measurable}\big\}.$$
This is the set of all {\it admissible controls}. Under some mild
conditions, for any $x\in X$ and $u(\cd)\in\cU$, {\it state
equation} (\ref{state2.5}) admits a unique {\it mild solution}
$y(\cd)\equiv y(\cd\,;x,u(\cd))$, i.e., the solution to the
following integral equation:
\bel{2.3}y(t)=e^{At}x+\int_0^te^{A(t-s)}f(s,y(s),u(s))ds,\qq
t\in[0,T],\ee
and the cost functional $J(x;u(\cd))$ is well-defined. Then one can
pose the following optimal control problem.

\ms

\bf Problem (C). \rm For any initial state $x\in X$, find a $\bar
u(\cd)\in\cU$ such that
\bel{2.4}J(x;\bar u(\cd))=\inf_{u(\cd)\in\cU}J(x;u(\cd)).\ee

\ms

Any $\bar u(\cd)\in\cU$ satisfying (\ref{2.4}) is called an {\it
optimal control}, the corresponding $\bar y(\cd)\equiv
y(\cd\,;x,\bar u(\cd))$ is called an {\it optimal state trajectory}
and $(\bar y(\cd),\bar u(\cd))$ is called an {\it optimal pair}.

\ms

With the above setting, we have the following standard result. To
simplify the presentation, we assume that the involved maps
$f,f^0,f^1$ have all the required measurability and smoothness. The
readers are referred to \cite{Li-Yong 1995} for details.

\ms

\bf Proposition 2.1. (Pontryagin's Minimum Principle) \sl Let
{\rm(H0)$'$} hold and let $(\bar y(\cd),\bar u(\cd))$ be an optimal
pair of Problem {\rm(C)}. Then the following {\it minimum condition}
holds:
\bel{min2.9}\ba{ll}
\ns\ds\lan\psi(t),f(t,y(t),u(t))\ran+f^0(t,y(t),u(t))=\min_{u\in
U}\big[\lan\psi(t),f(t,y(t),u)\ran+f^0(t,y(t),u)\big],\qq
t\in[0,T],\ea\ee
where $\psi(\cd)$ is the mild solution to the following {\it adjoint
equation}:
\bel{adjoint2.10}\left\{\2n\ba{ll}
\ns\ds\dot\psi(t)=-A^*\psi(t)-f_y(t,\bar y(t),\bar u(t))^*\psi(t)-f^0_y(t,
\bar y(t),\bar u(t)),\q t\in[0,T],\\
\ns\ds\psi(T)=f^1_y(\bar y(T)),\ea\right.\ee
i.e., the following holds:
\bel{}\ba{ll}
\ns\ds\psi(t)=e^{A^*(T-t)}f_y^1(\bar y(T))+\int_t^T
e^{A^*(s-t)}\[f_y(s,\bar y(s),\bar u(s))^*\psi(s)+f^0_y(s,\bar
y(s),\bar u(s))\]ds,\q t\in[0,T].\ea\ee

\rm

\ms

Note that (\ref{state2.5}) and (\ref{adjoint2.10}) form a system
with the minimum condition (\ref{min2.9}) bringing in the coupling.
Suppose there exists a map $\f:[0,T]\times X\times X\to U$ such that
$$\ba{ll}
\ns\ds\lan\psi,f(t,y,\f(t,y,\psi))\ran+f^0(t,y,\f(t,y,\psi))=\min_{u\in
U}\[\lan\psi,f(t,y,u)\ran+f^0(t,y,u)\].\ea$$
Then we obtain the following system (dropping the bar in $\bar
y(\cd)$)
\bel{2.7}\left\{\2n\ba{ll}
\ns\ds\dot y(t)=Ay(t)+f(t,y(t),\f(t,y(t),\psi(t))),\qq t\in[0,T],\\
[1mm]
\ns\ds\dot\psi(t)=-A^*\psi(t)-f_y(t,y(t),\f(t,y(t),\psi(t)))^*\psi(t)
-f^0_y(t,y(t),\f(t,y(t),\psi(t))),\qq t\in[0,T],\\
[1mm]
\ns\ds y(0)=x,\qq\psi(T)=f^1_y(y(T)).\ea\right.\ee
This is called the {\it optimality system} of Problem (C), which is
an FBEE of form (\ref{FBEE1}) with
$$\left\{\2n\ba{ll}
\ns\ds b(t,y,\psi)=f(t,y,\f(t,y,\psi)),\\
\ns\ds
g(t,y,\psi)=f_y(t,y,\f(t,y,\psi))^*\psi+f^0_y(t,y,\f(t,y,\psi)),\\
\ns\ds h(y)=f^1_y(y).\ea\right.$$
If $(y(\cd),\psi(\cd))$ is a mild solution of FBEE (\ref{2.7}), then
$y(\cd)$ will be a candidate of optimal trajectory and
$\f(\cd\,,y(\cd),\psi(\cd))$ will be a candidate of optimal control.

\ms

Let us now look an interesting special case of the above. To this
end, we let $U$ also be a real Hilbert space, and
$$\left\{\2n\ba{ll}
\ns\ds f(t,y,u)=F(t,y)+B(t)u,\qq f^1(y)=G(y),\\
\ns\ds f^0(t,y,u)=Q(t,y)+\lan S(t)y,u\ran+{\,1\over\,2}\lan
R(t)u,u\ran,\ea\right.$$
for some suitable maps $F(\cd\,,\cd)$, $B(\cd)$, $G(\cd)$,
$Q(\cd\,,\cd)$, $S(\cd)$, and $R(\cd)$. Then the state equation
becomes
\bel{2.8}\left\{\2n\ba{ll}
\ns\ds\dot y(t)=Ay(t)+F(t,y(t))+B(t)u(t),\qq t\in[0,T],\\
\ns\ds y(0)=x,\ea\right.\ee
and the cost functional takes the form:
\bel{2.9}\ba{ll}
\ns\ds J(x;u(\cd))=\int_0^T\[Q(y(t))+\lan
S(t)y(t),u(t)\ran+{1\over2}\lan R(t)u(t),u(t)\ran\]dt+G(y(T)).\ea\ee
Note that the right-hand side of the state equation is affine in
$u(\cd)$ and the integrand in the cost functional is up to quadratic
in $u(\cd)$. We therefore refer to the corresponding optimal control
problem as an {\it affine-quadratic optimal control problem} (AQ
problem, for short). For finite-dimensional case, general AQ problem
was studied in \cite{Wang-Yong 2014}. In current case, the adjoint
equation reads
$$\left\{\2n\ba{ll}
\ns\ds\dot\psi(t)=-A^*\psi(t)-F_y(t,y(t))^*\psi(t)-Q_y(t,y(t))-S(t)^*u(t),
\q t\in[0,T],\\
\ns\ds\psi(T)=G_y(y(T)).\ea\right.$$
By the minimum condition
$$\ba{ll}
\ns\ds\lan\psi(t),B(t)u(t)\ran+\lan
S(t)y(t),u(t)\ran+{1\over2}\lan R(t)u(t),u(t)\ran\\
\ns\ds=\min_{u\in U}\[\lan\psi(t),B(t)u\ran+\lan
S(t),u\ran+{1\over2}\lan R(t)u,u\ran\],\qq t\in[0,T],\ea$$
we obtain, assuming the invertibility of $R(t)$,
$$u(t)=-R(t)^{-1}\big[B(t)^*\psi(t)+S(t)y(t)\big],\qq t\in[0,T].$$
Therefore, the corresponding optimality system reads as follows:
$$\left\{\2n\ba{ll}
\ns\ds\dot y(t)=Ay(t)+F(t,y(t))-B(t)R(t)^{-1}S(t)y(t)-B(t)R(t)^{-1}B(t)^*\psi(t),\\
[1mm]
\ns\ds\dot\psi(t)=-A^*\psi(t)-\big[Q_y(t,y(t))-S(t)^*R(t)^{-1}S(t)y(t)\big]
-\big[F_y(t,y(t))^*-S(t)^*R(t)^{-1}B(t)^*\big]\psi(t),\\
[1mm]
\ns\ds y(0)=x,\qq\psi(T)=G_y(y(T)).\ea\right.$$
When
$$\left\{\2n\ba{ll}
\ns\ds y\mapsto F(t,y)\hb{ is linear},\\
\ns\ds y\mapsto Q(t,y),~y\mapsto G(y)\hb{ are convex},\ea\right.$$
the corresponding optimal control problem is referred to as {\it
linear-convex problem}, which was studied in \cite{You 1987a, You
1987b, You 1997}. See also \cite{You 2002} for some investigations
on finite-dimensional two-person zero-sum differential games of
linear state equation with non-quadratic payoff/cost functional
where the convexity of $y\mapsto Q(t,y)$ and $y\mapsto G(y)$ were
not assumed. Further, if
$$F(t,y)\equiv0,\q Q(t,y)={1\over2}\lan Q(t)y,y\ran,\q G(y)
={1\over2}\lan Gy,y\ran,$$
for some $Q:[0,T]\to\dbS(X)$ and $G\in\dbS(X)$, the problem is
reduced to a classical LQ problem. In this case, the optimality
system becomes the following linear FBEE:
\bel{FBEE-LQ}\left\{\2n\ba{ll}
\ns\ds\dot
y(t)=\big[A-B(t)R(t)^{-1}S(t)\big]y(t)-B(t)R(t)^{-1}B(t)^*\psi(t),\\
[1mm]
\ns\ds\dot\psi(t)=-\big[A\1n-\1n
B(t)R(t)^{-1}S(t)\big]^*\psi(t)\1n-\1n\big[Q(t)\1n-\1n
S(t)^*R(t)^{-1}S(t)\big]y(t),\\ [1mm]
\ns\ds y(0)=x,\qq\psi(T)=Gy(T).\ea\right.\ee
It is known that under the following conditions:
\bel{}R(t)\ges\d I,\q Q(t)-S(t)^*R(t)^{-1}S(t)\ges0,\q t\in[0,T],\q
G\ges0,\ee
the map $u(\cd)\mapsto J(x;u(\cd))$ for the current LQ problem is
uniformly convex, and the above linear FBEE (\ref{FBEE-LQ}) admits a
unique mild solution $(y(\cd),\psi(\cd))$ (\cite{Li-Yong 1995}).

\ms

Next, we note that under (H0)$'$, by Hille-Yosida's theorem, there
exist $M\ges1$ and $\o\in\dbR$ such that
\bel{H-Y1}\|(\l-A)^{-n}\|\les{M\over(\l-\o)^n},\qq\forall\l>\o,~n\ges1,\ee
and the {\it Yosida approximation} $A_\l$ of $A$ is well-defined:
\bel{Yosida}A_\l=\l A(\l-A)^{-1},\qq\qq\l>\o.\ee
By making a shifting and absorbing a relevant term into
$b(t,y,\psi)$ (see (\ref{FBEE1})), we may assume that $\o=0$ in the
above. Then by \cite{Pazy 1983}, we may assume the following:
\bel{2.19}\left\{\2n\ba{ll}
\ns\ds\|e^{A_\l t}\|\les M,\qq\forall t\ges0,\\
\ns\ds\lim_{\l\to\infty}\|A_\l x-Ax\|=0,\qq\forall x\in\cD(A),\\
\ns\ds\lim_{\l\to\infty}\sup_{t\in[0,T]}\|e^{A_\l
t}x-e^{At}x\|=0,\qq\forall x\in X.\ea\right.\ee

\ms

Now, let us look at (H0). It is clear that under condition
(\ref{case1}), one has
\bel{A<-d}\lan Ax,x\ran\les-\si_0\|x\|^2,\qq\forall x\in\cD(A),\ee
and under (\ref{case2}), one has
\bel{A=0}\lan Ax,x\ran=0,\qq\forall x\in\cD(A).\ee
The following simple result is concerned with the Yosida
approximation $A_\l$ of $A$, under (H0).

\ms

\bf Proposition 2.2. \sl If $(\ref{case1})$ holds, then
\bel{3.5}\lan A_\l
x,x\ran\les-{\l\si_0\over\l+\si_0}\|x\|^2,\qq\forall x\in
X,\q\l>0.\ee
If $(\ref{case2})$ holds, then
\bel{3.6}\lan(A_\l+A_\l^*)x,x\ran\les0,\qq\forall x\in X,\q\l>0.\ee

\ms

\it Proof. \rm Under (\ref{case1}), $A$ admits the following
spectral decomposition (\cite{Dunford-Schwartz}):
\bel{spectral}Ax=\int_{\si(A)}\m dE_\m x,\qq\forall x\in\cD(A),\ee
where $\m\mapsto E_\m$ is the projection-valued measure associated
with $A$, and $\si(A)\subseteq(-\infty,-\si_0]$ is the spectrum of
$A$.
Consequently,
$$A_\l=\l A(\l-A)^{-1}=\int_{\si(A)}{\l\m\over\l-\m}dE_\m.$$
Since the map $\m\mapsto{\l\m\over\l-\m}$ is increasing on
$(-\infty,-\si_0]$, we have
$$\lan A_\l x,x\ran=\int_{\si(A)}{\l\m\over\l-\m}d\|E_\m x\|^2\les-{\l\si_0\over\l+\si_0}
\|x\|^2,\qq\forall x\in X.$$
Now, let (\ref{case2}) hold. We let $\dbX=X+iX$ be the
complexification of $X$, i.e.,
$$\dbX=\big\{x+iy\bigm|x,y\in X\big\},$$
with the following definition of addition, scalar multiplication,
and inner product:
$$\ba{ll}
\ns\ds(x+iy)+(\wt x+i\wt y)=(x+\wt x)+i(y+\wt y),\qq\forall x,y,\wt
x,\wt y\in X,\\
\ns\ds(\a+i\b)(x+iy)=(\a x-\b y)+i(\a y+\b
x),\qq\forall\a,\b\in\dbR,~x,y\in X,\\
\ns\ds\lan x+iy,\wt x+i\wt y\ran=\lan x,\wt x\ran+\lan y,\wt
y\ran+i\big(\lan y,\wt x\ran-\lan x,\wt y\ran\big),\qq\forall x,\wt
x,y,\wt y\in X.\ea$$
Naturally extend $A$ to $\BA:\cD(\BA)\subseteq\dbX\to\dbX$ as
follows
$$\left\{\2n\ba{ll}
\ns\ds\cD(\BA)=\cD(A)+i\cD(A)\subseteq\dbX,\\
\ns\BA(x+iy)=Ax+iAy,\qq\forall x+iy\in\cD(\BA).\ea\right.$$
Then under (\ref{case2}), we have
$$\ba{ll}
\ns\ds\lan\BA(x+i y),\wt x+i\wt y\ran=\lan Ax+iAy,\wt x+i\wt
y\ran=\lan Ax,\wt x\ran+\lan Ay,\wt y\ran+i\big(\lan Ay,\wt
x\ran-\lan Ax,\wt y\ran\big)\\
\ns\ds=-\lan x,A\wt x\ran-\lan y,A\wt y\ran-i\big(y,A\wt x\ran-\lan
x,A\wt y\big)=-\lan\wt x+i\wt y,\BA(x+iy)\ran,\ea$$
which implies that
$$\BA^*=-\BA.$$
Hence,  $\BA$ admits the following spectral decomposition:
$$\BA z=\int_{\si(\BA)}\m dE_\m z,\qq\forall z\in\cD(\BA),$$
with $\si(\BA)\subseteq i\dbR$. Consequently, for any $\l>0$,
$$\ba{ll}
\ns\ds\lan(\BA_\l+\BA_\l^*)z,z\ran
=\int_{\si(\BA)}\[{\l\m\over\l-\m}-{\l\m\over\l+\m}\]d\|E_\m
z\|^2\\
\ns\ds\qq\qq\qq=-\int_{\si(\BA)}{2\l|\m|^2\over\l^2+|\m|^2}\,d\|E_\m
z\|^2\les0,\qq z=x+iy\in\dbX.\ea$$
Note that for any $x\in X$, if
$$(\l-\BA)^{-1}x=\wt z=\wt x+i\wt y,$$
with $\wt x,\wt y\in X$, then
$$x=(\l-\BA)(\wt x+i\wt y)=(\l-A)\wt x+i(\l-A)\wt y.$$
Hence, we must have
$$\wt x=(\l-A)^{-1}x,\qq\wt y=0.$$
Consequently,
$$\BA_\l x=A_\l x,\qq\forall x\in X.$$
Likewise,
$$\BA_\l^*x=A^*_\l x,\qq\forall x\in X.$$
Hence, (\ref{3.6}) follows. \endpf

\ms

To conclude this section, let us introduce some assumptions on the
coefficients of FBEE (\ref{FBEE1}).

\ms

{\bf(H1)} The maps $b,g:[0,T]\times X\times X\to X$ and $h:X\to X$
are continuous, and the map
$(y,\psi)\mapsto(b(t,y,\psi),g(t,y,\psi),h(y))$ is locally
Lipschitz.

\ms

{\bf(H2)} In addition to (H1), let the map
$(y,\psi)\mapsto(b(t,y,\psi),g(t,y,\psi),h(y))$ be uniformly
Lipschitz and of uniformly linear growth, i.e., there exists a
constant $L>0$ such that
\bel{2.1}\left\{\2n\ba{ll}
\ns\ds\|b(t,0,0)\|\les L,\qq\qq t\in[0,T],\\ [1mm]
\ns\ds\|b(t,y,\psi)-b(t,\bar y,\bar\psi)\|\les L\|y-\bar
y\|+L\|\psi-\bar\psi\|,\qq\forall t\in[0,T],~y,\bar
y,\psi,\bar\psi\in X,\ea\right.\ee
\bel{2.2}\left\{\2n\ba{ll}
\ns\ds\|g(t,0,0)\|\les L,\qq\qq t\in[0,T],\\ [1mm]
\ns\ds\|g(t,y,\psi)-g(t,\bar y,\bar\psi)\|\les L\|y-\bar
y\|+L\|\psi-\bar\psi\|,\qq\forall t\in[0,T],~y,\bar
y,\psi,\bar\psi\in X,\ea\right.\ee
and
\bel{}\|h(y)-h(\bar y)\|\les L\|y-\bar y\|,\qq\forall y,\bar y\in
X.\ee

\ms

{\bf(H3)} In addition to (H1), let the map
$(y,\psi)\mapsto(b(t,y,\psi),g(t,y,\psi),h(y))$ be Fr\'echet
differentiable, with continuous Fr\'echet derivatives.

\ms

Note that (H3) is neither stronger nor weaker than (H2), since the
Fr\'echet derivatives $b_x,b_\psi,g_x,g_\psi,h_y$, if they exist,
are not necessarily uniformly bounded. We let $\cG_1$, $\cG_2$,
$\cG_3$ be the set of all $(b,g,h)$ satisfying (H1), (H2), and (H3),
respectively. Any $(b,g,h)\in\cG_1$ uniquely {\it generates} an FBEE
(1.1) (without mentioning the well-posedness). Hence, any
$(b,g,h)\in\cG_1$ is called the {\it generator} of an FBEE of form
(\ref{FBEE1}).

\ms

\section{Linear FBEEs}

In this section, we consider the following linear FBEE:
\bel{linear}\left\{\2n\ba{ll}
\ns\ds\3n\ba{ll}\dot y(t)=Ay(t)+B_{11}(t)y(t)+B_{12}(t)\psi(t)+b_0(t),\\
\ns\ds\dot\psi(t)=-A^*\psi(t)-B_{21}(t)y(t)-B_{22}(t)\psi(t)-g_0(t),\ea\q t\in[0,T],\\
\ns\ds y(0)=x,\qq\psi(T)=Hy(T)+h_0,\ea\right.\ee
with
\bel{Bbg}\left\{\2n\ba{ll}
\ns\ds B_{ij}(\cd)\in L^\infty(0,T;\cL(X)),\qq i,j=1,2,\\
\ns\ds b_0(\cd),g_0(\cd)\in L^1(0,T;X),\q H\in\cL(X),\q h_0\in
X.\ea\right.\ee
The above is a special case of (\ref{FBEE1}). A pair
$(y(\cd),\psi(\cd))$ is called a mild solution to (\ref{linear}) if
the following holds:
$$\left\{\2n\ba{ll}
\ns\ds
y(t)=e^{At}x+\int_0^te^{A(t-s)}\big[B_{11}(s)y(s)+B_{12}(s)\psi(s)+b_0(s)\big]ds,\\
\ns\ds\psi(t)=e^{A^*(T-t)}\big[Hy(T)+h_0\big]+\int_t^Te^{A^*(s-t)}\big[B_{21}(s)y(s)
+B_{22}(s)\psi(s)+g_0(s)\big]ds,\ea\q t\in[0,T].\right.$$
Our first result is the
following.

\ms

\bf Proposition 3.1. \sl Let {\rm(H0)$'$} and {\rm(\ref{Bbg})} hold.
Then FBEE {\rm(\ref{linear})} admits a mild solution if the
following operator:
$$\ba{ll}
\ns\ds\psi(\cd)\mapsto\psi(\cd)-\int_0^T\(\F_{22}(T,\cd)H\F_{11}(T,s)
+\int_{s\vee\cd}^T\F_{22}(r,\cd)B_{21}(r)\F_{11}(r,s)dr\)B_{12}(s)\psi(s)ds\ea$$
is invertible on $C([0,T];X)$, where $\F_{11}(\cd\,,\cd)$ and
$\F_{22}(\cd\,,\cd)$ are evolution operators generated by
$A+B_{11}(\cd)$ and $A^*+B_{22}(\cd)$, respectively.

\ms

\it Proof. \rm By the variation of constants formula, we have
$$y(t)=\F_{11}(t,0)x+\int_0^t\F_{11}(t,s)\big[B_{12}(s)\psi(s)+b_0(s)\big]ds,$$
and
$$\ba{ll}
\ns\ds\psi(t)=\F_{22}(T,t)\big[Hy(T)+h_0\big]+\int_t^T\F_{22}(s,t)\big[B_{21}(s)y(s)
+g_0(s)\big]ds\\
\ns\ds\qq=\F_{22}(T,t)\Big\{H\[\F_{11}(T,0)x+\int_0^T\F_{11}(T,s)
\(B_{12}(s)\psi(s)+b_0(s)\)ds\]+h_0\Big\}\\
\ns\ds\qq\q+\int_t^T\F_{22}(s,t)\[B_{21}(s)\(\F_{11}(s,0)x+\int_0^s\F_{11}(s,r)
\big[B_{12}(r)\psi(r)+b_0(r)\big]dr\)+g_0(s)\]ds\\
\ns\ds\qq=\2n\int_0^T\2n\(\F_{22}(T,t)H\F_{11}(T,s)+\2n\int_{s\vee
t}^T\F_{22}(r,t)B_{21}(r)\F_{11}(r,s)dr\)B_{12}(s)\psi(s)ds\\
\ns\ds\qq\q+\(\F_{22}(T,t)H\F_{11}(T,0)+\int_t^T\F_{22}(s,t)B_{21}(s)\F_{11}(s,0)ds\)x
\\
\ns\ds\qq\q+\int_0^T\(\F_{22}(T,t)H \F_{11}(T,s)+\int_{s\vee
t}\F_{22}(r,t)B_{21}(r)\F_{11}(r,s)dr\)b_0(s)ds\\
\ns\ds\qq\q+\int_t^T\F_{22}(s,t)g_0(s)ds+\F_{22}(T,t)h_0,\qq
t\in[0,T].\ea$$
The above is a Fredholm integral equation for $\psi(\cd)$ of the
second kind. By our assumption, it has a unique solution. Then our
result follows. \endpf

\ms

Next, we consider a special case: $A^*=-A$. For such a case, we have
that
\bel{AA}\dbA\equiv\begin{pmatrix}A&0\\
0&-A^*\end{pmatrix}=\begin{pmatrix}A&0\\ 0&A\end{pmatrix}\ee
generates a $C_0$-group $e^{\dbA t}\equiv\begin{pmatrix}\sc
e^{At}&\sc0\\ \sc0&\sc e^{At}\end{pmatrix}$ on $X\times X$. Hence,
if we denote
\bel{BB}\dbB(t)=\begin{pmatrix}B_{11}(t)&B_{12}(t)\\
-B_{21}(t)&-B_{22}(t)\end{pmatrix},\ee
then $\dbA+\dbB(\cd)$ generates an evolution operator
$\h\F(\cd\,,\cd)$ on $X\times X$. The following result concerns the
well-posedness of the corresponding linear FBEE.

\ms

\bf Proposition 3.2. \sl Let $A^*=-A$ and {\rm(\ref{Bbg})} hold.
Then linear FBEE {\rm(\ref{linear})} admits a unique mild solution
$(y(\cd),\psi(\cd))$ for any $h_0\in X$ if and only if
\bel{-1}\left[(-H,I)\h\F(T,0)\begin{pmatrix}0\\
I\end{pmatrix}\right]^{-1}\in\cL(X\times X).\ee

\ms

\it Proof. \rm Suppose (\ref{linear}) admits a unique mild solution.
Then we have
\bel{}\begin{pmatrix}y(t)\\ \psi(t)\end{pmatrix}=\h\F(t,0)\begin{pmatrix}x\\
\psi(0)\end{pmatrix}+\int_0^t\h\F(t,s)\begin{pmatrix}b_0(s)\\-g_0(s)
\end{pmatrix}ds,\qq t\in[0,T],\ee
with $\psi(0)$ undetermined. By the condition at $t=T$, we have
$$\ba{ll}
\ns\ds h_0=-Hy(T)+\psi(T)=(-H,I)\begin{pmatrix}y(T)\\
\psi(T)\end{pmatrix}\\
\ns\ds\q=(-H,I)\h\F(T,0)\begin{pmatrix}x\\
\psi(0)\end{pmatrix}+\int_0^T(-H,I)\h\F(t,s)\begin{pmatrix}b_0(s)\\-g_0(s)
\end{pmatrix}ds\\
\ns\ds\q=(-H,I)\h\F(T,0)\begin{pmatrix}0\\
I\end{pmatrix}\psi(0)+(-H,I)\h\F(T,0)\begin{pmatrix}I\\
0\end{pmatrix}x+\int_0^T(-H,I)\h\F(t,s)\begin{pmatrix}b_0(s)\\-g_0(s)
\end{pmatrix}ds.\ea$$
Hence, in order for any $h_0\in X$ the above uniquely determines
$\psi(0)$, we need (\ref{-1}). Conversely, if (\ref{-1}) holds, one
obtains
\bel{psi(0)}\ba{ll}
\ns\ds\psi(0)=\left[(-H,I)\h\F(T,0)\begin{pmatrix}0\\
I\end{pmatrix}\right]^{-1}\[h_0-(-H,I)\h\F(T,0)\begin{pmatrix}I\\
0\end{pmatrix}x-\int_0^T(-H,I)\h\F(t,s)\begin{pmatrix}b_0(s)\\-g_0(s)
\end{pmatrix}ds\].\ea\ee
From this we obtain the mild solution $(y(\cd),\psi(\cd))$ of FBEE
(\ref{linear}). \endpf

\ms

We note that in principle, condition (\ref{-1}) is checkable,
although it might be practically complicated. We also note that, in
the above, the condition that $A^*=-A$, or $e^{At}$ is a group,
plays an essential role. It seems that if $e^{At}$ is not a group,
the arguments used above will not work (since $\h\F(\cd\,,\cd)$ in
the above might not be defined).

\ms

\rm

We now return to the general linear FBEE (\ref{linear}) (without
assuming (H0)). Suppose $(y(\cd),\psi(\cd))$ is a strong solution to
linear FBEE (\ref{linear}). Inspired by the well-known invariant
imbedding idea (\cite{Bellman-Wing 1975,Ma-Protter-Yong 1994,
Yong-Zhou 1999,Ma-Wu-Zhang-Zhang 2014}), we suppose that the
following relation holds:
$$\psi(t)=\dbP(t)y(t)+p(t),\qq t\in[0,T],$$
for some Fr\'echet differentiable functions $\dbP:[0,T]\to\cL(X)$
and $p:[0,T]\to X$. Then, formally, we should have
$$\ba{ll}
\ns\ds-A^*[\dbP(t)y(t)+p(t)]-B_{21}(t)y(t)-B_{22}(t)[\dbP(t)y(t)+p(t)]-g_0(t)
=\dot\psi(t)\\
\ns\ds=\dot\dbP(t)y(t)+\dbP(t)\(Ay(t)+B_{11}(t)y(t)
+B_{12}(t)[\dbP(t)y(t)+p(t)]+b_0(t)\)+\dot p(t)\\
\ns\ds=\(\dot\dbP(t)+\dbP(t)A+\dbP(t)B_{11}(t)+\dbP(t)B_{12}(t)\dbP(t)\)y(t)
+\dbP(t)B_{12}(t)p(t)+\dbP(t)b_0(t)+\dot
p(t).\ea$$
Hence,
$$\ba{ll}
\ns\ds0=\(\dot\dbP(t)+\dbP(t)A+A^*\dbP(t)+\dbP(t)B_{11}(t)+B_{22}(t)\dbP(t)
+\dbP(t)B_{12}(t)\dbP(t)+B_{21}(t)\)y(t)\\
\ns\ds\qq+\dot
p(t)+A^*p(t)+\(\dbP(t)B_{12}(t)+B_{22}(t)\)p(t)+\dbP(t)b_0(t)+g_0(t).\ea$$
This suggests that we choose $\dbP(\cd)$ satisfying the following:
\bel{Riccati}\left\{\2n\ba{ll}
\ns\ds\dot\dbP(t)+\dbP(t)A+A^*\dbP(t)+\dbP(t)B_{11}(t)+B_{22}(t)\dbP(t)
+\dbP(t)B_{12}(t)\dbP(t)+B_{21}(t)=0,\qq t\in[0,T],\\
\ns\ds\dbP(T)=H,\ea\right.\ee
and choose $p(\cd)$ satisfying
\bel{p}\left\{\2n\ba{ll}
\ns\ds\dot
p(t)+A^*p(t)+\(\dbP(t)B_{12}(t)+B_{22}(t)\)p(t)+\dbP(t)b_0(t)+g_0(t)=0,\qq t\in[0,T],\\
\ns\ds p(T)=h_0.\ea\right.\ee
Equation (\ref{Riccati}) is called a {\it differential Riccati
equation}. Any $\dbP:[0,T]\to\cL(X)$ is called a {mild solution} of
(\ref{Riccati}) if the following integral equation is satisfied:
\bel{Riccati2}\ba{ll}
\ns\ds\dbP(t)=e^{A^*(T-t)}He^{A(T-t)}+\int_t^Te^{A^*(s-t)}\[\dbP(s)B_{11}(s)
+B_{22}(s)\dbP(s)\\
\ns\ds\qq\qq\qq\qq\qq+\dbP(s)B_{12}(t)\dbP(s)+B_{21}(s)\]e^{A(s-t)}ds,\q
t\in[0,T].\ea\ee
Note that if $A$ is bounded, (\ref{Riccati}) and (\ref{Riccati2})
are equivalent. Further, recalling that $\F_{11}(\cd\,,\cd)$ and
$\F_{22}(\cd\,,\cd)$ are the evolution operators generated by
$A+B_{11}(\cd)$ and $A^*+B_{22}(\cd)$, respectively, one sees that
(\ref{Riccati2}) is also equivalent to the following:
\bel{Riccati3}\ba{ll}
\ns\ds\dbP(t)=\F_{22}(T,t)H\F_{11}(T,t)+\int_t^T\F_{22}(s,t)
\[\dbP(s)B_{12}(t)\dbP(s)+B_{21}(s)\]\F_{11}(s,t)ds,\q t\in[0,T].\ea\ee
Having the above derivation, we now present the following result.

\ms

\bf Proposition 3.3. \sl Let {\rm(H0)$'$} and {\rm(\ref{Bbg})} hold.
Let Riccati equation {\rm(\ref{Riccati})} admit a unique mild
solution $\dbP:[0,T]\to\cL(X)$. Then linear FBEE {\rm(\ref{linear})}
admits a unique mild solution $(y(\cd),\psi(\cd))$.

\ms

\it Proof. \rm For any $\l>0$, consider the
following:
$$\ba{ll}
\ns\ds\dbP_\l(t)=e^{A_\l^*(T-t)}He^{A_\l(T-t)}+\int_t^Te^{A_\l^*(s-t)}
\h\dbQ(s)e^{A_\l(s-t)}ds,\q t\in[0,T],\ea$$
where, with the mild solution $\dbP(\cd)$ of the Riccati equation
(\ref{Riccati}),
$$\h\dbQ(s)=\dbP(s)B_{11}(s)
+B_{22}(s)\dbP(s)+\dbP(s)B_{12}(t)\dbP(s)+B_{21}(s).$$
Clearly, $\dbP_\l(\cd)$ is uniformly bounded (by noting
(\ref{2.19})). Moreover, for any $x\in X$,
$$\ba{ll}
\ns\ds\|\dbP_\l(t)x-\dbP(t)x\|\le\|e^{A^*_\l(T-t)}He^{A_\l(T-t)}x-e^{A^*(T-t)}H
e^{A(T-t)}x\|\\
\ns\ds\qq+\int_t^T\|e^{A_\l^*(s-t)}\h\dbQ(s)e^{A_\l(s-t)}x-
e^{A^*(s-t)}\h\dbQ(s)e^{A(s-t)}x\|ds\\
\ns\ds\les\|e^{A^*_\l(T-t)}H\big[e^{A_\l(T-t)}x-e^{A(T-t)}x\big]\|
+\|e^{A^*_\l(T-t)}He^{A(T-t)}x-e^{A^*(T-t)}H
e^{A(T-t)}x\|\\
\ns\ds\qq+\int_t^T\(\|e^{A_\l^*(s-t)}\h\dbQ(s)\big[e^{A_\l(s-t)}x-
e^{A(s-t)}x\big]\|+\|e^{A_\l^*(s-t)}\h\dbQ(s)e^{A(s-t)}x-
e^{A^*(s-t)}\h\dbQ(s)e^{A(s-t)}x\|\)ds\\
\ns\ds\les
K\|e^{A_\l(T-t)}x-e^{A(T-t)}x\|+\|e^{A^*_\l(T-t)}He^{A(T-t)}x-e^{A^*(T-t)}H
e^{A(T-t)}x\|\\
\ns\ds\qq+\int_t^T\(K\|e^{A_\l(s-t)}x-
e^{A(s-t)}x\|+\|e^{A_\l^*(s-t)}\h\dbQ(s)e^{A(s-t)}x-
e^{A^*(s-t)}\h\dbQ(s)e^{A(s-t)}x\|\)ds\to0.\ea$$
Hereafter, $K>0$ represents a generic constant which can be
different from line to line. Note that $\dbP_\l(\cd)$ also solves
the following Lyapunov equation:
\bel{}\left\{\2n\ba{ll}
\ns\ds\dot\dbP_\l(t)+\dbP_\l(t)A_\l+A_\l^*\dbP_\l(t)+\h\dbQ(t)=0,\qq t\in[0,T],\\
\ns\ds\dbP_\l(T)=H.\ea\right.\ee
Now, we let $p_\l(\cd)$ be the solution of the following:
\bel{p_l}\left\{\2n\ba{ll}
\ns\ds\dot p_\l(t)\1n+\1n
A_\l^*p_\l(t)\1n+\1n\big[\dbP_\l(t)B_{12}(t)\1n+\1n
B_{22}(t)\big]p_\l(t)\1n+\1n\dbP_\l(t)b_0(t) \1n+\1n g_0(t)=0,\qq t\in[0,T],\\
\ns\ds p_\l(T)=h_0,\ea\right.\ee
It is clear that
$$\|p_\l(\cd)\|_\infty<\infty.$$
We estimate
$$\ba{ll}
\ns\ds\|p_\l(t)-p(t)\|\les\|e^{A_\l^*(T-t)}h_0-e^{A^*(T-t)}h_0\|\\
\ns\ds\qq\qq+\int_t^T\2n\Big\{\Big\|e^{A_\l^*(s-t)}\(\dbP_\l(s)B_{12}(s)+B_{22}(s)\)
p_\l(s)-e^{A^*(s-t)}\(\dbP(s)B_{12}(s)+B_{22}(s)\)p(s)\Big\|\\
\ns\ds\qq\qq+\Big\|e^{A_\l^*(s-t)}\(\dbP_\l(s)b_0(s)+g_0(s)\)-e^{A^*(s-t)}
\(\dbP(s)b_0(s)+g_0(s)\)\Big\|\Big\}ds\\
\ns\ds\les\|e^{A_\l^*(T-t)}h_0-e^{A^*(T-t)}h_0\|+\int_t^T\Big\{\Big\|e^{A_\l^*(s-t)}\(\dbP_\l(s)B_{12}(s)
+B_{22}(s)\)\(p_\l(s)-p(s)\)\Big\|\\
\ns\ds\qq\qq+\Big\|\(e^{A^*_\l(s-t)}-e^{A^*(s-t)}\)\(\dbP(s)B_{12}(s)+B_{22}(s)\)p(s)\Big\|\\
\ns\ds\qq\qq+\Big\|e^{A_\l^*(s-t)}\(\dbP_\l(s)-\dbP(s)\)\(B_{12}(s)+B_{22}(s)\)p(s)\Big\|\\
\ns\ds\qq\qq+\Big\|\(e^{A_\l^*(s-t)}-e^{A^*(s-t)}\)\(\dbP(s)b_0(s)+g_0(s)\)\Big\|
+\Big\|e^{A_\l^*(s-t)}
\(\dbP_\l(s)-\dbP(s)\)b_0(s)\Big\|\Big\}ds\\
\ns\ds\les\|e^{A_\l^*(T-t)}h_0-e^{A^*(T-t)}h_0\|\\
\ns\ds\qq\q\int_t^T\Big\{K\|p_\l(s)-p(s)\|+\Big\|\big(e^{A^*_\l(s-t)}-e^{A^*(s-t)}\big)
\big(\dbP(s)B_{12}(s)+B_{22}(s)\big)p(s)\Big\|\\
\ns\ds\qq\q+\Big\|\big(e^{A_\l^*(s-t)}-e^{A^*(s-t)}\big)\big(\dbP(s)b_0(s)+g_0(s)\big)
\Big\|\\
\ns\ds\qq\q+K\(\Big\|\big(\dbP_\l(s)-\dbP(s)\big)\big(B_{12}(s)+B_{22}(s)\big)p(s)\Big\|
\1n+\1n\Big\|\big(\dbP_\l(s)-\dbP(s)\big)b_0(s)\Big\|\)\Big\}ds.\ea$$
Then by Gronwall's inequality, we have
$$\ba{ll}
\ns\ds\|p_\l(t)-p(t)\|\les K\|e^{A_\l^*(T-t)}h_0-e^{A^*(T-t)}h_0\|\\
\ns\ds\qq\qq+K\int_t^T\Big\{\Big\|\(e^{A^*_\l(s-t)}-e^{A^*(s-t)}\)\(\dbP(s)B_{12}(s)+B_{22}(s)\)p(s)\Big\|\\
\ns\ds\qq\qq+\Big\|\(e^{A_\l^*(s-t)}-e^{A^*(s-t)}\)\(\dbP(s)b_0(s)+g_0(s)\)\Big\|\\
\ns\ds\qq\qq+\Big\|\big(\dbP_\l(s)-\dbP(s)\big)\big(B_{12}(s)+B_{22}(s)\big)p(s)\Big\|
+\Big\|\big(\dbP_\l(s)-\dbP(s)\big)b_0(s)\Big\|\Big\}ds\to0.\ea$$
Now, let $y_\l(\cd)$ be the solution to the following:
$$\left\{\2n\ba{ll}
\ns\ds\dot y_\l(t)=A_\l
y_\l(t)+B_{11}(t)y_\l(t)+B_{12}(t)\big[\dbP_\l(t)y_\l(t)+p_\l(t)\big]+b_0(t),\q t\in[0,T],\\
\ns\ds y_\l(0)=x.\ea\right.$$
By the convergence of $\dbP_\l(\cd)\to\dbP(\cd)$ and $p_\l(\cd)\to
p(\cd)$, we have
$$\lim_{\l\to\infty}\|y_\l(\cd)-y(\cd)\|_\infty=0.$$
Define
$$\psi_\l(t)=\dbP_\l(t)y_\l(t)+p_\l(t).$$
Then one has
$$\left\{\2n\ba{ll}
\ns\ds\dot y_\l(t)=A_\l
y_\l(t)+B_{11}(t)y_\l(t)+B_{12}(t)\psi_\l(t)+b_0(t),\\
\ns\ds y_\l(0)=x,\ea\right.$$
and
$$\ba{ll}
\ns\ds\dot\psi_\l(t)=\dot\dbP_\l(t)y_\l(t)+\dbP_\l(t)\dot
y_\l(t)+\dot p_\l(t)\\
\ns\ds=-\(\dbP_\l(t)A_\l+A_\l^*\dbP_\l(t)+\dbP(t)B_{11}(t)
+B_{22}(t)\dbP(t)+\dbP(t)B_{12}(t)\dbP(s)+B_{21}(t)\)y_\l(t)\\
\ns\ds\qq+\dbP_\l(t)\Big\{A_\l
y_\l(t)+B_{11}(t)y_\l(t)+B_{12}(t)\[\dbP_\l(t)y_\l(t)+p_\l(t)\]
+b_0(t)\Big\}\\
\ns\ds\qq-\[A_\l^*p_\l(t)+\(\dbP_\l(t)B_{12}(t)+B_{22}(t)\)p_\l(t)+\dbP_\l(t)b_0(t)
+g_0(t)\]\\
\ns\ds=-A_\l^*\big[\dbP_\l(t)y_\l(t)+p_\l(t)\big]-B_{21}(t)y_\l(t)
-B_{22}(t)\big[\dbP(t)y_\l(t)+p_\l(t)\big]-g_0(t)\\
\ns\ds\qq+\big[\dbP_\l(t)-\dbP(t)\big]B_{11}(t)y_\l(t)+\big[\dbP_\l(t)B_{12}(t)\dbP_\l(t)
-\dbP(t)B_{12}(t)\dbP(t)\big]y_\l(t)\\
\ns\ds=-A_\l^*\psi_\l(t)-B_{21}(t)y_\l(t)
-B_{22}(t)\psi_\l(t)-g_0(t)+R_\l(t),\ea$$
where
$$R_\l(t)=\(\big[\dbP_\l(t)-\dbP(t)\big]B_{11}(t)+\big[\dbP_\l(t)B_{12}(t)\dbP_\l(t)
-\dbP(t)B_{12}(t)\dbP(t)\big]\)y_\l(t).$$
Since
$$\ba{ll}
\ns\ds\|R_\l(t)\|\les
K\|y_\l(t)-y(t)\|+\Big\|\(\big[\dbP_\l(t)-\dbP(t)\big]B_{11}(t)+\big[\dbP_\l(t)B_{12}(t)\dbP_\l(t)
-\dbP(t)B_{12}(t)\dbP(t)\big]\)y(t)\Big\|\\
\ns\ds\les K\|y_\l(t)-y(t)\|+\|\big[\dbP_\l(t)-\dbP(t)\big]B_{11}(t)y(t)\|\\
\ns\ds\q+\|\big[\dbP_\l(t)-\dbP(t)\big]B_{12}(t)\dbP(t)y(t)\|+\|\dbP_\l(t)B_{12}(t)
\big[\dbP_\l(t)-\dbP(t)\big]y(t)\|\\
\ns\ds\les K\|y_\l(t)-y(t)\|+\|\big[\dbP_\l(t)-\dbP(t)\big]B_{11}(t)y(t)\|\\
\ns\ds\q+\|\big[\dbP_\l(t)-\dbP(t)\big]B_{12}(t)\dbP(t)y(t)\|
+K\|\big[\dbP_\l(t)-\dbP(t)\big]y(t)\|\to0,\ea$$
we see that
$$\psi_\l(\cd)\to\psi(\cd),$$
and $(y(\cd),\psi(\cd))$ is a mild solution to FBEE (\ref{linear}).
\endpf

\ms

Now, a natural question is when the Riccati equation (\ref{Riccati})
admits a mild solution. Let us rewrite the Riccait equation as
follows:
\bel{Riccati4}\ba{ll}
\ns\ds\dbP(t)=\F_{22}(T,t)H\F_{11}(T,t)+\int_t^T\F_{22}(s,t)B_{21}(s)\F_{11}(s,t)ds\\
\ns\ds\qq\qq+\int_t^T\F_{22}(s,t)
\dbP(s)B_{12}(s)\dbP(s)\F_{11}(s,t)ds,\qq t\in[0,T].\ea\ee
A trivial case is $B_{12}(\cd)=0$ for which the above equation is
linear and it always has a solution, under (H0)$'$ and (\ref{Bbg}).
In general, we have the following result.

\ms

\bf Proposition 3.4. \sl Suppose {\rm(H0)$'$} and $(\ref{Bbg})$
hold.

\ms

{\rm(i)} Equation {\rm(\ref{Riccati4})} admits at most one solution
$\dbP(\cd)\in C([0,T];\cL(X))$.

\ms

{\rm(ii)} Suppose in addition that
\bel{*=}B_{22}(t)=B_{11}(t)^*,\q t\in[0,T],\ee
and
\bel{ge0}H\in\dbS^+(X),\qq -B_{12}(\cd),B_{21}(\cd)\in
L^\infty(0,T;\dbS^+(X)).\ee
%
%
Then Riccati equation {\rm(\ref{Riccati4})} admits a unique solution
$\dbP(\cd)\in C([0,T];\dbS^+(X))$.

\rm

\ms

\it Proof. \rm (i) Suppose $\dbP(\cd),\wt\dbP(\cd)\in
C([0,T];\cL(X))$ are two solutions to (\ref{Riccati4}). Then
$\h\dbP(\cd)\equiv\dbP(\cd)-\wt\dbP(\cd)$ satisfies the following:
$$\ba{ll}
\ns\ds\h\dbP(t)=\int_t^T\F_{22}(s,t)\[\dbP(s)B_{12}(s)\dbP(s)-\wt\dbP(s)B_{12}(s)
\wt\dbP(s)\]\F_{11}(s,t)ds\\
\ns\ds\qq=\int_t^T\F_{22}(s,t)\[\h\dbP(s)B_{12}(s)\dbP(s)+\wt\dbP(s)B_{12}(s)
\h\dbP(s)\]\F_{11}(s,t)ds.\ea$$
Here, we note that $\dbP(\cd)=\h\dbP(\cd)+\wt\dbP(\cd)$. Hence
$$\ba{ll}
\ns\ds\dbP(s)B_{12}(s)\dbP(s)-\wt\dbP(s)B_{12}(s)\wt\dbP(s)
=\[\h\dbP(s)+\wt\dbP(s)\]B_{12}(s)\dbP(s)-\wt\dbP(s)B_{12}(s)\[\dbP(s)-\h\dbP(s)\]\\
\ns\ds=\h\dbP(s)B_{12}(s)\dbP(s)+\wt\dbP(s)B_{12}(s)\h\dbP(s).\ea$$
Then by Gronwall's inequality, we obtain that
$$\wt\dbP(\cd)=\dbP(\cd).$$

(ii) Under our conditions, we have
$$A^*+B_{22}(\cd)=\big[A+B_{11}(\cd)\big]^*.$$
Hence,
$$\F_{22}(\cd\,,\cd)=\F_{11}(\cd\,,\cd)^*,$$
and (\ref{Riccati}) can be written as
\bel{Riccati*}\left\{\2n\ba{ll}
\ns\ds\dot\dbP(t)+\dbP(t)\big[A+B_{11}(t)+B_{12}(t)\dbP(t)\big]
+\big[A+B_{11}(t)+B_{12}(t)\dbP(t)\big]^*\dbP(t)\\
\ns\ds\qq+\dbP(t)\big[-B_{12}(t)\big]\dbP(t)+B_{21}(t)=0,\qq t\in[0,T],\\
\ns\ds\dbP(T)=H,\ea\right.\ee
We now let
$$\dbP_0(t)=0,\qq t\in[0,T],$$
and let $\dbP_{n+1}(\cd)$ be the mild solution of the following
Lyapunov equation:
\bel{Riccati*}\left\{\2n\ba{ll}
\ns\ds\dot\dbP_{n+1}(t)+\dbP_{n+1}(t)\big[A+B_{11}(t)+B_{12}(t)\dbP_n(t)\big]
+\big[A+B_{11}(t)+B_{12}(t)\dbP_n(t)\big]^*\dbP_{n+1}(t)\\
\ns\ds\qq+\dbP_n(t)\big[-B_{12}(t)\big]\dbP_n(t)+B_{21}(t)=0,\qq t\in[0,T],\\
\ns\ds\dbP_{n+1}(T)=H,\ea\right.\ee
Observe the following:
$$\ba{ll}
\ns\ds0=\dot\dbP_{n+1}(t)-\dot\dbP_n(t)+\big[\dbP_{n+1}(t)-\dbP_n(t)\big]
\big[A+B_{11}(t)\big]+\big[A+B_{11}(t)\big]^*\big[\dbP_{n+1}(t)-\dbP_n(t)\big]\\
\ns\ds\qq+\dbP_{n+1}(t)B_{12}(t)\dbP_n(t)+\dbP_n(t)B_{12}(t)\dbP_{n+1}(t)
-\dbP_n(t)B_{12}(t)\dbP_{n-1}(t)-\dbP_{n-1}(t)B_{12}(t)\dbP_n(t)\\
\ns\ds\qq-\dbP_n(t)B_{12}(t)\dbP_n(t)+\dbP_{n-1}(t)B_{12}(t)\dbP_{n-1}(t)\\
\ns\ds\q=\dot\dbP_{n+1}(t)-\dot\dbP_n(t)+\big[\dbP_{n+1}(t)-\dbP_n(t)\big]
\big[A+B_{11}(t)+B_{12}(t)\dbP_n(t)\big]\\
\ns\ds\qq+\big[A+B_{11}(t)+B_{12}(t)\dbP_n(t)\big]^*\big[\dbP_{n+1}(t)-\dbP_n(t)\big]
+\big[\dbP_n(t)-\dbP_{n-1}(t)\big]B_{12}(t)\big[\dbP_n(t)-\dbP_{n-1}(t)\big].\ea$$
This implies that
$$\dbP_{n+1}(t)\les\dbP_n(t),\qq t\in[0,T],\qq n\ges1.$$
On the other hand, from (\ref{Riccati*}), one has
$$\dbP_n(t)\ges0,\qq\forall t\in[0,T],~n\ges1.$$
Hence, by \cite{Riesz-Nagy 1955, Muscat 2014} for any $t\in[t,T]$,
there exists a $\dbP(t)\in\dbS^+(X)$ such that
\bel{limit}\lim_{n\to\infty}\|\dbP_n(t)x-\dbP(t)x\|=0,\qq\forall
x\in X.\ee
Note that for any $x\in X$,
$$\ba{ll}
\ns\ds\dbP_{n+1}(t)x=\F_{11}(T,t)^*H\F_{11}(T,t)x
+\int_t^T\F_{11}(s,t)^*B_{21}(s)\F_{11}(s,t)xds\\
\ns\ds\qq\qq\qq+\int_t^T\F_{11}(s,t)^*\[\dbP_{n+1}(s)B_{12}(s)\dbP_n(s)
+\dbP_n(s)B_{12}(s)\dbP_{n+1}(s)\\
\ns\ds\qq\qq\qq\qq-\dbP_n(s)B_{12}(s)\dbP_n(s)\]\F_{11}(s,t)xds,\qq
t\in[0,T].\ea$$
Thus, making use of (\ref{limit}), we obtain that $\dbP(\cd)$ is a
mild solution to (\ref{Riccati4}). \endpf

\ms

Let us look at linear FBEE (\ref{FBEE-LQ}) resulting from
linear-quadratic optimal control problem. We rewrite (\ref{FBEE-LQ})
below:
\bel{FBEE-LQ*}\left\{\2n\ba{ll}
\ns\ds\dot
y(t)=Ay(t)-B(t)R(t)^{-1}S(t)y(t)-B(t)R(t)^{-1}B(t)^*\psi(t),\\
[1mm]
\ns\ds\dot\psi(t)=-A^*\psi(t)-\big[Q(t)-S(t)^*R(t)^{-1}S(t)\big]y(t)
+S(t)^*R(t)^{-1}B(t)^*\psi(t),\\
[1mm]
\ns\ds y(0)=x,\qq\psi(T)=Gy(T).\ea\right.\ee
Note that
$$R(t)\ges\d I,\qq\forall t\in[0,T],$$
which leads to the existence of $R(t)^{-1}$. Hence,
$$\left\{\2n\ba{ll}
\ns\ds B_{11}(t)=-B(t)R(t)^{-1}S(t),\qq
B_{12}(t)=-B(t)R(t)^{-1}B(t)^*,\\
\ns\ds B_{21}(t)=\big[Q(t)-S(t)^*R(t)^{-1}S(t)\big],\qq
B_{22}(t)=-S(t)^*R(t)^{-1}B(t)^*.\ea\right.$$
Then in the case that
$$G\ges0,\qq Q(t)-S(t)^*R(t)^{-1}S(t)\ges0,\qq t\in[0,T],$$
the corresponding Riccati equation admits a unique solution and
linear FBEE (\ref{FBEE-LQ}) admits a solution.

\section{Decoupling Method --- A Brief Description}

We note that in the previous section, the essence of the approach by
means of Riccati equation is to use the ansatz
$$\psi(t)=\dbP(t)y(t)+p(t),\qq t\in[0,T].$$
Inspired by this, we now look at nonlinear cases. For FBEE
(\ref{FBEE1}), suppose
$$\psi(t)=\dbK(t,y(t)),\qq t\in[0,T],$$
for some $\dbK(\cd\,,\cd)$. Let $\{\z_n\}_{n\ge1}$ be an orthonormal
basis of $X$. Then
$$\dbK(t,y)=\sum_{n=1}^\infty\lan\dbK(t,y),\z_n\ran\z_n\equiv\sum_{n=1}^\infty
k^n(t,y)\z_n,$$
with $k^n:[0,T]\times X\to\dbR$. Suppose $(t,y)\mapsto\dbK(t,y)$ is
Fr\'echet differentiable. Then, so is $(t,y)\mapsto
k^n(t,y)=\lan\dbK(t,y),\z_n\ran$, and for any $z\in X$,
$$\ba{ll}
\ns\ds{\dbK(t,y+\d z)-\dbK(t,y)\over\d}=\sum_{n=1}^\infty{k^n(t,y+\d
z)-k^n(t,y)\over\d}\z_n\\
\ns\ds\qq\qq\qq\qq\qq\to\sum_{n=1}^\infty\lan
k^n_y(t,y),z\ran\z_n\equiv\[\sum_{n=1}^\infty\z_n\otimes
k_y^n(t,y)\]z.\ea$$
This means
\bel{K_y}\dbK_y(t,y)=\sum_{n=1}^\infty\z_n\otimes
k_y^n(t,y),\qq(t,y)\in[0,t]\times X.\ee
Note that $\dbK_y(t,y)$ is independent of the choice of
$\{\z_n\}_{n\ge1}$.

\ms

We now present the following result.

\ms

\bf Proposition 4.1. \sl Let {\rm(H0)$'$} and {\rm(H2)} hold. Let
$\dbK:[0,T]\times X\to X$ be Fr\'echet differentiable satisfying the
following:
\bel{K}\left\{\2n\ba{ll}
\ns\ds\dbK_t(t,y)\1n+\1n\dbK_y(t,y)\big[Ay\1n+\1n
b(t,y,\dbK(t,y))\big] \1n+\1n A^*\dbK(t,y)\1n+\1n
g(t,y,\dbK(t,y))=0,\q(t,y)\in[0,T]\times\cD(A),\\
\ns\ds\dbK(T,y)=h(y),\qq y\in X.\ea\right.\ee
Let $y(\cd)$ be a classical solution to the following:
\bel{closed}\left\{\2n\ba{ll}
\ns\ds\dot y(t)=Ay(t)+b(t,y(t),\dbK(t,y(t)),\qq t\in[0,T],\\
\ns\ds y(0)=x,\ea\right.\ee
and
\bel{psi=K}\psi(t)=\dbK(t,y(t)),\qq t\in[0,T].\ee
Then $(y(\cd),\psi(\cd))$ is a strong solution of FBEE
{\rm(\ref{FBEE1})}.

\ms

\it Proof. \rm Note that as a part of requirement for
$\dbK(\cd\,,\cd)$ being a solution to (\ref{K}), one has
$$\dbK(t,y)\in\cD(A^*),\qq\forall(t,y)\in[0,T]\times X.$$
By (\ref{psi=K}), we have
$$\ba{ll}
\ns\ds\dot\psi(t)=\dbK_t(t,y(t))+\dbK_y(t,y(t))\[Ay(t)+b(t,y(t),\dbK(t,y(t))\]\\
\ns\ds\qq=-A^*\dbK(t,y(t))-g(t,y(t),\dbK(t,y(t))=-A^*\psi(t)-g(t,y(t),\psi(t)),\ea$$
and
$$\psi(T)=\dbK(T,y(T))=h(y(T)).$$
Hence, our claim follows. \endpf

\ms

It is seen that thanks to the map $\dbK(\cd\,,\cd)$, the original
FBEE (\ref{FBEE1}) is decoupled into (\ref{closed}) and
(\ref{psi=K}). Because of this, we introduce the following notion:

\ms

\bf Definition 4.2. \rm A map $\dbK:[0,T]\times X\to X$ is called a
{\it decoupling field} of FBEE (\ref{FBEE1}) if it is a solution to
(\ref{K}).

\ms

Now the natural question is when one can solve equation (\ref{K}).
The linear case has been treated in the previous section. To look at
the nonlinear case, let us further assume that $A^*=A$ and it has a
sequence of eigenvalues
$$0>-\si_0=\l_1\ges\l_2\ges\l_2\cds,$$
with the corresponding eigenfunctions $\{\z_n\}_{n\ges1}$ which form
an orthonormal basis for $X$. Then
$$\dbK(t,y)=\sum_{n=1}^\infty\lan\dbK(t,y),\z_n\ran\z_n\equiv\sum_{n=1}^\infty
k^n(t,y)\z_n.$$
Hence,
$$\ba{ll}
\ns\ds0=\dbK_t(t,y)\1n+\1n\dbK_y(t,y)\big[Ay\1n+\1n
b(t,y,\dbK(t,y))\big] \1n+\1n A^*\dbK(t,y)\1n+\1n
g(t,y,\dbK(t,y))\\
\ns\ds\q=\sum_{n=1}^\infty\(k^n_t(t,y)+\lan
k^n_y(t,y),Ay+b(t,y,\dbK(t,y))\ran+\l_nk^n(t,y)+\lan
g(t,y,\dbK(t,y)),\z_n\ran\)\z_n.\ea$$
Therefore, we obtain a coupled system of countably many equations:
\bel{kn}\left\{\2n\ba{ll}
\ns\ds k^n_t(t,y)+\lan
k^n_y(t,y),Ay+b(t,y,\dbK(t,y))\ran+\l_nk^n(t,y)\\
\ns\ds\qq\qq+\lan g(t,y,\dbK(t,y)),\z_n\ran=0,\qq(t,y)\in[0,T]\times
X,\\
\ns\ds k^n(T,y)=\lan h(y),\z_n\ran,\qq y\in X.\ea\right.\ee
Let us look at a special case. Suppose
$$\dbK(t,y)=k^1(t,y)\z_1,\qq(t,y)\in[0,T]\times X,~\z_1\in\cD(A^*).$$
Then
\bel{k}\left\{\2n\ba{ll}
\ns\ds k^1_t(t,y)+\lan
k^1_y(t,y),Ay+b(t,y,k^1(t,y)\z_1)\ran+\l_1k^1(t,y)\\
\ns\ds\qq\qq+\lan g(t,y,k^1(t,y)\z_1),\z_1\ran
=0,\qq(t,y)\in[0,T]\times X,\\
\ns\ds k^1(T,y)=h^1(y),\qq y\in X,\ea\right.\ee
and
\bel{k=0}\lan g(t,y,k^1(t,y)\z_1),\z_n\ran=0,\qq n>1.\ee
To guarantee (\ref{k=0}), we assume that
\bel{g4.7}g(t,y,\span\{\z_1\})\subseteq\span\{\z_1\}.\ee

\ms

We note that (\ref{k}) is a first order Hamilton-Jacobi equation in
the Hilbert space $X$, involving an unbounded operator $A$.
Therefore, it is possible to study the existence of viscosity
solution of it. When the viscosity solution has certain regularity,
one might be able to obtain a decoupling field
$\dbK(t,y)=k^1(t,y)\z_1$ for our FBEE. Apparently, this is merely a
very special case for the general FBEEs, and it already looks
complicated. Hence, there is a very long way to go in this direction
to establish a satisfactory theory (for nonlinear FBEEs). We hope to
report some further results in this direction in our future
publications.

\section{Lyapunov Operators and a Priori Estimates}

We now look at the solvability by another method, called {\it method
of continuity}. We first look at the following linear FBEE:
\bel{FBEE2}\left\{\2n\ba{ll}
\ns\ds\dot y(t)=Ay(t)+B_{11}(t)y(t)+B_{12}(t)\psi(t)+b_0(t),\qq t\in[0,T],\\
\ns\ds\dot\psi(t)=-A^*\psi(t)-B_{21}(t)y(t)-B_{22}(t)\psi(t)-g_0(t),\qq t
\in[0,T],\\
\ns\ds\hb{$y(0)$ and $\psi(T)$ are given},\ea\right.\ee
with $B_{ij}:[0,T]\to\cL(X)$. Let
\bel{AB}\dbA=\begin{pmatrix}\sc A&\sc0\\
\sc0&\sc-A^*\end{pmatrix},\q\dbB(t)=\begin{pmatrix}\sc B_{11}(t)&\sc B_{12}(t)\\
\sc-B_{21}(t)&\sc-B_{22}(t)\end{pmatrix}.\ee
Then FBEE (\ref{FBEE2}) can be written as
\bel{}\left\{\2n\ba{ll}
\ns\ds\begin{pmatrix}\sc\dot y(t)\\
\sc\dot\psi(t)\end{pmatrix}=\big[\dbA+\dbB(t)\big]\begin{pmatrix}\sc y(t)\\
\sc\psi(t)\end{pmatrix}+\begin{pmatrix}\sc b_0(t)\\
\sc-g_0(t)\end{pmatrix},\qq t\in[0,T],\\ [4mm]
\ns\ds\hb{$y(0)$ and $\psi(T)$ are given}.\ea\right.\ee
We introduce the following {\it Lyapunov differential equation} for
operator-valued function $\Pi(\cd)$:
\bel{Pi}\dot\Pi(t)+\Pi(t)[\dbA-\dbM(t)]+[\dbA-\dbM(t)]^*\Pi(t)+\dbQ(t)=0,\q
t\in[0,T],\ee
where
$$\Pi(t)=\begin{pmatrix}\sc P(t)&\sc\G(t)^*\\
\sc\G(t)&\sc\bar P(t)\end{pmatrix},\q\dbM(t)=\begin{pmatrix}\sc
M(t)&\sc0\\
\sc 0&\sc-\bar M(t)^*\end{pmatrix},\q
\dbQ(t)=\begin{pmatrix}\sc Q_0(t)&\sc\Th(t)^*\\
\sc\Th(t)&\sc\bar Q_0(t)\end{pmatrix},$$
with $M,\bar M,\Th:[0,T]\to\cL(X)$ and $Q_0,\bar
Q_0:[0,T]\to\dbS(X)$ to be properly chosen later. We may
equivalently write (\ref{Pi}) as follows:
\bel{L1}\dot P(t)+P(t)[A-M(t)]+[A^*-M(t)^*]P(t)+Q_0(t)=0,\ee
\bel{L2}\dot{\bar P}(t)-\bar P(t)[A^*-\bar M(t)^*]-[A-\bar M(t)]\bar
P(t)+\bar Q_0(t)=0,\ee
\bel{Psi}\dot\G(t)+\G(t)[A-M(t)]-[A-\bar M(t)]\G(t)+\Th(t)=0,\ee
and
\bel{Psi*}\dot\G(t)^*-\G(t)^*[A^*-\bar
M(t)^*]+[A^*-M(t)^*]\G(t)^*+\Th(t)^*=0.\ee
Let us first look at (\ref{L1}) and (\ref{L2}). Operator-valued
functions $P(\cd)$ and $\bar P(\cd)$ are mild solutions to
(\ref{L1}) and (\ref{L2}), respectively, if the following hold:
\bel{L1a*}\ba{ll}
\ns\ds P(t)\1n=\1n e^{A^*(T-t)}P(T)e^{A(T-t)}-\2n\int_t^T\2n
e^{A^*(s-t)}[P(s)M(s)\1n+\2n M(s)^*P(s)\1n-\1n
Q_0(s)]e^{A(s-t)}ds,~t\1n\in\1n[0,T],\ea\ee
and
\bel{L2a*}\ba{ll}
\ns\ds\bar P(t)\1n=\1n e^{At}\bar P(0)e^{A^*t}-\2n\int_0^t\2n
e^{A(t-s)}[\bar P(s)\bar M(s)^*\2n+\1n\bar M(s)\bar P(s)\1n+\1n\bar
Q_0(s)]e^{A^*(t-s)}ds,~t\1n\in\1n[0,T].\ea\ee
We use the above definition simply because when $A$ is bounded,
(\ref{L1}) is equivalent to (\ref{L1a*}), and (\ref{L2}) is
equivalent to (\ref{L2a*}). Further, if we let $\F(\cd\,,\cd)$ and
$\bar\F(\cd\,,\cd)$ be the evolution operators generated by
$A-M(\cd)$ and $A-\bar M(\cd)$, respectively, then $P(\cd)$ and
$\bar P(\cd)$ admit the following representation:
\bel{L1a**}P(t)=\F(T,t)^*P(T)\F(T,t)\1n+\2n\int_t^T\2n
\F(s,t)^*Q_0(s)\F(s,t)ds,\qq t\in[0,T],\ee
and
\bel{L2a**}\bar P(t)=\bar\F(t,0)\bar P(0)\bar\F(t,0)^*-\int_0^t\2n
\bar\F(t,s)\bar Q_0(s)\bar\F(t,s)^*ds,\qq t\in[0,T].\ee
This yields that if
$$P(T),-\bar P(0),Q_0(t),\bar Q_0(t)\ges0,\qq t\in[0,T],$$
then
\bel{P>0}P(t)\ges0,\q\bar P(t)\les0,\qq t\in[0,T].\ee

\ms

Now, let us look at (\ref{Psi}) and (\ref{Psi*}), which are
equivalent. We assume that (H0) holds. Therefore, we have two cases
to discuss.

\ms

\it Case 1. \rm Let (\ref{case1}) hold. In this case, since $A$ is
dissipative, the appearance of the term $\G(t)A-A\G(t)$ makes
(\ref{Psi}) and (\ref{Psi*}) difficult to solve in general. To
overcome this, we require that for all $t\in[0,T]$,
$\G(t):\cD(A)\to\cD(A)$ and
\bel{PsiA=APsi}\G(t)Ax=A\G(t)x,\qq t\in[0,T],~x\in\cD(A).\ee
Then both (\ref{Psi}) and (\ref{Psi*}) are reduced to the following:
\bel{G1}\dot\G(t)-\G(t)M(t)+\bar M(t)\G(t)+\Th(t)=0,\ee
which admits a unique solution as long as, say, $M(\cd)$, $\bar
M(\cd)$ and $\Th(\cd)$ are bounded and $\G(T)\in\cL(X)$ is given.
Actually, if $\Psi(\cd\,,\cd)$ and $\bar\Psi(\cd\,,\cd)$ are
evolution operators generated by $-M(\cd)$ and $\bar M(\cd)$,
respectively, then
\bel{G}\G(t)=\bar\Psi(T,t)\G(T)\Psi(T,t)+\int_t^T\bar\Psi(s,t)\Th(s)\Psi(s,t)
ds,\q t\in[0,T].\ee
Note that under (\ref{case1}), $A$ admits a spectral decomposition
$$A=\int_{\si(A)}\m dE_\m,$$
with $\si(A)\subseteq(-\infty,-\si_0]$ being the spectral of $A$, if
\bel{3.24}\left\{\2n\ba{ll}
\ns\ds\G(T)=\int_{\si(A)}\g(\m)dE_\m,\q\Th(t)=\int_{\si(A)}\th(t,\m)dE_\m,\\
\ns\ds M(t)=\int_{\si(A)}m(t,\m)dE_\m,\q\bar M(t)=\int_{\si(A)}\bar
m(t,\m)dE_\m,\ea\right.\ee
for some suitable maps $\g:\si(A)\to\dbR$ and $\th,m,\bar
m:[0,T]\times\si(A)\to\dbR$, then
$$\G(t)=\int_{\si(A)}\(e^{\int_t^T[\bar m(s,\m)-m(s,\m)]ds}\g(\m)
+\int_t^Te^{\int_t^\t[\bar m(s,\m)-m(s,\m)]ds}
\th(\t,\m)d\t\)dE_\m.$$
Hence, (\ref{PsiA=APsi}) holds in this case. In particular, if
\bel{3.25}\left\{\2n\ba{ll}
\ns\ds\g(\m)=\g,\q\th(t,\m)=\th(t),\\
\ns\ds m(t,\m)=m(t),\q\bar m(t,\m)=\bar
m(t),\ea\right.\q(t,\m)\in[0,T]\times\si(A),\ee
we have
$$\G(t)=\(e^{\int_t^T[\bar m(s)-m(s)]ds}\g+\int_t^Te^{\int_t^\t[\bar m(s)-m(s)]ds}
\th(\t)d\t\)I.$$
Also, as a special case of (\ref{3.24}), if
\bel{3.26}\left\{\2n\ba{ll}
\ns\ds\G(t)=\g(t)A_\l,\qq\Th(t)=\th(t)A_\l,\\
\ns\ds M(t)=m(t)A_\l,\qq\bar M(t)=\bar m(t)A_\l,\ea\right.\qq
t\in[0,T],\ee
for some suitable scalar functions $\g(\cd),\th(\cd),m(\cd),\bar
m(\cd)$, then
$$\G(t)=\(e^{\int_t^T[\bar
m(s)-m(s)]ds}\g+\int_t^Te^{\int_t^\t[\bar m(s)-m(s)]ds}
\th(\t)d\t\)A_\l,$$
for which (\ref{PsiA=APsi}) will also hold.

\ms

\it Case 2. \rm Let (\ref{case2}) hold. Then $e^{At}$ is a group.
Consequently, $e^{-At}$ is well-defined. Similar to the case of
$P(\cd)$, a map $\G(\cd)$ is called a mild solution to equation
(\ref{Psi}), if the following holds:
\bel{3.29}\ba{ll}
\ns\ds\G(t)=e^{-A(T-t)}\G(T)e^{A(T-t)}\\
\ns\ds\qq\q+\2n\int_t^T\2n e^{-A(s-t)}[\bar
M(s)\G(s)\1n-\1n\G(s)M(s)\1n+\1n\Th(s)]e^{A(s-t)}ds,\q
t\in[0,T].\ea\ee
By recalling the evolution operators $\F(\cd\,,\cd)$ and
$\bar\F(\cd\,,\cd)$ generated by $A-M(\cd)$ and $A-\bar M(\cd)$,
(noting that in the current case, $\bar\F(s,t)^{-1}$ exists) we have
\bel{3.26}\G(t)=\bar\F(T,t)^{-1}\G(T)\F(T,t)+\int_t^T\bar\F(s,t)^{-1}
\Th(s)\F(s,t)ds,\q t\in[0,T].\ee
We point out that in the current case, (\ref{PsiA=APsi}) is not
needed. However, if (\ref{3.24}) holds and we are working in a
complex Hilbert space, we will still have (\ref{PsiA=APsi}) and
$\G(\cd)$ can also be given by (\ref{G}).

\ms

In what follows, when we say a mild solution $\Pi(\cd)$ of
(\ref{Pi}), we mean that $P(\cd)$ and $\bar P(\cd)$ are given by
(\ref{L1a*}) and (\ref{L2a*}), respectively, and $\G(\cd)$ is
defined by by (\ref{G}) such that (\ref{PsiA=APsi}) holds for the
case $A=A^*$ and $\G(\cd)$ is defined by (\ref{3.26}) for the case
$A^*=-A$ (since we prefer to stay with a real Hilbert space).

\ms

The following is the main result of this section and it will play an
important role below.

\ms

\bf Proposition 5.1. \sl Let $(y(\cd),\psi(\cd))$ be a mild solution
of linear FBEE $(\ref{FBEE2})$ and $\Pi(\cd)$ be a mild solution of
Lyapunov differential equation $(\ref{Pi})$. Then
\bel{}\ba{ll}
\ns\ds\lan\Pi(T)\begin{pmatrix}y(T)\\
\psi(T)\end{pmatrix},\begin{pmatrix}y(T)\\
\psi(T)\end{pmatrix}\ran-\lan\Pi(0)\begin{pmatrix}y(0)\\
\psi(0)\end{pmatrix},\begin{pmatrix}y(0)\\
\psi(0)\end{pmatrix}\ran\\
\ns\ds=\2n\int_0^T\2n\[\1n\lan\1n\(\Pi(t)[\dbB(t)\1n+\1n\dbM(t)]\1n
+\1n[\dbB(t)\1n+\1n\dbM(t)]^*\Pi(t)
\1n-\1n\dbQ(t)\)\1n\begin{pmatrix}y(t)\\
\psi(t)\end{pmatrix},\1n\begin{pmatrix}y(t)\\
\psi(t)\end{pmatrix}\ran\\
\ns\ds\qq\qq+2\lan\Pi(t)\begin{pmatrix}b_0(t)\\
-g_0(t)\end{pmatrix},\begin{pmatrix}y(t)\\
\psi(t)\end{pmatrix}\ran\]dt.\ea\ee

\it Proof. \rm For $\l>0$, let $\F_\l(\cd\,,\cd)$ and
$\bar\F_\l(\cd\,,\cd)$ be the evolution operators generated by
$A_\l-M(\cd)$ and $A_\l-\bar M(\cd)$, respectively. Define
\bel{L1b*}P_\l(t)\1n=\1n\F_\l(T,t)^*P(T)\F_\l(T,t)+\int_t^T\F_\l(s,t)^*
Q_0(s)\F_\l(s,t)ds,\qq t\in[0,T],\ee
and
\bel{L2b*}\bar P_\l(t)=\bar\F_\l(t,0)\bar
P(0)\bar\F_\l(t,0)^*-\2n\int_0^t\2n\bar\F_\l(t,s)\bar
Q_0(s)\bar\F_\l(t,s)^*ds,\qq t\in[0,T].\ee
For the case $A^*=A$, we have (\ref{G}) with (\ref{PsiA=APsi}) which
leads to
$$\G(t)A_\l=A_\l\G(t),\qq t\in[0,T].$$
For the case $A^*=-A$, we define
\bel{}\G_\l(t)=\1n\bar
\F_\l(T,t)^{-1}\G(T)\F_\l(T,t)\1n+\3n\int_t^T\3n\bar\F_\l(s,t)^{-1}\Th(s)
\F_\l(s,t)ds,\qq t\in[0,T].\ee
A direct computation shows that
$$\dot\Pi_\l(t)+\Pi_\l(t)[\dbA_\l-\dbM(t)]+[\dbA_\l-\dbM(t)]^*\Pi_\l(t)
+\dbQ(t)=0,$$
where
$$\dbA_\l=\begin{pmatrix}\sc
A_\l&\sc0\\\sc0&\sc-A^*_\l\end{pmatrix},\q\dbQ(t)=\begin{pmatrix}\sc Q_0(t)&
\sc\Th(t)^*\\
\sc\Th(t)&\sc\bar Q_0(t)\end{pmatrix},$$
and
$$\Pi_\l(t)=\begin{pmatrix}
\sc P_\l(t)&\sc\G_\l(t)^*\\
\sc\G_\l(t)&\sc\bar P_\l(t)\end{pmatrix},$$
with $\G_\l(\cd)=\G(\cd)$ for the case $A^*=A$. At the same time, we
let
$$\left\{\2n\ba{ll}
\ns\ds\begin{pmatrix}\sc\dot y_\l(t)\\
\sc\dot\psi_\l(t)\end{pmatrix}=\big[\dbA_\l+\dbB(t)\big]
\begin{pmatrix}\sc y_\l(t)\\
\sc\psi_\l(t)\end{pmatrix}+\begin{pmatrix}\sc b_0(t)\\ \sc-g_0(t)
\end{pmatrix},\qq t\in[0,T],\\ [2mm]
\ns\ds y_\l(0)=y(0),\qq\psi_\l(T)=\psi(T).\ea\right.$$
Then
$$\ba{ll}
\ns\ds\lan\Pi_\l(T)\begin{pmatrix}y_\l(T)\\
\psi_\l(T)\end{pmatrix},\begin{pmatrix}y_\l(T)\\
\psi_\l(T)\end{pmatrix}\ran-\lan\Pi_\l(0)\begin{pmatrix}y_\l(0)\\
\psi_\l(0)\end{pmatrix},\begin{pmatrix}y_\l(0)\\
\psi_\l(0)\end{pmatrix}\ran\\
\ns\ds=\2n\int_0^T\2n\[\lan-\(\Pi_\l(t)[\dbA_\l\1n-\1n\dbM(t)]\1n
+\1n[\dbA_\l\1n-\1n\dbM(t)]^*\Pi_\l(t)
\1n+\1n\dbQ(t)\)\begin{pmatrix}y_\l(t)\\
\psi_\l(t)\end{pmatrix},\begin{pmatrix}y_\l(t)\\
\psi_\l(t)\end{pmatrix}\ran\\
\ns\ds\qq+\1n\lan\big\{\Pi_\l(t)[\dbA_\l+\dbB(t)]\1n+\1n[\dbA_\l
+\dbB(t)]^*\Pi_\l(t)\big\}
\begin{pmatrix}y_\l(t)\\
\psi_\l(t)\end{pmatrix},\begin{pmatrix}y_\l(t)\\
\psi_\l(t)\end{pmatrix}\ran\\
\ns\ds\qq+2\lan\Pi_\l(t)\begin{pmatrix}b_0(t)\\
-g_0(t)\end{pmatrix},\begin{pmatrix}y_\l(t)\\
\psi_\l(t)\end{pmatrix}\ran\]dt\\
\ns\ds=\2n\int_0^T\2n\[\lan\1n\(\Pi_\l(t)[\dbB(t)\1n+\1n\dbM(t)]\1n
+\1n[\dbB(t)\1n+\1n\dbM(t)]^*\Pi_\l(t)\1n-\1n\dbQ(t)\)
\begin{pmatrix}y_\l(t)\\
\psi_\l(t)\end{pmatrix},\begin{pmatrix}y_\l(t)\\
\psi_\l(t)\end{pmatrix}\ran\\
\ns\ds\qq+2\lan\Pi_\l(t)\begin{pmatrix}b_0(t)\\
-g_0(t)\end{pmatrix},\begin{pmatrix}y_\l(t)\\
\psi_\l(t)\end{pmatrix}\ran\]dt.\ea$$
Passing to the limit, we obtain our result. \endpf

\ms

\ms
 Next, we let
$$\cG_0=\Big\{(b_0,g_0,h_0)\bigm|b_0(\cd),g_0(\cd)\in L^2(0,T;X),~
h_0\in X\Big\}.$$
For any $x\in X$ and $(b,g,h)\in\cG_1$, $(b_0,g_0,h_0)\in\cG_0$, and
$\rho\in[0,1]$, consider the following FBEE: It is easy to see that
when
\bel{rho1}\left\{\2n\ba{ll}
\ns\ds\3n\ba{ll}\dot y^\rho(t)=Ay^\rho(t)+\rho
b(t,y^\rho(t),\psi^\rho(t))
+b_0(t),\\
\ns\ds\dot\psi^\rho(t)=-A^*\psi^\rho(t)-\rho
g(t,y^\rho(t),\psi^\rho(t))
-g_0(t),\ea\qq t\in[0,T],\\
\ns\ds y^\rho(0)=x,\qq\psi^\rho(T)=\rho
h(y^\rho(T))+h_0.\ea\right.\ee
It is easy to see that for $\rho=0$, (\ref{rho1}) is a trivial
decoupled FBEE which admits a unique mild solution, and for
$\rho=1$, (\ref{rho1}) is essentially the same as (although it looks
a little more general than) FBEE (1.1). We will show that under
certain conditions, there exists an absolute constant $\e>0$ such
that when (\ref{rho1}) is (uniquely) solvable for some
$\rho\in[0,1)$, it must be (uniquely) solvable for (\ref{rho1}) with
$\rho$ replaced by $(\rho+\e)\land1$. Then by repeating the same
argument, we obtain the (unique) solvability of (\ref{FBEE1}) over
$[0,T]$. Such an argument is called a {\it method of continuation}
(see \cite{Yong 1997}). In doing so, the key is to establish an {\it
a priori} estimate for the mild solutions to (\ref{rho1}), uniform
in $\rho\in[0,1]$. To this end, we need to make some preparations.

\ms

For any $\l>0$, we introduce the following approximate system of
FBEE (\ref{rho1}):
\bel{rho2}\left\{\2n\ba{ll}
\ns\ds\3n\ba{ll}\dot y^\rho_\l(t)=A_\l y^\rho_\l(t)+\rho
b(t,y^\rho_\l(t),\psi^\rho_\l(t))+b_0(t),\\
\ns\ds\dot\psi^\rho_\l(t)=-A_\l^*\psi^\rho_\l(t)-\rho
g(t,y^\rho_\l(t),\psi^\rho_\l(t))-g_0(t),\ea\qq t\in[0,T],\\
\ns\ds y^\rho_\l(0)=x,\qq\psi^\rho_\l(T)=\rho
h(y^\rho_\l(T))+h_0.\ea\right.\ee
Suppose for the initial condition $x\in X$, the generator
$(b,g,h)\in\cG_4$ and $(b_0,g_0,h_0)\in\cG_0$, FBEE (\ref{rho2})
admits a solution $(y_\l^\rho(\cd),\psi_\l^\rho(\cd))$. Also, let
$(\bar y_\l^\rho(\cd),\bar \psi_\l^\rho(\cd))$ be a solution of
$(\ref{rho2})$ with $(b,g,h)$, $(b_0,g_0,h_0)$, and $x$ respectively
replaced by $(\bar b,\bar g,\bar h)\in\cG_1$, $(\bar b_0,\bar
g_0,\bar h_0)\in\cG_0$, and $\bar x\in X$. Define
\bel{4.4}\h y(\cd)=\bar y_\l^\rho(\cd)-y_\l^\rho(\cd),\q
\h\psi(\cd)=\bar\psi_\l^\rho(\cd)-\psi_\l^\rho(\cd).\ee
Denote
\bel{4.8}\left\{\2n\ba{ll}
\ns\ds\wt b_y(t)\1n=\2n\int_0^1b_y(t,y_\l^\rho(t)+\a\h y(t),
\psi_\l^\rho(t)+\a\h\psi(t))d\a,\\
\ns\ds\wt b_\psi(t)\1n=\2n\int_0^1b_\psi(t,y_\l^\rho(t)+\a\h
y(t),\psi_\l^\rho(t)+\a\h\psi(t)))d\a,\\
\ns\ds\wt g_y(t)\1n=\2n\int_0^1g_y(t,y_\l^\rho(t)+\a\h y(t),
\psi_\l^\rho(t)+\a\h\psi(t))d\a,\\
\ns\ds\wt g_\psi(t)\1n=\2n\int_0^1g_\psi(t,y_\l^\rho(t)+\a\h
y(t),\psi_\l^\rho(t)+\a\h\psi(t)))d\a,\\
\ns\ds\wt h_y=\int_0^1h_y(y_\l^\rho(T)+\a\h
y(T))d\a,\ea\right.\ee
and set
\bel{5.30}\left\{\2n\ba{ll}
\ns\ds\d b(t)=\bar b(t,\bar
y_\l^\rho(t),\bar\psi_\l^\rho(t))-b(t,\bar
y_\l^\rho(t),\bar\psi_\l^\rho(t)),\q\d b_0(t)=\bar
b_0(t)-b_0(t),\\
\ns\ds\d g(t)=\bar g(t,\bar
y_\l^\rho(t),\bar\psi_\l^\rho(t))-g(t,\bar
y_\l^\rho(t),\bar\psi_\l^\rho(t)),\q\d g_0(t)=\bar
g_0(t)-g_0(t),\\
\ns\ds\d h=\bar h(\bar y_\l^\rho(T))-h(\bar y_\l^\rho(T)),\q\d
h_0=\bar h_0-h_0,\q\h x=\bar x-x.\ea\right.\ee
Then $(\h y(\cd),\h\psi(\cd))$ satisfies
\bel{rho h}\left\{\ba{ll}
\ns\ds\3n\ba{ll}\dot{\2n\h y}(t)=A_\l\h y(t)+\rho\wt b_y(t)\h y(t)
+\rho\wt b_\psi(t)\h\psi(t)+\rho\d b(t)+\d b_0(t),\\
\ns\ds\dot{\3n\h\psi}(t)=-A_\l^*\h\psi(t)\1n-\1n\rho\wt g_y(t)\h
y(t)
\1n-\1n\rho\wt g_\psi(t)\h\psi(t)\1n-\1n\rho\d g(t)\1n-\1n\d g_0(t),
\ea\qq t\in[0,T],\\
\ns\ds\2n\h y(t)=\h x,\qq\h\psi(T)=\rho\wt h_y\h y(T)+\d h+\d
h_0.\ea\right.\ee
For the above linear FBEE, we have the following result.

\ms

\bf Proposition 5.2. \sl Let $(b,g,h)\in\cG_3$ and
\bel{l(bg)}\left\{\2n\ba{ll}
\ns\ds L_{by}(t)\1n\deq\sup_{(y,\psi)\in X\times X}\[\max\si\(
{b_y(t,y,\psi)+b_y(t,y,\psi)^*\over2}\)\]^+,\\
\ns\ds L_{g\psi}(t)\deq\sup_{(y,\psi)\in X\times X}\[\max\si\({
g_\psi(t,y,\psi)+g_\psi(t,y,\psi)^*\over2}\)\]^+,\ea\right.\q
t\in[0,T],\ee
where $\si(\L)$ is the spectrum of the operator $\L\in\cL(X)$. Let
$(y^\rho_\l(\cd),\psi^\rho_\l(\cd))$ be a solution to FBEE
$(\ref{rho2})$, and $(\bar y^\rho_\l(\cd),\bar\psi^\rho_\l(\cd))$ be
a solution of {\rm(\ref{rho2})} corresponding to $(\bar b,\bar g,
\bar h)\in\cG_1$, $(b_0,g_0,h_0)\in\cG_0$ and $\bar x\in X$. Then
\bel{|hy|}\ba{ll}
\ns\ds\|\h
y(\cd)\|_\infty\les\rho\int_0^Te^{\rho\int_s^TL_{by}(\t)d\t} \|\wt
b_\psi(s)\h\psi(s)\|ds+K\[\|\h x\|+\int_0^T\(\|\d b(s)\|+\|\d
b_0(s)\|\)ds\],\ea\ee
and
\bel{|hpsi|}\ba{ll}
\ns\ds\|\h\psi(\cd)\|_\infty\les\rho\[e^{\rho\int_0^TL_{g\psi}(\t)d\t}\|\wt
h_y\h y(T)\|
+\int_t^T\3n e^{\rho\int_0^sL_{g\psi}(\t)d\t}\|\wt g_y(s)\h y(s)\|ds\]\\
\ns\ds\qq\qq\qq\qq+K\[\|\d h\|+\|\d h_0\|+\int_0^T\(\|\d g(s)\|+\|\d
g_0(s)\|\)ds\].\ea\ee

\ms

\rm

The proof is straightforward and for reader's convenience, a proof
is presented in the appendix.

\ms

We note that in the above proposition, it is only assumed that
$(b,g,h)\in \cG_3$ (the set of all generators satisfying (H3)).
Therefore, the Fr\'echet derivatives $b_y,b_\psi$, and so on are not
necessarily bounded. However, it is still possible that
$$\int_0^T\(\|\wt b_\psi(s)\h\psi(s)\|+\|\wt g_y(s)\h y(s)\|\)ds<\infty.$$
On the other hand, in the case $(b,g,h)\in\cG_4$, we have
$$\ba{ll}
\ns\ds\|\h
y(\cd)\|_\infty\les\rho\int_0^Te^{\rho\int_s^TL_{by}(\t)d\t} \|\wt
b_\psi(s)\h\psi(s)\|ds+K\[\|\h x\|+\int_0^T\(\|\d b(s)\|+\|\d
b_0(s)\|\)ds\]\\
\ns\ds\qq\qq\les\rho\(\int_0^Te^{\rho\int_s^TL_{by}(\t)d\t}\|\wt
b_\psi(s)\|ds\) \|\h\psi(\cd)\|_\infty+K\[\|\h x\|+\int_0^T\(\|\d
b(s)\|+\|\d b_0(s)\|\)ds\],\ea$$
and
$$\ba{ll}
\ns\ds\|\h\psi(\cd)\|_\infty\les\rho\[e^{\rho\int_0^TL_{g\psi}(\t)d\t}\|\wt
h_y\h y(T)\|
+\int_0^T\3n e^{\rho\int_0^sL_{g\psi}(\t)d\t}\|\wt g_y(s)\h y(s)\|ds\]\\
\ns\ds\qq\qq\qq\qq+K\[\|\d h\|+\|\d h_0\|+\int_0^T\(\|\d g(s)\|+\|\d
g_0(s)\|\)ds\]\\
\ns\ds\les\1n\rho^2\1n\[e^{\rho\1n\int_t^T\2n L_{g\psi}(\t)d\t}\|\wt
h_y\|\1n+\2n\int_0^T\3n e^{\rho\1n\int_0^s\2n L_{g\psi}(\t)d\t}\|\wt
g_y(s)\|ds\]\[\1n\int_0^T\2n e^{\rho\1n\int_s^T\2n
L_{by}(\t)d\t}\|\wt
b_\psi(s)\|ds\]\|\h\psi(\cd)\|_\infty\\
\ns\ds\q+K\[\|\h x\|\1n+\1n\|\d h\|\1n+\1n\|\d
h_0\|\1n+\2n\int_0^T\3n\(\|\d b(s)\|\1n+\1n\|\d b_0(s)\|\1n+\1n\|\d
g(s)\|\1n+\1n\|\d g_0(s)\|\)ds\].\ea$$
Hence, when the following holds:
\bel{<1}\ba{ll}
\ns\ds\rho^2\1n\[e^{\rho\1n\int_0^T\2n L_{g\psi}(\t)d\t}\|\wt h_y\|
\1n+\2n\int_0^T\3n e^{\rho\1n\int_0^s\2n L_{g\psi}(\t)d\t}\|\wt
g_y(s)\|ds\]\[\1n\int_0^T\2n e^{\rho\1n\int_s^T\2n
L_{by}(\t)d\t}\|\wt b_\psi(s)\|ds\]\2n<\1n1,\ea\ee
FBEE (\ref{rho2}) admits a unique mild solution
$(y_\l^\rho(\cd),\psi_\l^\rho(\cd))$, by means of contraction
mapping theorem. It is not hard to see that condition (\ref{<1})
holds when one of the following holds:

\ms

$\bullet$ The parameter $\rho=0$, this is a trivial case, for which
the FBEE is linear and decoupled.

\ms

$\bullet$ The time duration $T$ is small enough.

\ms

$\bullet$ The coupling is weak enough in the sense that the
Lipschitz constant of $b(t,y,\psi)$ with respect to $\psi$ (the
bound of $b_\psi(\cd)$), and/or the Lipschitz constants of
$g(t,y,\psi)$ and $h(y)$ with respect to $y$ (the bounds of
$g_y(\cd)$ and $h_y(\cd)$) are small enough. An extreme case is that
$b(t,y,\psi)$ is independent of $\psi$, or $g(t,y,\psi)$ and $h(y)$
are independent of $y$, which corresponds to the decoupled case.

\ms

From Proposition 5.2, we see that due to the coupling, in general,
one can only obtain an estimate of $\h y(\cd)$ in terms of
$\h\psi(\cd)$, and an estimate of $\h\psi(\cd)$ in terms of $\h
y(\cd)$. In order to obtain an a priori estimate on the whole $(\h
y(\cd),\h\psi(\cd))$, we need either have an estimate for
$$\int_0^T\|\wt b_\psi(s)\h\psi(s)\|^2ds$$
independent of $\h y(\cd)$, or have an estimate for
$$\|\wt h_y\h y(T)\|^2+\int_0^T\|\wt g_y(s)\h y(s)\|^2ds$$
independent of $\h\psi(\cd)$. We now search conditions under which
this is possible. To this end, we introduce the following notions.

\ms

\bf Definition 5.3. \rm A continuous function
$\Pi(\cd)\equiv\begin{pmatrix}\sc P(\cd) &\sc\G(\cd)^*\\
\sc\G(\cd)&\sc\bar P(\cd)\end{pmatrix}:[0,T]\to\dbS(X\times X)$ is
called a type (I) {\it Lyapunov operator} of the generator
$(b,g,h)\in\cG_3$ if there exist $\dbQ:[0,T]\to\dbS(X\times X)$ and
$\dbM:[0,T]\to\cL(X\times X)$ with
$$\dbQ(t)=\begin{pmatrix}\sc Q_0(t)&\sc\Th(t)^*\\
\sc\Th(t)&\sc\bar
Q_0(t)\end{pmatrix},\q\dbM(t)=\begin{pmatrix}\sc M(t)&\sc0\\
\sc0&\sc-\bar M(t)^*\end{pmatrix},\q t\in[0,T],$$
such that $\Pi(\cd)$ is a mild solution to the Lyapunov differential
equation
\bel{Pi*}\dot\Pi(t)+\Pi(t)[\dbA-\dbM(t)]+[A-\dbM(t)]^*\Pi(t)+\dbQ(t)=0,\q
t\in[0,T].\ee
and for some constants $\m,K>0$, the following are satisfied:
\bel{I}\left\{\2n\ba{ll}
\ns\ds\;\Pi(0)+\begin{pmatrix} -K&0\cr0&0\end{pmatrix}\les0,\\
[4mm]
\ns\ds\begin{pmatrix}I&\rho h_y(y)^*\\
0&I\end{pmatrix}\1n\Pi(T)\1n\begin{pmatrix}I&0\\
\rho h_y(y)&I\end{pmatrix}\2n+\2n\begin{pmatrix}-\m
h_y(y)^*h_y(y)&0\cr0&K\end{pmatrix}\2n\ges\1n0,\qq\forall y\1n\in\1n
X,~\rho\1n\in\1n[0,1],\\ [5mm]
\ns\ds\1n\begin{pmatrix}\dbH^\rho(t,\Pi(t),y,\psi)-\dbQ(t)
+\m\begin{pmatrix}
g_y(t,y,\psi)^*g_y(t,y,\psi)&0\cr0&0\end{pmatrix}&\Pi(t)\cr
\Pi(t)&-K\end{pmatrix}\2n\les\1n0,\\ [5mm]
\ns\ds\qq\qq\qq\qq\qq\qq\qq\forall(t,y,\psi)\in[0,T]\times X\times
X,~\rho\in[0,1],\ea\right.\ee
where
\bel{Hrho}\ba{ll}
\ns\ds\dbH^\rho(t,\Pi,y,\psi)=\rho\big[\Pi\dbB(t,y,\psi)+\dbB(t,y,\psi)^*
\Pi\big]+\Pi\dbM(t)+\dbM(t)^*\Pi,\\
[1mm]
\ns\ds\qq\qq\qq\qq\forall(t,\Pi,y,\psi)\in[0,T]\times\dbS(X)\times
X\times X,~\rho\in[0,1],\ea\ee
and
$$\dbB(t,y,\psi)=\begin{pmatrix}\sc b_y(t,y,\psi)&\sc b_\psi(t,y,\psi)\\
                                \sc-g_y(t,y,\psi)&\sc-g_\psi(t,y,\psi)\end{pmatrix},
\qq\forall(t,y,\psi)\in[0,T]\times X\times X.$$
If (\ref{I}) is replaced by the following:
\bel{II}\left\{\2n\ba{ll}
\ns\ds\Pi(0)+\begin{pmatrix}-K&0\cr0&0\end{pmatrix}\les0,\\
[3mm]
\ns\ds\begin{pmatrix}I&\rho h_y(y)^*\\
0&I\end{pmatrix}\Pi(T)\begin{pmatrix}I&0\\
\rho
h_y(y)&I\end{pmatrix}+\begin{pmatrix}0&0\cr0&K\end{pmatrix}\ges0,\q\forall
y\in X,~\rho\in[0,1],\\ [5mm]
\ns\ds\begin{pmatrix}\dbH^\rho(t,\Pi(t),y,\psi)-\dbQ(t)
+\m\begin{pmatrix}0&0\cr0&b_\psi(t,y,\psi)^*
b_\psi(t,y,\psi)\end{pmatrix}&\Pi(t)\cr\Pi(t)&-K\end{pmatrix}\2n\les\1n0,\\
[5mm]
\ns\ds\qq\qq\qq\qq\qq\qq\qq\forall(t,y,\psi)\in[0,T]\times X\times
X,~\rho\in[0,1],\ea\right.\ee
then $\Pi(\cd)$ is called a type (II) {\it Lyapunov operator} of
$(b,g,h)$.

\ms

If $\Pi(\cd)$ is either a type (I) or Type (II) Lyapunov operator of
$(b,g,h)$, we simply call it a Lyapunov operator of $(b,g,h)$.

\ms

The existence of a Lyapunov operator gives some kind of
compatibility of the coefficients in FBEE (\ref{FBEE1}), which will
guarantee the well-posedness of the FBEE. We will carefully discuss
properties and existence of Lyapunov operators a little later.
First, we present the following result gives the (uniform) stability
of mild solutions to (\ref{rho1}) when the generator
$(b,g,h)\in\cG_3$ admits a Lyapunov operator.

\ms

\bf Proposition 5.4. \sl Let $(b,g,h)\in\cG_3$ admit a Lyapunov
operator $\Pi(\cd)$ of either type (I) or (II). For any
$\rho\in[0,1]$, and $x\in X$, let $(y^\rho(\cd),\psi^\rho(\cd))$ be
a mild solution of FBEE $(\ref{rho1})$ with some
$(b_0(\cd),g_0(\cd),h_0)\in\cG_0$, and let $(\bar
y^\rho(\cd),\bar\psi^\rho(\cd))$ be a mild solutions of FBEE
$(\ref{rho1})$ corresponding to another generator $(\bar b,\bar
g,\bar h)\in\cG_1$, and some $(\bar b_0(\cd),\bar g_0(\cd),\bar
h_0)\in\cG_0$, $\bar x\in X$. Then
\bel{stability1}\2n\2n\ba{ll}
\ns\ds\|\bar
y^\rho(\cd)-y^\rho(\cd)\|_\infty+\|\bar\psi^\rho(\cd)-\psi^\rho(\cd)\|_\infty
\les
K\Big\{\|\bar x-x\|^2\2n+\|\bar h_0-h_0\|^2
+\|\bar h(\bar y^\rho(T))-h(\bar y^\rho(T))\|^2\\
\ns\ds\qq+\int_0^T\Big(\|\bar b(s,\bar y^\rho(s),\bar\psi^\rho(s))-
b(t,\bar y^\rho(t),\bar\psi^\rho(s))|^2+\|\bar b_0(s)-b_0(s)\|^2\\
\ns\ds\qq\q+\|\bar g(s,\bar y^\rho(s),\bar\psi^\rho(s))-g(s,\bar
y^\rho(s),\bar\psi^\rho(s))\|^2+\|\bar
g_0(s)-g_0(s)\|^2\Big)ds\Big\},\ea\ee
uniformly in $\rho\in[0,1]$. In particular, if $(\bar b,\bar g,\bar
h)=(b,g,h)$, then
\bel{stability2}\2n\2n\ba{ll}
\ns\ds\|\bar y^\rho(\cd)-y^\rho(\cd)\|_\infty^2+\|\bar\psi^\rho(\cd)-\psi^\rho(\cd)
\|_\infty^2\\
\ns\ds\les\1n K\Big\{\|\bar x-x\|^2\2n+\1n\|\bar
h_0-h_0\|^2\2n+\2n\int_0^T\2n\Big(\|\bar b_0(s)\1n-\1n
b_0(s)\|^2\2n+\1n\|\bar g_0(t)\1n-\1n g_0(s)\|^2\Big)ds\Big\},\ea\ee
uniformly in $\rho\in[0,1]$.

\ms

\it Proof. \rm Recall notations in (\ref{4.4})--(\ref{5.30}) and
noting Proposition 5.1, we have (suppressing $s$ when it has no
ambiguity)
$$\ba{ll}
\ns\ds\lan\begin{pmatrix}\sc I&\sc\rho(\wt h_y)^*\\
\sc0&\sc I\end{pmatrix}\Pi(T)\begin{pmatrix}\sc I&\sc0\\
\sc\rho\wt h_y&\sc I\end{pmatrix}\begin{pmatrix}\sc\h
y(T)\cr\sc\rho\d h+\d h_0\end{pmatrix},\begin{pmatrix}\sc\h
y(T)\cr\sc\rho\d h+\d
h_0\end{pmatrix}\ran-\lan\Pi(0)\begin{pmatrix}\sc\h x\cr\sc\h
\psi(0)\end{pmatrix},\begin{pmatrix}\sc\h x\cr\sc\h\psi(0)
\end{pmatrix}\ran\\ [4mm]
\ns\ds=\lan\Pi(T)\begin{pmatrix}\sc\h y(T)\cr\sc\rho\wt h_y\h
y(T)+\rho\d h+\d h_0\end{pmatrix},\begin{pmatrix}\sc\h
y(T)\cr\sc\rho\wt h_y\h y(T)+\rho\d h+\d
h_0\end{pmatrix}\ran-\lan\Pi(0)\begin{pmatrix}\sc\h x\cr\sc\h
\psi(0)\end{pmatrix},\begin{pmatrix}\sc\h x\cr\sc\h\psi(0)
\end{pmatrix}\ran\\ [4mm]
\ns\ds=\lan\Pi(T)\begin{pmatrix}\sc\h
y(T)\cr\sc\h\psi(T)\end{pmatrix},\begin{pmatrix}\sc\h
y(T)\cr\sc\h\psi(T)\end{pmatrix}\ran-\lan \Pi(0)\begin{pmatrix}\sc\h
x\cr\sc\h\psi(0)\end{pmatrix},\begin{pmatrix}\sc\h x\cr\sc\h\psi(0)
\end{pmatrix}\ran\\ [4mm]
\ns\ds=\2n\int_0^T\3n
\Big\{\1n\lan\1n\(\wt\dbH^\rho-\dbQ\)\begin{pmatrix}\sc\h
y\cr\sc\h\psi\end{pmatrix},\begin{pmatrix}\sc\h y\cr\sc\h
\psi\end{pmatrix}\ran+2\lan\Pi\begin{pmatrix}\sc\h
y\cr\sc\h\psi\end{pmatrix},\begin{pmatrix}\sc\rho\d b+\d
b_0\cr\sc\rho\d g+\d g_0\end{pmatrix}\ran\Big\}ds\\ [4mm]
\ns\ds=\int_0^T\lan\begin{pmatrix}\sc\wt\dbH^\rho-\dbQ&\sc\Pi\cr\sc
\Pi&\sc0\end{pmatrix}\begin{pmatrix}\sc\h y\cr\sc\h\psi\cr\sc\rho\d
b+\d b_0\cr\sc\rho\d g+\d g_0\end{pmatrix},\begin{pmatrix}\sc\h
y\cr\sc\h\psi\cr\sc\rho\d b+\d b_0\cr\sc\rho\d g+\d
g_0\end{pmatrix}\ran ds,\ea$$
where
$$\wt\dbH^\rho(t,\Pi,y,\psi)=\rho\big[\Pi\wt\dbB+\wt\dbB^*
\Pi\big]+\Pi\dbM(t)+\dbM(t)^*\Pi,$$
and
$$\wt\dbB=\begin{pmatrix}\wt b_y&\wt b_\psi\\ -\wt g_y&-\wt g_\psi\end{pmatrix},$$
with $\wt b_y$, etc. given by (\ref{4.8}). Consequently, in the case
that $\Pi(\cd)$ is a type (I) Lyapunov operator of $(b,g,h)$, we
have
$$\ba{ll}
\ns\ds\lan\Pi(T)\begin{pmatrix}\sc\h y(T)\cr\sc\rho\wt h_y\h
y(T)+\rho\d h+\d h_0\end{pmatrix},\begin{pmatrix}\sc\h
y(T)\cr\sc\rho\wt h_y\h y(T)+\rho\d h+\d h_0\end{pmatrix}\ran-\lan
\Pi(0)\begin{pmatrix}\sc\h x\cr\sc\h
\psi(0)\end{pmatrix},\begin{pmatrix}\sc\h x\cr\sc\h\psi(0)\end{pmatrix}\ran\\
\ns\ds\ges\m\|\wt h_y\h y(T)\|^2-K\[\|\h x\|^2+\|\d h\|^2+\|\d
h_0\|^2\],\ea$$
and
$$\ba{ll}
\ns\ds\int_0^T\lan\begin{pmatrix}\sc\wt\dbH^\rho-\dbQ&\sc\Pi\cr\sc
\Pi&\sc0\end{pmatrix}\begin{pmatrix}\sc\h y\cr\sc\h\psi\cr\sc\rho\d
b+\d b_0\cr\sc\rho\d g+\d g_0\end{pmatrix},\begin{pmatrix}\sc\h
y\cr\sc\h\psi\cr\sc\rho\d b+\d b_0\cr\sc\rho\d g+\d
g_0\end{pmatrix}\ran ds\\
\ns\ds\les\int_0^T\(-\m\|\wt g_y(s)\h y(s)\|^2+K\|\rho\d b(s)+\d
b_0(s)\|^2+K\|\rho\d g(s)+\d g_0(s)\|^2\)ds.\ea$$
Hence,
$$\ba{ll}
\ns\ds\|\wt h_y\h y(T)\|^2+\int_0^T\|\wt g_y(s)\h y(s)\|^2ds\les K\[\|\h x\|^2+\|\d h\|^2+\|\d h_0\|^2\\
\ns\ds\qq\qq\qq+\int_0^T\(\|\d b(s)\|^2+\|\d b_0(s)\|^2+\|\d
g(s)\|^2+\|\d g_0(s)\|^2\)ds\].\ea$$
Combining the above with (\ref{|hy|}) and (\ref{|hpsi|}), we obtain
(\ref{stability1}).

\ms

On the other hand, in the case that $\Pi(\cd)$ is a type (II)
Lyapunov operator for $(b,g,h)$, we have
$$\ba{ll}
\ns\ds\lan\Pi(T)\begin{pmatrix}\h y(T)\cr\rho\wt h_y\h y(T)+\rho\d
h+\d h_0\end{pmatrix},\begin{pmatrix}\h y(T)\cr\rho\wt h_y\h
y(T)+\rho\d h+\d h_0\end{pmatrix}\ran-\lan \Pi(0)\begin{pmatrix}\h
x\cr\h \psi(0)\end{pmatrix},\begin{pmatrix}\h
x\cr\h\psi(0)\end{pmatrix}\ran\\ [4mm]
\ns\ds\ges-K\[\|\h x\|^2+\|\d h\|^2+\|\d h_0\|^2\],\ea$$
and
$$\ba{ll}
\ns\ds\int_0^T\lan\begin{pmatrix}\sc\wt\dbH^\rho-\dbQ&\sc\Pi\cr\sc
\Pi&\sc0\end{pmatrix}\begin{pmatrix}\sc\h y\cr\sc\h\psi\cr\sc\rho\d
b+\d b_0\cr\sc\rho\d g+\d g_0\end{pmatrix},\begin{pmatrix}\sc\h
y\cr\sc\h\psi\cr\sc\rho\d b+\d b_0\cr\sc\rho\d g+\d
g_0\end{pmatrix}\ran ds\\
\ns\ds\les\int_0^T\(-\m\|\wt b_\psi(s)\h\psi(s)\|^2+K\|\rho\d
b(s)+\d b_0(s)\|^2+K\|\rho\d g(s)+\d g_0(s)\|^2\)ds.\ea$$
Hence,
$$\ba{ll}
\ns\ds\int_0^T\|\wt b_\psi(s)\h\psi(s)\|^2ds\les K\[\|\h x\|^2+\|\d h\|^2+\|\d h_0\|^2\\
\ns\ds\qq\qq\qq\qq\qq\qq+\int_0^T\(\|\d b(s)\|^2+\|\d
b_0(s)\|^2+\|\d g(s)\|^2+\|\d g_0(s)\|^2\)ds\].\ea$$
Then, combining the above with (\ref{|hy|}) and (\ref{|hpsi|}), we
again obtain (\ref{stability1}). \endpf

\section{Well-Posedness of FBEEs via Lyapunov Operators}

We now state and prove the following theorem concerning the
well-posedness of FBEE (\ref{FBEE1}).

\ms

\bf Theorem 6.1. \sl Let $(b,g,h)\in\cG_2\cap\cG_3$ admit a type
{\rm(I)} or {\rm(II)} Lyapunov operator $\Pi(\cd)$. Then FBEE
$(\ref{FBEE1})$ admits a unique mild solution $(y(\cd),\psi(\cd))$.
Moreover, the following estimate holds:
\bel{2.38}\ba{ll}
\ns\ds\|y(\cd)\|_\infty+\|\psi(\cd)\|_\infty\les
K\[\|x\|+\|h(0)\|+\int_0^T\(\|b(s,0,0)\|+\|g(s,0,0)\|\)ds\].\ea\ee
Further, if $(\bar y(\cd),\bar\psi(\cd))$ is a mild solution of FBEE
$(1.1)$ corresponding to $(\bar b,\bar g,\bar h)\in\cG_2\cap\cG_3$,
then the following stability estimate holds:
\bel{stability3}\ba{ll}
\ns\ds\|\bar y(\cd)-y(\cd)\|_\infty+\|\bar
\psi(\cd)-\psi(\cd)\|_\infty\les K\Big\{\|\bar x-x\|+\|\bar h(\bar y(T))-h(\bar y(T))\|\\
\ns\ds\qq\qq\qq+\int_0^T\(\|\bar b(s,\bar y(s),\bar\psi(s))-b(s,\bar
y(s),\bar\psi(s))\|+\|\bar g(s,\bar y(s),\bar \psi(s)-g(s,\bar
y(s),\bar\psi(s)\|\)ds\Big\}.\ea\ee

\it Proof. \rm Let $(b_0,g_0,h_0)\in\cG_0$. Let $\rho\in[0,1)$.
Suppose the following (coupled) FBEE admits a unique mild solution
$(y^\rho(\cd),\psi^\rho(\cd))$:
\bel{lambda3}\left\{\2n\ba{ll}
\ns\ds\dot y^\rho(t)=Ay^\rho(t)+\rho
b(t,y^\rho(t),\psi^\rho(t))+b_0(t),\\
\ns\ds\dot\psi^\rho(t)=-A^*\psi^\rho(t)+\rho g(t,y^\rho(t),\psi^\rho(t))
+g_0(t),\\
\ns\ds y^\rho(0)=x,\qq\psi^\rho(T)=\rho
h(y^\rho(T))+h_0,\ea\right.\ee
and the following estimate holds:
\bel{}\ba{ll}
\ns\ds\|y^\rho(\cd)\|_\infty+\|\psi^\rho(\cd)\|_\infty\les K\Big\{\|x\|+\|h(0)\|+\|h_0\|\\
\ns\ds\qq\qq\qq\qq\qq\qq+\int_0^T\2n\(\|b(s,0,0)\|+\|b_0(s)\|+\|g(s,0,0)\|+\|g_0(s)\|\)ds\Big\}.\ea\ee
Now, let $\e>0$ such that $\rho+\e\in[0,1]$. Consider the following
coupled FBEE:
\bel{3.3}\left\{\2n\ba{ll}
\ns\ds
\dot y^{\rho+\e}(t)=Ay^{\rho+\e}(t)+(\rho+\e)b(t,y^{\rho+\e}(t),\psi^{\rho+\e}(t))+b_0(t),\\
\ns\ds
\dot\psi^{\rho+\e}(t)=-A^*\psi^{\rho+\e}(t)+(\rho+\e)g(t,y^{\rho+\e}(s),\psi^{\rho+\e}(t))
+g_0(t),\\
\ns\ds y^{\rho+\e}(t)=x,\qq\psi^{\rho+\e}(T)=(\rho+\e)
h(y^{\rho+\e}(T))+h_0,\ea\right.\ee
To obtain the (unique) solvability of the above problem, we
introduce the following sequence of problems:
\bel{}\left\{\2n\ba{ll}
\ns\ds y^{\rho+\e,0}(\cd)=\psi^{\rho+\e,0}(\cd)=0,\\
\ns\ds\dot y^{\rho+\e,k+1}(t)=Ay^{\rho+\e,k+1}(t)+\rho
b(t,y^{\rho+\e,k+1}(t),\psi^{\rho+\e,k+1}(t))\\
\ns\ds\qq\qq\qq\q+\e
b(t,y^{\rho+\e,k}(t),\psi^{\rho+\e,k}(t))+b_0(t),\\
\ns\ds\dot\psi^{\rho+\e,k+1}(t)=-A^*\psi^{\rho+\e,k+1}(t)+\rho
g(t,y^{\rho+\e,k+1}(t),\psi^{\rho+\e,k+1}(t))\\
\ns\ds\qq\qq\qq\q+\e
g(t,y^{\rho+\e,k}(t),\psi^{\rho+\e,k}(t))+g_0(t),\\
\ns\ds y^{\rho+\e,k+1}(0)=x,\\
\ns\ds\psi^{\rho+\e,k+1}(T)=\rho h(y^{\rho+\e,k+1}(T))+\e
h(y^{\rho+\e,k}(T))+h_0.\ea\right.\ee
By our assumption, inductively, for each $k\ges0$, as long as
$(y^{\rho+\e,k}(\cd),\psi^{\rho+\e,k}(\cd))\in C([t,T];X\times X)$,
the above FBEE admits a unique mild solution
$(y^{\rho+\e,k+1}(\cd),$ $\psi^{\rho+\e,k+1}(\cd))\in
C([t,T];X\times X)$. Further,
\bel{6.7*}\ba{ll}
\ns\ds\|y^{\rho+\e,k+1}(\cd)\|_\infty+\|\psi^{\rho+\e,k+1}(\cd)\|_\infty\\
\ns\ds\les\2n K\Big\{\|x\|+\|h(0)\|+\|\e
h(y^{\rho+\e,k+1}(T))+h_0\|\\
\ns\ds\qq\q+\int_0^T\2n\Big(\|b(s,0,0)\|+\|\e
b(s,y^{\rho+\e,k}(s),\psi^{\rho+\e,k}(s))+
b_0(s)\|\)ds\\
\ns\ds\qq\q+\int_0^T\2n\Big(\|g(s,0,0)\|+\|\e
g(s,y^{\rho+\e,k}(s),\psi^{\rho+\e,k}(s))+g_0(t)\|\)ds\Big\}\\
\ns\ds\les K\Big\{\|x\|+\|h(0)\|+\|h_0\|+\e
\|y^{\rho+\e,k}(T)\|\\
\ns\ds\qq\q+\int_0^T\3n\Big[\|b(s,0,0)\|+\|b_0(s)\|+\e
\(\|y^{\rho+\e,k}(s)|+\|\psi^{\rho+\e,k}(s)\|\)\]ds\\
\ns\ds\qq\q+\int_0^T\3n\[\|g(s,0,0)\|+\|g_0(s)\|
+\1n\e\(\|y^{\rho+\e,k}(s)\|+\|\psi^{\rho+\e,k}(s)\|\)\]ds\Big\}\\
\ns\ds\les
K\e\(\|y^{\rho+\e,k}(\cd)\|_\infty+\|\psi^{\rho+\e,k}(\cd)\|_\infty\)\\
\ns\ds\qq+K\Big\{\|x\|
+\|h(0)\|+\|h_0\|+\int_0^T\3n\(\|b(s,0,0)\|+\|b_0(s)\|+
\|g(s,0,0)\|+\|g_0(s)\|\)ds\Big\}.\ea\ee
Now, since $(b,g,h)$ admits a type (I) or type (II) Lyapunov
operator $P(\cd)$, by Proposition 5.3, we obtain
$$\ba{ll}
\ns\ds\|y^{\rho+\e,k+1}(\cd)-y^{\rho+\e,k}(\cd)\|_\infty+\|\psi^{\rho+\e,k+1}(\cd)-\psi^{\rho+\e,k}(\cd)
\|_\infty\\
\ns\ds\les\e K\Big\{\|h(y^{\rho+\e,k}(T))-h(y^{\rho+\e,k-1}(T))\|\\
\ns\ds\qq+\int_0^T\Big[\|
b(s,y^{\rho+\e,k}(s),\psi^{\rho+\e,k}(s))-b(s,y^{\rho+\e,k-1}(s),\psi^{\rho+\e,k-1}(s))\|\\
\ns\ds\qq\q+\|g(s,y^{\rho+\e,k}(s),\psi^{\rho+\e,k}(s))-g(s,y^{\rho+\e,k-1}(s),\psi^{\rho+\e,k-1}(s))
\|\Big]ds\Big\}\\
\ns\ds\les\e
K_0\(\|y^{\rho+\e,k}(\cd)-y^{\rho+\e,k-1}(\cd)\|_\infty+\|\psi^{\rho+\e,k}(\cd)-\psi^{\rho+\e,k-1}(\cd)
\|_\infty\).\ea$$
Here, $K_0>0$ is an absolute constant (independent of $k\ges1$).
Thus, taking $\e>0$ small enough so that $\e K_0\les{1\over2}$, we
obtain
$$\lim_{k\to\infty}\(\|y^{\rho+\e,k}(\cd)-y^{\rho+\e}(\cd)\|_\infty+\|\psi^{\rho+\e,k}(\cd)-
\psi^{\rho+\e}(\cd)\|_\infty\)=0,$$
with
$$\left\{\2n\ba{ll}
\ns\ds
y^{\rho+\e}(\cd)=\sum_{k=1}^\infty\Big[y^{\rho+\e,k}(\cd)-y^{\rho+\e,k-1}(\cd)\Big],\\
\ns\ds
\psi^{\rho+\e}(\cd)=\sum_{k=1}^\infty\Big[\psi^{\rho+\e,k}(\cd)
-\psi^{\rho+\e,k-1}(\cd)\Big].\ea\right.$$
which is the unique mild solution of FBEE (\ref{3.3}). Further, let
$k\to\infty$ in (\ref{6.7*}), we obtain
$$\ba{ll}
\ns\ds\|y^{\rho+\e}(\cd)\|_\infty+\|\psi^{\rho+\e}(\cd)\|_\infty\les
K\e\(\|y^{\rho+\e}(\cd)\|_\infty+\|\psi^{\rho+\e}(\cd)\|_\infty\)\\
\ns\ds\qq\qq+K\Big\{
\|x\|+\|h(0)\|+\|h_0\|+\int_t^T\3n\Big(\|b(s,0,0)\|+\|b_0(s)\|+\|g(s,0,0)\|+\|g_0(s)\|\Big)^2\1n
ds\Big]\Big\}.\ea$$
Note that the constant $K$ in front of $\e$ above is universal. Then
choose an $\e>0$ satisfying $K\e\le{1\over2}$ so that the first term
on the right hand side can be absorbed into the left hand, leading
to the following:
$$\ba{ll}
\ns\ds\|y^{\rho+\e}(\cd)\|_\infty+\|\psi^{\rho+\e}(\cd)\|_\infty\les
K\Big\{\|x\|+\|h(0)\|
+\|h_0\|\\
\ns\ds\qq\qq+\int_t^T\Big(\|b(s,0,0)\|+\|b_0(s)\|+\|g(s,0,0)\|+\|g_0(s)\|\)ds\Big\}.\ea$$
Continuing the above procedure, we obtain the solvability of the
following coupled FBEE:
$$\left\{\2n\ba{ll}
\ns\ds
\dot y(t)=Ay(t)+b(t,y(t),\psi(t))+b_0(t),\\
\ns\ds\dot\psi(t)=-A^*\psi(t)-g(t,y(t),\psi(t))-g_0(t),\\
\ns\ds y(0)=x,\qq\psi(T)=h(y(T))+h_0,\ea\right.$$
with the mild solution $(y(\cd),\psi(\cd))$ satisfying
$$\ba{ll}
\ns\ds\|y(\cd)\|_\infty+\|\psi(\cd)\|_\infty\les
K\Big\{\|x\|+\|h(0)\|+\int_0^T\Big(\|b(s,0,0)\|+\|b_0(s)\|+
\|g(s,0,0)\|+\|g_0(s)\|\Big)ds\Big\}.\ea$$
Thus, in particular, by taking $(b_0,g_0,h_0)=0$, we obtain the
solvability of FBEE (1.1) with estimate
\bel{5.15}\ba{ll}
\ns\ds\|y(\cd)\|_\infty\1n+\1n\|\psi(\cd)\|_\infty\1n\les\1n
K\Big\{\|x\|\1n+\1n\|h(0)\|\1n+\2n\int_0^T\2n\Big(\|b(s,0,0)\|\1n+\1n\|g(s,0,0)\|\Big)ds\Big\}.\ea\ee
Now, let $(\bar y(\cd),\bar\psi(\cd))$ be a mild solution to (1.1)
corresponding to $(\bar b,\bar g,\bar h)\in\cG_2\cap\cG_3$. Then,
$(\bar y(\cd)-y(\cd),\bar\psi(\cd)-\psi(\cd))$ satisfies a linear
FBEE with the generator admitting a type (I) or type (II) Lyapunov
operator $\Pi(\cd)$, the same as that for the generator $(b,g,h)$.
Hence, applying (\ref{5.15}), we obtain the following stability
estimate:
\bel{5.16}\ba{ll}
\ns\ds\|\bar
y(\cd)-y(\cd)\|_\infty+\|\bar\psi(\cd)-\psi(\cd)\|_\infty\les
K\Big\{\|\bar x-x\|+\|\bar h(\bar y(T))-h(\bar
y(T))\|\\
\ns\ds\qq\q+\int_0^T\(\|\bar b(s,\bar y(s),\bar\psi(s))-b(s,\bar
y(s),\bar\psi(s))\|+\|\bar g(s,\bar y(s),\bar\psi(s))-g(s,\bar
y(s),\bar\psi(s))\|\)ds\Big\}.\ea\ee
This proves the theorem. \endpf

\section{Construction of Lyapunov Operators and Solvable FBEEs}

In this section, we will construct some Lyapunov operators, through
which we obtain well-posedness of corresponding FBEEs. First of all,
we prove the following result which is practically more convenient
to use than the definition.

\ms

\bf Theorem 7.1. \sl Let {\rm(H0)} hold and let
$(b,g,h)\in\cG_2\cap\cG_3$. Let
$\Pi(\cd)\equiv\begin{pmatrix}\sc P(\cd)&\sc\G(\cd)^*\\
\sc\G(\cd)&\sc\bar P(\cd)\end{pmatrix}$ be a mild solution to linear
Lyapunov differential equation $(\ref{Pi})$ for some
\bel{MQ}\dbM(\cd)\equiv\begin{pmatrix}\sc M(\cd)&\sc0\\
\sc0&\sc-\bar
M(\cd)^*\end{pmatrix},\qq\dbQ(\cd)\equiv\begin{pmatrix}\sc
Q_0(\cd)&\Th(\cd)^*\\ \Th(\cd)&\sc\bar Q_0(\cd)\end{pmatrix}.\ee
Then $\Pi(\cd)$ is both a Lyapunov operator of types {\rm(I)} and
{\rm(II)} for $(b,g,h)$ if the following hold:
\bel{t,T}\bar P(t)\les-\d,\qq P(T)\ges\d,\ee
\bel{T}P(T)+h_y(y)^*\G(T)\1n+\1n\G(T)^*h_y(y)+h_y(y)^*\bar
P(T)h_y(y)\ges\d,\qq\forall y\in X,\ee
and
\bel{t<s<T}\left\{\2n\ba{ll}
\ns\ds\Pi(t)\dbM(t)+\dbM(t)^*\Pi(t)-\dbQ(t)\les-\d,\qq t\in[0,T],\\
[1mm]
\ns\ds\Pi(t)\dbB(t,y,\psi)+\dbB(t,y,\psi)^*\Pi(t)+\Pi(t)\dbM(t)+\dbM(t)^*
\Pi(t)-\dbQ(t)\les-\d,\q(t,y,\psi)\1n\in\1n[0,T]\1n\times \1n
X\1n\times\1n X.\ea\right.\ee
for some $\d>0$, with
$$\dbB(t,y,\psi)=\begin{pmatrix}\sc b_y(t,y,\psi)&\sc b_\psi(t,y,\psi)\\
\sc-g_y(t,y,\psi)&\sc-g_\psi(t,y,\psi)\end{pmatrix},\qq(t,y,\psi)\in[0,T]
\times X\times X.$$

\ms

\rm

Note that $\d>0$ appears in (\ref{t,T})--(\ref{t<s<T}) does not have
to be the same. But, we can always make them the same by shrinking
$\d$ if necessary.

\ms

\it Proof. \rm First of all, in order $\Pi(\cd)$ to be a type (I)
Lyapunov operator of the generator $(b,g,h)$, one needs (\ref{I}).
Hence, at $t=0$, one needs
$$0\ges\Pi(0)+\begin{pmatrix}\sc
-K&\sc0\cr\sc0&\sc0\end{pmatrix}=\begin{pmatrix}\sc
P(0)-K&\sc\G(0)^*\cr\sc\G(0)&\sc\bar P(0)\end{pmatrix},$$
for some $K>0$, which will be ensured by the following:
\bel{}\bar P(0)\les-\d,\ee
for some $\d>0$. Next, at $t=T$, one needs
\bel{6.5}\ba{ll}
\ns\ds0\les\begin{pmatrix}\sc I&\sc\rho h_y(y)^*\\
\sc0&\sc I\end{pmatrix}\Pi(T)\begin{pmatrix}\sc I&\sc0\\
\sc\rho h_y(y)&\sc I\end{pmatrix}+\begin{pmatrix}\sc-\m h_y(y)^*
h_y(y)&\sc0\cr\sc0&\sc K\end{pmatrix}\\ [4mm]
\ns\ds\q=\1n\rho^2\1n\begin{pmatrix}\sc0&\sc h_y(y)^*\\
\sc0&\sc0\end{pmatrix}\1n\Pi(T)\1n\begin{pmatrix}\sc0&\sc0\\
\sc h_y(y)&\sc0\end{pmatrix}\2n+\1n\rho\[\1n\begin{pmatrix}\sc0&\sc h_y(y)^*\\
\sc0&\sc0\end{pmatrix}\1n\Pi(T)\1n+\1n\Pi(T)\1n\begin{pmatrix}\sc0&\sc0\\
\sc h_y(y)&\sc0\end{pmatrix}\1n\]\\ [4mm]
\ns\ds\qq\q+\Pi(T)+\begin{pmatrix}\sc-\m h_y(y)^*
h_y(y)&\sc0\cr\sc0&\sc K\end{pmatrix},\qq\qq\forall y\1n\in\1n
X,~\rho\1n\in\1n[0,1],\ea\ee
for some $\m,K>0$. If we are able to show the following (which will
be done below)
\bel{bar P(T)<0}\bar P(T)\les0,\ee
then
\bel{bar P(T)<0*}\begin{pmatrix}0&h_y(y)^*\\
0&0\end{pmatrix}\Pi(T)\begin{pmatrix}0&0\\
h_y(y)&0\end{pmatrix}=\begin{pmatrix}h_y(y)^*\bar P(T)h_y(y)&0\\
0&0\end{pmatrix}\les0.\ee
Hence, (\ref{6.5}) is true if and only if it is true for $\rho=0,1$,
i.e.,
$$0\les\Pi(T)+\begin{pmatrix}-\m h_y(y)^*
h_y(y)&0\cr0&K\end{pmatrix}=\begin{pmatrix}P(T)-\m
h_y(y)^*h_y(y)&\G(T)^*\cr\G(T)&\bar P(T)+K\end{pmatrix},$$
and
$$\ba{ll}
\ns\ds0\les\begin{pmatrix}0&h_y(y)^*\\
0&0\end{pmatrix}\Pi(T)\begin{pmatrix}0&0\\
h_y(y)&0\end{pmatrix}+\[\begin{pmatrix}0&h_y(y)^*\\
0&0\end{pmatrix}\Pi(T)+\Pi(T)\begin{pmatrix}0&0\\
h_y(y)&0\end{pmatrix}\]\\ [4mm]
\ns\ds\qq+\Pi(T)+\begin{pmatrix}-\m h_y(y)^* h_y(y)&0\cr0&
K\end{pmatrix}\\ [4mm]
\ns\ds\q=\2n\begin{pmatrix}P(T)+h_y(y)^*\G(T)+\G(T)^* h_y(y)+\1n
h_y(y)^*\bar P(T)h_y(y)-\m h_y(y)^*h_y(y)&\G(T)^*\1n+h_y(y)^*\bar
P(T)\cr\G(T)+\bar P(T)h_y(y)&\bar P(T)+K\end{pmatrix}.\ea$$
By choosing $\m>0$ small enough and $K>0$ large enough, we see that
the above is implied by the following:
$$\left\{\2n\ba{ll}
\ns\ds\begin{pmatrix}P(T)&\G(T)^*\cr\G(T)&\bar
P(T)+K\end{pmatrix}\ges\d,\\ [4mm]
\ns\ds\begin{pmatrix}P(T)+h_y(y)^*\G(T)+\G(T)^* h_y(y)+\1n
h_y(y)^*\bar P(T)h_y(y)&\G(T)^*+h_y(y)^*\bar P(T)\cr\G(T)+\bar
P(T)h_y(y)&\bar P(T)+K\end{pmatrix}\ges\d,\ea\right.$$
for some $\d>0$. This will further be implied by
$$\ba{ll}
\ns\ds P(T)\ges\d,\q P(T)\1n+\1n
h_y(y)^*\G(T)\1n+\1n\G(T)^*h_y(y)\1n+\1n h_y(y)^*\bar
P(T)h_y(y)\1n\ges\1n\d,\ea$$
for some $\d>0$. Hence, to summarize, at $t=0,T$, it suffices to
have (\ref{t,T})--(\ref{T}).

\ms

Now, we look at $t\in(0,T)$. One needs
$$\begin{pmatrix}\dbH^\rho(t,y,\psi)-\dbQ(t)+\m\begin{pmatrix}
g_y(t,y,\psi)^*g_y(t,y,\psi)&0\cr0&0\end{pmatrix}&\Pi(t)\cr
\Pi(t)&-K\end{pmatrix}\2n\les\1n0,$$
for some $\m,K>0$. The left-hand side is affine in $\rho$. Hence,
the above is true for all $\rho\in[0,1]$ if and only if it is true
for $\rho=0,1$, i.e.,
$$\3n\ba{ll}
\ns\ds\begin{pmatrix}\Pi(t)\dbM(t)+\dbM(t)^*\Pi(t)-\dbQ(t)
+\m\begin{pmatrix}
g_y(t,y,\psi)^*g_y(t,y,\psi)&0\cr0&0\end{pmatrix}&\Pi(t)\cr
\Pi(t)&-K\end{pmatrix}\2n\les\1n0,\\ [2mm]
\ns\ds\begin{pmatrix}\Pi(t)\dbB(t,y,\psi)\1n+\1n\dbB(t,y,\psi)^*\1n\Pi(t)
\1n+\1n\Pi(t)\dbM(t)\1n+\1n\dbM(t)^*\1n\Pi(t)-\dbQ(t)
\1n+\1n\m\1n\begin{pmatrix}g_y(t,y,\psi)^*\1n
g_y(t,y,\psi)&0\cr0&0\end{pmatrix}&\Pi(t)\cr
\Pi(t)&-K\end{pmatrix}\3n\les\1n0,\ea$$
for some $\m,K>0$. These are implied by the following: (by letting
$\m>0$ small enough)
\bel{6.7}\begin{pmatrix}\Pi(t)\dbM(t)+\dbM(t)^*\Pi(t)-\dbQ(t)
&\Pi(t)\cr\Pi(t)&-K\end{pmatrix}\2n\les-\d,\ee
and
\bel{6.8}\begin{pmatrix}\Pi(t)\dbB(t,y,\psi)+\dbB(t,y,\psi)^*\Pi(t)
+\Pi(t)\dbM(t)+\dbM(t)^*\Pi(t)-\dbQ(t)&\Pi(t)\cr
\Pi(t)&-K\end{pmatrix}\2n\les-\d,\ee
for some $\d>0$. It is clear that (\ref{6.7})--(\ref{6.8}) hold for
some large $K>0$ if (\ref{t<s<T}) holds. Further, we note that the
first condition in (\ref{t<s<T}) implies
$$\left\{\2n\ba{ll}
\ns\ds-P(t)M(t)-M(t)^*P(t)+Q_0(t)\ges\d,\\ [2mm]
\ns\ds\bar P(t)\bar M(t)^*+\bar M(t)\bar P(t)+\bar
Q_0(t)\ges\d.\ea\right.$$
Hence, for any $t\in[t,T]$, using $P(T)\ges\d$ and $\bar
P(0)\les-\d$, we obtain
\bel{L1a}\ba{ll}
\ns\ds P(t)=e^{A^*(T-t)}P(T)e^{A(T-t)}+\1n\int_t^T\2n
e^{A^*(s-t)}[-P(s)M(s)\1n-\2n M(s)^*P(s)\1n+\1n
Q_0(s)]e^{A(s-t)}ds\ges0,\ea\ee
and
\bel{L2a}\ba{ll}
\ns\ds\bar P(t)=e^{At}\bar P(0)e^{A^*t}-\1n\int_0^t\2n
e^{A(t-s)}[\bar P(s)\bar M(s)^*\2n+\1n\bar M(s)\bar P(s)\1n+\1n\bar
Q_0(s)]e^{A^*(t-s)}ds\les0.\ea\ee
In particular, (\ref{bar P(T)<0}) holds. This proves that under
(\ref{t,T})--(\ref{t<s<T}), $\Pi(\cd)$ is a type (I) Lyapunov
operator for the generator $(b,g,h)$.

\ms

We can similarly prove that under (\ref{t,T})--(\ref{t<s<T}),
$\Pi(\cd)$ is also a type (II) Lyapunov operator for the generator
$(b,g,h)$. \endpf

\ms

\bf Corollary 7.2. \sl Let {\rm(H0)} hold and let
$(b,g,h)\in\cG_2\cap\cG_3$.
Let $\Pi(\cd)\equiv\begin{pmatrix}\sc P(\cd)&\sc\G(\cd)^*\\
\sc\G(\cd)&\sc\bar P(\cd)\end{pmatrix}$ be a mild solution to
$(\ref{Pi})$ for some $\dbM(\cd)$ and $\dbQ(\cd)$ of form
{\rm(\ref{MQ})} such that {\rm(\ref{t,T})--(\ref{T})} hold. Suppose
the following hold:
\bel{t<s<T*}\left\{\2n\ba{ll}
\ns\ds\Pi(t)\dbM(t)+\dbM(t)^*\Pi(t)-\dbQ(t)\les-\d-\e,\qq t\in[0,T],\\
[1mm]
\ns\ds\Pi(t)\dbB(t,y,\psi)+\dbB(t,y,\psi)^*\Pi(t)\les\e,\qq\q(t,y,\psi)\1n
\in\1n(0,T)\1n\times \1n X\1n\times\1n X,\ea\right.\ee
for some $\d,\e>0$. Then $\Pi(\cd)$ is both a Lyapunov operator of
types {\rm(I)} and {\rm(II)} for $(b,g,h)$. In particular, this is
the case if the following holds:
\bel{t<s<T**}\left\{\2n\ba{ll}
\ns\ds\dbQ(t)\ges\d+\e,\qq t\in[0,T],\\
\ns\ds\Pi(t)\dbB(t,y,\psi)\1n+\1n\dbB(t,y,\psi)^*\Pi(t)\1n\les\e,
\q(t,y,\psi)\1n\in\1n(0,T)\1n\times\1n X\1n\times\1n X,\ea\right.\ee
for some $\d,\e>0$.

\ms

\it Proof. \rm It is clear that (\ref{t<s<T*}) implies
(\ref{t<s<T}). In particular, by letting $\dbM(\cd)=0$, we see that
(\ref{t<s<T**}) implies (\ref{t<s<T}). \endpf

\ms

We note that condition (\ref{t<s<T*}) is more convenient to check
than (\ref{t<s<T}). Next, we look at some concrete special cases of
Theorem 7.1, which will be more practically useful. We first present
the following result.

\ms

\bf Lemma 7.3. \sl Let {\rm(H0)} hold and $p_1,\bar p_0,q_0,\bar
q_0,\th,\g\in\dbR$ and $m,\bar m>-\si_0$. Let
\bel{MM}\left\{\2n\ba{ll}
\ns\ds M(t)=mI,\q\bar M(t)=\bar mI,\\
\ns\ds Q_0(t)=q_0I,\q\bar Q_0(t)=\bar q_0I,\q\Th(t)=\th
I.\ea\right.\ee
The mild solution $\Pi(\cd)$ of {\rm(\ref{Pi})} satisfying
\bel{7.16}P(T)=p_1I,\q\bar P(0)=-\bar p_0I,\q\G(T)=\g I,\ee
is given by the following: In the case $A^*=A$, for any $t\in[0,T]$,
\bel{Pi(A=A)}\ba{ll}
\ns\ds\Pi(t)\1n=\1n\begin{pmatrix}\sc p_1e^{2(A-m)(T-t)}&\sc\big[\g
e^{(\bar m-m)(T-t)} +\th{e^{(\bar m-m)(T-t)}-1\over\bar m-m}\big]I\\
\sc\big[\g e^{(\bar m-m)(T\1n-\1n s)}\1n+\th{e^{(\bar
m-m)(T-t)}-1\over
\bar m-m}\big]I&\sc-\bar p_0e^{2(A-\bar m)t}\end{pmatrix}\\
[4mm]
\ns\ds\qq\qq+{1\over2}\1n\begin{pmatrix}\sc q_0(A-m)^{-1}[e^{2(A-m)(T-t)}-I]&\sc0\\
\sc0&\sc-\bar q_0(A-\bar m)^{-1}[e^{2(A-\bar
m)t}-I]\end{pmatrix}\1n,\ea\ee
and in the case $A^*=-A$, for $t\in[0,T]$,
\bel{Pi(A=-A)}\Pi(t)\1n=\2n\begin{pmatrix}\sc
\big[p_1e^{-2m(T-t)}+{q_0\over2m}(1-e^{-2m(T-t)})\big]I&\sc\big[\g
e^{(\bar m-m)(T\1n-\1n s)}\1n+\th{e^{(\bar m-m)(T-t)}-1\over\bar
m-m}\big]I\\ \\ \sc\big[\g e^{(\bar m-m)(T-t)}+\th{e^{(\bar
m-m)(T-t)}-1\over\bar m-m}\big]I&\sc-\big[\bar p_0e^{-2\bar
mt}+{\bar q_0\over2\bar m}(1-e^{-2\bar
mt})\big]I\end{pmatrix}\2n.\ee
In the above, the following convention is adopted:
\bel{convention}{1-e^{-\a\b}\over\a}\equiv\b,\qq\hb{if }\a=0.\ee
In particular, if
\bel{=}\bar m=m,\ee
then in the case $A^*=A$, for $m\in\dbR(A)$ and $t\in[0,T]$,
\bel{Pi(A=A)**}\ba{ll}
\ns\ds\Pi(t)\1n=\begin{pmatrix}\sc p_1e^{2(A-m)(T-t)}&\sc
[\g+\th(T-t)]I\\
\sc[\g+\th(T-t)]I&\sc-\bar p_0e^{2(A-m)t}\end{pmatrix}\\
[4mm]
\ns\ds\qq\qq+{1\over2}\1n\begin{pmatrix}\sc
q_0(A-m)^{-1}[e^{2(A-m)(T-t)}-I]&
\sc0\\
\sc0&\sc-\bar q_0(A-m)^{-1}[e^{2(A-m)t}-I]\end{pmatrix}\1n,\ea\ee
and in the case $A^*=-A$, for $t\in[0,T]$,
\bel{Pi(A=-A)**}\ba{ll}
\ns\ds\Pi(s)\1n=\begin{pmatrix}\[p_1e^{-2m(T-s)}+q_0{1-e^{2m(T-s)}\over2m}\]I&\sc[\g+\th(T-t)]I\\
\sc[\g+\th(T-t)]I&\sc-\[\bar p_0e^{-2m(s-t)}+\bar
q_0{1-e^{-2mt}\over2m}\]I\end{pmatrix}.\ea\ee
Further, if $m=\bar m=0$, then, for $A^*=A$, $t\in[0,T]$,
\bel{Pi(A=A)*}\Pi(t)=\begin{pmatrix}\sc
(p_1+{q_0\over2}A^{-1})e^{2A(T-t)}+{q_0\over2}A^{-1}&\sc[\g
+\th(T-t)]I\\
\sc[\g+\th(T-t)]I&\sc-(\bar p_0+{\bar q_0\over2}A^{-1})e^{2At}+{\bar
q_0\over2}A^{-1}
\end{pmatrix},\ee
and for $A^*=-A$, $t\in[0,T]$,
\bel{Pi(A=-A)*}\Pi(t)=\begin{pmatrix}\sc[p_1+q_0(T-t)]I&\sc[\g
+\th(T-t)]I\\
\\ \sc[\g+\th(T-t)]I&\sc-[\bar p_0+\bar q_0t]I\end{pmatrix}.\ee

\it Proof. \rm With the choice of (\ref{MM}), we have
\bel{6.23P}\ba{ll}
\ns\ds P(t)=e^{(A-m)^*(T-t)}P(T)e^{(A-m)(T-t)}+\int_t^Te^{(A-m)^*(s-t)}Q_0(s)e^{(A-m)(s-t)}ds\\
\ns\ds\qq=p_1e^{(A-m)^*(T-t)}e^{(A-m)(T-t)}+q_0\int_t^Te^{(A-m)^*(s-t)}e^{(A-m)(s-t)}ds,\q
t\in[0,T],\ea\ee
and
\bel{6.24barP}\ba{ll}
\ns\ds\bar P(t)=e^{(A-\bar m)t}\bar P(0)e^{(A-\bar
m)^*t}-\int_0^t\2n
e^{(A-\bar m)(t-s)}\bar Q_0(s)e^{(A-\bar m)^*(t-s)}ds\\
\ns\ds\qq=-\bar p_0e^{(A-\bar m)t}e^{(A-\bar m)^*t}-\bar
q_0\int_0^te^{(A-\bar m)(t-s)}e^{(A-\bar m)^*(t-s)}ds,\q
t\in[0,T].\ea\ee
Also, in the current case, $\G(\cd)$ satisfies
$$\dot\G(t)+(\bar m-m)\G(t)+\th I=0,$$
which leads to (with $\G(T)=\g I$)
\bel{6.25G}\ba{ll}
\ns\ds\G(t)=\[\g e^{(\bar m-m)(T-t)}+\th\int_t^Te^{(\bar
m-m)(s-t)}dr\]I\\
\ns\ds\qq=\[\g e^{(\bar m-m)(T-t)}+\th{e^{(\bar m-m)(T-t)}-1\over
\bar m-m}\]I=\G(t)^*,\q t\in[0,T].\ea\ee
When $\bar m-m=0$, the above is understood as
\bel{6.26G}\G(t)=\big[\g+\th(T-t)\big]I,\qq t\in[0,T].\ee
We now look at two cases.

\ms

In the case $A^*=A$, (\ref{6.23P}) and (\ref{6.24barP}) become
\bel{6.27P}P(t)=p_1e^{2(A-m)(T-t)}+{q_0\over2}(A-m)^{-1}\big[e^{2(A-m)(T-t)}
-I\big],\q t\in[0,T],\ee
and
\bel{6.28barP}\bar P(t)=-\bar p_0e^{2(A-\bar m)t}-{\bar
q_0\over2}(A-\bar m)^{-1}\big[e^{2(A-\bar m)t}-I\big],\q
t\in[0,T].\ee

\ms

In the case $A^*=-A$, (\ref{6.23P}) and (\ref{6.24barP}) become
\bel{6.29P}\ba{ll}
\ns\ds P(t)\1n=\1n\(p_1e^{-2m(T-t)}\1n+\1n
q_0\int_t^Te^{-2m(s-t)}ds\)I\1n
=\1n\(p_1e^{-2m(T-t)}\1n+\1n{q_0\over2m}(1-e^{-2m(T-t)})\)I,\q
t\in[0,T],\ea\ee
with the above understood as follows when $m=0$,
\bel{6.30P}P(t)=\(p_1+q_0(T-t)\)I,\qq t\in[0,T],\ee
and
\bel{6.31barP}\ba{ll}
\ns\ds\bar P(t)=-\(\bar p_0e^{-2\bar mt}+\bar q_0\int_0^te^{-2\bar
m(t-s)}ds\)I=-\(\bar p_0e^{-2\bar mt}+{\bar q_0\over2\bar
m}(1-e^{-2\bar mt})\)I,\qq t\in[0,T],\ea\ee
with the above understood as follows when $\bar m=0$,
\bel{6.32barP}\bar P(t)=-\(\bar p_0+\bar q_0t\)I,\qq t\in[0,T].\ee
The rest conclusions are clear. \endpf

\ms

Combining Theorem 7.1 or Corollary 7.2 with Lemma 7.3, we can
present many concrete cases for which the corresponding FBEEs are
well-posed. For the simplicity of presentation, we only consider
below the case that (\ref{=}) holds. First, we present a simple
lemma.

\ms

\bf Lemma 7.4. \rm Let
$$f(\k)=\a e^{-\k}+\b{1-e^{-\k}\over\k},\qq\k>0,$$
with $\a,\b>0$. Then $\k\mapsto f(\k)$ is decreasing on $(0,\infty)$
and
$$0=\lim_{\k\to\infty}f(\k)=\inf_{\k>0}f(\k)<\sup_{\k>0}f(\k)=\lim_{\k\to0}f(\k)=\a+\b.$$

\it Proof. \rm We note that
$$f'(\k)=-\a e^{-\k}+\b{e^{-\k}\k-(1-e^{-k})\over\k}=-\({\a\over e^\k}+{e^\k-1-\k\over\k^2e^\k}\)<0.$$
Then our conclusion follows immediately. \endpf

\ms

\bf Theorem 7.5. \sl Let {\rm(H0)} hold and let
$(b,g,h)\in\cG_2\cap\cG_3$. Suppose there are constants $p_1,\bar
p_0,q_0,\bar q_0,\d,\bar\d,\e>0$, $m>-\si_0$, and $\g,\th\in\dbR$
such that
\bel{6.33}\ba{ll}
\ns\ds p_1I+\g[h_y(y)+h_y(y)^*]-\[\bar p_0e^{-2(\si_0+m)(T-t)}+{\bar
q_0(1-e^{-2(\si_0+m)(T-t)})
\over2(\si_0+m)}\]h_y(y)^*h_y(y)\ges\d,\ea\ee
\bel{6.34}\ba{ll}
\ns\ds\begin{pmatrix}q_0{\si_0+me^{-2(\si_0+m)(T-t)}\over\si_0+m}-2mp_1e^{-2(\si_0+m)(T-t)}
&\th\\
\th&\bar q_0{\si_0+me^{-2(\si_0+m)t}\over\si_0+m}-2m\bar
p_0e^{-2(\si_0+m)t}
\end{pmatrix}\ges\bar\d+\e,\\ [5mm]
\ns\ds\qq\qq\qq\qq\qq\qq\qq\qq\qq\qq\qq\qq\qq\qq\qq\qq\forall
t\in[0,T].\ea\ee
Then the corresponding FBEE is well-posed if one of the following
holds:

\ms

{\rm(i)} In the case $A^*=A$, it holds that
\bel{PiB+BPi}\ba{ll}
\ns\ds\begin{pmatrix}p_1e^{2(A-m)(T-t)}&0\\
0&-\bar p_0e^{2(A-m)t}\end{pmatrix}\dbB(t,y,\psi)
+\dbB(t,y,\psi)^*\begin{pmatrix}p_1e^{2(A-m)(T-t)}&0\\
0&-\bar p_0e^{2(A-m)t}\end{pmatrix}\\
\ns\ds+[\g+\th(T-t)]\[\begin{pmatrix}0&I\\
I&0\end{pmatrix}\dbB(t,y,\psi)
+\dbB(t,y,\psi)^*\begin{pmatrix}0&I\\ I&0\end{pmatrix}\]\\
[4mm]
\ns\ds+{1\over2}\[\begin{pmatrix}q_0(A-m)^{-1}[e^{2(A-m)(T-t)}-I]&0\\
0&-\bar q_0(A-m)^{-1}[e^{2(A-m)t}-I]\end{pmatrix}\dbB(t,y,\psi)\\
\ns\ds+\dbB(t,y,\psi)^*\2n\begin{pmatrix}
q_0(A-m)^{-1}[e^{2(A-m)(T-t)}-I] &0\\ 0&-\bar
q_0(A-m)^{-1}[e^{2(A-m)t}-I]\end{pmatrix}\1n\]\2n\les\1n\e.\ea\ee

\ms

{\rm(ii)} In the case $A^*=-A$, it holds that
\bel{PiB+BPi*}\ba{ll}
\ns\ds\begin{pmatrix}p_1e^{-2m(T-t)}&0\\
0&-\bar p_0e^{-2mt}\end{pmatrix}\dbB(t,y,\psi)+\dbB(t,y,\psi)^*\begin{pmatrix}p_1e^{-2m(T-t)}&0\\
0&-\bar p_0e^{-2mt}\end{pmatrix}\\
\ns\ds+[\g+\th(T-t)]\[\begin{pmatrix}0&I\\
I&0\end{pmatrix}\dbB(t,y,\psi)+\dbB(t,y,\psi)^*\begin{pmatrix}0&
I\\ I&\sc0\end{pmatrix}\]\\
\ns\ds+\begin{pmatrix}q_0{1-e^{-2m(T-t)}\over2m}I&0\\
0&-\bar
q_0{1-e^{-2mt}\over2m}I\end{pmatrix}\dbB(t,y,\psi)+\dbB(t,y,\psi)^*\begin{pmatrix}
q_0{1-e^{-2m(T-t)}\over2m}I&
0\\
0&-\bar q_0{1-e^{-2mt}\over2m}I\end{pmatrix}\les\e.\ea\ee

\ms

\it Proof. \rm (i) We consider the case $A^*=A$. First of all, by
taking
\bel{6.38}P(T)=p_1I>0,\qq\bar P(0)=-\bar p_0I<0,\ee
we see that (\ref{t,T}) holds. Next, to get (\ref{T}), we look at
the following (recalling (\ref{6.31barP})):
$$\ba{ll}
\ns\ds P(T)+h_y(y)^*\G(T)\1n+\1n\G(T)^*h_y(y)+h_y(y)^*\bar
P(T)h_y(y)\\
\ns\ds=p_1I+\g[h_y(y)+h_y(y)^*]-h_y(y)^*\(\bar p_0e^{2(A-m)T}+{\bar
q_0\over2}(A-m)^{-1}\big[e^{2(A-m)T}-I\big]\)h_y(y).\ea$$
Let us estimate the quadratic term in $h_y(y)$ on the right hand
side of the above. To this end, we observe the following: For any
$\t>\si_0$ (recall $m>-\si_0$),
$$\ba{ll}
\ns\ds \bar p_0e^{-2(\t+m)T}+{\bar
q_0T(1-e^{-2(\t+m)T})\over2(\t+m)T}\equiv\bar p_0e^{-\k}+{\bar
q_0T(1-e^{-\k})\over\k}\equiv f(\k),\ea$$
with $\k=2(\t+m)T\ges2(\si_0+m)T>0$. By Lemma 7.4, we have
$$\ba{ll}
\ns\ds\sup_{\k\ges2(\si_0+m)T}f(\k)=f(2(\si_0+m)T)=\bar
p_0e^{-2(\si_0+m)T}+{\bar
q_0(1-e^{-2(\si_0+m)T})\over2(\si_0+m)}.\ea$$
By the spectral decomposition of $A$, making use of (\ref{6.33}),
one has
$$\ba{ll}
\ns\ds P(T)+h_y(y)^*\G(T)\1n+\1n\G(T)^*h_y(y)+h_y(y)^*\bar
P(T)h_y(y)\\
\ns\ds=p_1I+\g[h_y(y)+h_y(y)^*]-h_y(y)^*\(\bar p_0e^{2(A-m)T}+{\bar
q_0\over2}(A-m)^{-1}
\big[e^{2(A-m)T}-I\big]\)h_y(y)\\
\ns\ds\ges p_1I+\g[h_y(y)+h_y(y)^*]-\[\bar p_0e^{-2(\si_0+m)T}+{\bar
q_0(1-e^{-2(\si_0+m)T}) \over2(\si_0+m)}\]h_y(y)^*h_y(y)\ges\d,\ea$$
which gives (\ref{T}). Next, note that (in the case $\bar m=m$)
\bel{6.40}-\Pi(t)\dbM(t)-\dbM(t)^*\Pi(t)+\dbQ(t)=\begin{pmatrix}\sc-2mP(t)
+q_0I&\sc\th I\\
\sc\th I& \sc2m\bar P(t)+\bar q_0I\end{pmatrix}.\ee
We now look at the following (noting (\ref{6.27P}) and
(\ref{6.28barP})):
$$\ba{ll}
\ns\ds-2mP(t)+q_0I=\1n-2mp_1e^{2(A-m)(T-t)}\1n-\1n mq_0(A\1n-\1n
m)^{-1}\big[e^{2(A-m)(T-t)}\1n-\1n I\big]
\2n+\1n q_0I\\
\ns\ds=-2mp_1e^{2(A-m)(T-t)}+q_0\big[A-me^{2(A-m)(T-t)}\big](A-m)^{-1},\ea$$
and
$$\ba{ll}
\ns\ds2m\bar P(t)+q_0I=-2m\bar p_0e^{2(A-m)t}
-m\bar q_0(A-m)^{-1}\big[e^{2(A-m)t}-I\big]\2n+\1n\bar q_0I\\
\ns\ds=-2m\bar p_0e^{2(A-m)t}+\bar
q_0\big[A-me^{2(A-m)t}\big](A-m)^{-1}.\ea$$
Similar to the above, for any $\t>\si_0$ and $m>-\si_0$, we have
($\t+m>0$)
$$\ba{ll}
\ns\ds-2mp_1e^{-2(\t+m)(T-t)}+{q_0(\t+me^{-2(\t+m)(T-t)})\over\t+m}\\
\ns\ds=-2mp_1e^{-2(\t+m)(T-t)}+q_0+{2mq_0(T-t)(e^{-2(\t+m)(T-t)}-1)
\over2(\t+m)(T-t)}\\
\ns\ds=q_0-\[2mp_1e^{-\k}-{2mq_0(T-t)(1-e^{-\k})\over\k}\]\equiv
q_0- f(\k),\ea$$
with $\k=2(\t+m)(T-t)>2(\si_0+m)(T-t)$. By Lemma 7.4 again, we have
\bel{6.41}\ba{ll}
\ns\ds-2mp_1e^{-2(\t+m)(T-t)}+{q_0(\t+me^{-2(\t+m)(T-t)})\over\t+m}\\
\ns\ds\qq\ges q_0-\sup_{\k\ges2(\si_0+m)(T-t)}f(\k)=q_0-f(2(\si_0+m)(T-t))\\
\ns\ds\qq=q_0-2mp_1e^{-2(\si_0+m)(T-t)}+{2mq_0(e^{-2(\si_0+m)(T-t)}-1)\over
2(\si_0+m)}\\
\ns\ds\qq=q_0\[1-{2m(1-e^{-2(\si_0+m)(T-t)})\over2(\si_0+m)}\]
-2mp_1e^{-2(\si_0+m)(T-t)}\\
\ns\ds\qq=q_0{\si_0+me^{-2(\si_0+m)(T-t)}\over\si_0+m}-2mp_1e^{-2(\si_0+m)(T-t)}\\
\ns\ds\qq\ges
q_0-\lim_{\k\to0}f(\k)=q_0[1+2\si_0(T-t)]+2\si_0p_1.\ea\ee
Similarly,
\bel{6.42}\ba{ll}
\ns\ds-2m\bar p_0e^{-2(\t+m)t}+{\bar
q_0(\t+me^{-2(\t+m)t})\over\t+m}\ges\bar
q_0\[1-{2m(1-e^{-2(\si_0+m)t})\over2(\si_0+m)}\]-2m\bar p_0
e^{-2(\si_0+m)t}\\
\ns\ds\qq\ges\bar q_0{\si_0+me^{-2(\si_0+m)t}\over\si_0+m}-2m\bar
p_0 e^{-2(\si_0+m)t}\ges\bar q_0[(1+2\si_0t]+2\si_0\bar p_0.\ea\ee
Consequently, using the spectral decomposition of $A$, we have
$$\ba{ll}
\ns\ds-2mP(t)+q_0I\ges q_0{\si_0+me^{-2(\si_0+m)(T-t)}\over\si_0+m}
-2mp_1e^{-2(\si_0+m)(T-t)},\\
[2mm]
\ns\ds2m\bar P(t)+\bar q_0I\ges\bar
q_0{\si_0+me^{-2(\si_0+m)t}\over\si_0+m}-2m\bar p_0
e^{-2(\si_0+m)t}.\ea$$
Hence,
$$\ba{ll}
\ns\ds-\Pi(t)\dbM(t)-\dbM(t)^*\Pi(t)+\dbQ(t)=\begin{pmatrix}-2mp_1e^{2(A-m)(T-t)}&\th I\\
\th I&-2m\bar p_0e^{2(A-m)t}\end{pmatrix}\\ [4mm]
\ns\ds\qq+\begin{pmatrix}q_0[A-me^{2(A-m)(T-t)}](A-m)^{-1}&0\\
0&\bar q_0[A-me^{2(A-m)t}](A-m)^{-1}\end{pmatrix}\\ [4mm]
\ns\ds\ges\begin{pmatrix}\big[q_0{\si_0+me^{-2(\si_0+m)(T-t)}\over\si_0+m}\big]I&\th I\\
\th I&\big[\bar q_0{\si_0+me^{-2(\si_0+m) t}\over\si_0+m}\big]I
\end{pmatrix}\\ [4mm]
\ns\ds\qq-2m\begin{pmatrix}p_1e^{-2(\si_0+m)(T-t)}I&0\\
0&\bar p_0e^{-2(\si_0+m)(T-t)}I\end{pmatrix}\ges\d,\qq\forall
t\in[0,T],\ea$$
for some $\d>0$, provided (\ref{6.34}) holds. Then by
(\ref{PiB+BPi}), together with the representation of $\Pi(\cd)$ from
Lemma 7.3, we have
$$\Pi(t)\dbB(t,y,\psi)+\dbB(t,y,\psi)^*\Pi(t)\les\e.$$
Hence, Corollary 7.2 applies.

\ms

(ii) We now consider the case $A^*=-A$. Again, we still have
(\ref{t,T}) by (\ref{6.33}). Next, for the current case, recalling
(\ref{6.28barP}),
$$\ba{ll}
\ns\ds P(T)+h_y(y)^*\G(T)\1n+\1n\G(T)^*h_y(y)+h_y(y)^*\bar
P(T)h_y(y)\\
\ns\ds=p_1I+\g[h_y(y)+h_y(y)^*]-\(\bar p_0e^{-2m(T-t)}+{\bar
q_0\over2m}(1-e^{-2m(T-t)})\)h_y(y)^*h_y(y)\ges\d.\ea$$
This leads to (\ref{6.33}) with $\si_0=0$. In the current case, we
still have (\ref{6.40}). We observe the following (noting
(\ref{6.29P}) and (\ref{6.31barP})):
$$\ba{ll}
\ns\ds-2mP(t)+q_0I=\1n-2m\(p_1e^{-2m(T-t)}+{q_0\over2m}(1-e^{-2m(T-t)})\)I
+\1n q_0I\\
\ns\ds=-2mp_1e^{-2m(T-t)}+q_0e^{-2m(T-t)}=(q_0-2mp_1)e^{-2m(T-t)},\ea$$
and
$$\ba{ll}
\ns\ds2m\bar P(t)+\bar q_0I=-2m\(\bar p_0e^{-2mt}+{\bar
q_0\over2m}(1-e^{-2mt})\)I
+\bar q_0I\\
\ns\ds=-2m\bar p_0e^{-2mt}+\bar q_0e^{-2mt}=(\bar q_0-2m\bar
p_0)e^{-2mt},\ea$$
Hence,
$$\ba{ll}
\ns\ds-\Pi(t)\dbM(t)-\dbM(t)^*\Pi(t)+\dbQ(t)\\
\ns\ds=\begin{pmatrix}(q_0-2mp_1)e^{-2m(T-t)}I&\th I\\
\th I&(\bar q_0-2m\bar p_0)e^{-2mt}I\end{pmatrix}\ges\d,\q\forall
t\in[0,T],\ea$$
for some $\d>0$, provided
$$\begin{pmatrix}(q_0-2mp_1)e^{-2m(T-t)}&\th\\
\th&(\bar q_0-2m\bar p_0)e^{-2mt}\end{pmatrix}>0,\q\forall
t\in[0,T],$$
which is implied by (\ref{6.34}) with $\si_0=0$. The rest proof is
obvious. \endpf

\ms

\bf Corollary 7.6. \sl Let {\rm(H0)} hold, and let
$(b,g,h)\in\cG_2\cap\cG_3$. Let $p_1,\bar p_0,q_0,\bar
q_0,\d,\bar\d,\e>0$, $\g\in\dbR$ such that
\bel{6.41}p_1+\g[h_y(y)+h_y(y)^*]-(\bar p_0+\bar
q_0T)h_y(y)^*h_y(y)\ges\d,\ee
and
\bel{6.42}\begin{pmatrix}q_0[1+2\si_0(T-t)]+2\si_0p_1&\th\\
\th&\bar q_0[1+2\si_0t]+2\si_0\bar p_0\end{pmatrix}\ges\bar\d+\e,\qq
t\in[0,T].\ee
Then the corresponding FBEE is well-posed if one of the following
holds:

\ms

{\rm(i)} For the case $A^*=A$, the following holds:
\bel{6.43}\ba{ll}
\ns\ds\begin{pmatrix}p_1e^{2(A+\si_0)(T-t)}&
0\\
0&-\bar p_0e^{2(A+\si_0)t}\end{pmatrix}\dbB(t,y,\psi)
+\dbB(t,y,\psi)^*\begin{pmatrix}p_1e^{2(A+\si_0)(T-t)}&0\\
0&-\bar p_0e^{2(A+\si_0)t}\end{pmatrix}\\
\ns\ds+[\g+\th(T-t)]\[\begin{pmatrix}0&I\\
I&0\end{pmatrix}\dbB(t,y,\psi)
+\dbB(t,y,\psi)^*\begin{pmatrix}0&I\\ I&0\end{pmatrix}\]\\
[4mm]
\ns\ds+\[\begin{pmatrix}q_0(T-t)\eta\big((A+\si_0)(T-t)\big)&0\\
0&-\bar
q_0t\eta\big((A+\si_0)t\big)\end{pmatrix}\dbB(t,y,\psi)\\
\ns\ds+\dbB(t,y,\psi)^*\2n\begin{pmatrix}\sc q_0
(T-t)\eta\big((A+\si_0)(T-t)\big)&0\\
0&-\bar
q_0t\eta\big((A+\si_0)t\big)\end{pmatrix}\1n\]\2n\les\1n\e,\ea\ee
where
\bel{eta}\eta(\k)=\left\{\2n\ba{ll}
\ns\ds{e^\k-1\over\k},\qq\k\ne0,\\
\ns\ds1,\qq\qq\k=0.\ea\right.\ee

{\rm(ii)} For the case $A^*=-A$, the following holds:
\bel{6.44}\ba{ll}
\ns\ds\begin{pmatrix}\sc[p_1+q_0(T-t)]I&\sc[\g+\th(T-t)]I\\
\sc[\g+\th(T-t)]I&\sc-[\bar p_0+\bar q_0t]I\end{pmatrix}\dbB(t,y,\psi)\\
[4mm]
\ns\ds+\dbB(t,y,\psi)^*\begin{pmatrix}\sc[p_1+q_0(T-t)]I&\sc[\g+\th(T-t)]I\\
\sc[\g+\th(T-t)]I&\sc-[\bar p_0+\bar
q_0t]I\end{pmatrix}\les\1n\e.\ea\ee

\ms

\it Proof. \rm By letting $m\to-\si_0$ in (\ref{6.33}) and
(\ref{6.34}), we have (\ref{6.41})--(\ref{6.42}). Thus, when
(\ref{6.41})--(\ref{6.42}) hold, for $m$ sufficiently closes to
$-\si_0$, (\ref{6.33})--(\ref{6.34}) hold.

\ms

(i) For the case $A^*=A$, we note that
$${1\over2}(A-m)^{-1}[e^{2(A-m)(T-t)}-I]=\int_{\si(A)}{e^{2(\m-m)(T-t)}-1
\over2(\m-m)}dE_\m.$$
By the definition of $\eta(\cd)$, we have that
$${e^{2(\m-m)(T-t)}-1\over2(\m-m)}=(T-t)\eta\big(2(\m-m)(T-t)\big),$$
and
$$\ba{ll}
\ns\ds\lim_{m\to-\si_0}{1\over2}(A-m)^{-1}[e^{A-m)(T-t)}-I]
=(T-t)\int_{\si(A)}\eta\big((\m+\si_0)(T-t)\big)dE_\m\\
\ns\ds\qq\qq\qq\qq\qq\qq\qq\qq=(T-t)\eta\big((A+\si_0)(T-t)\big).\ea$$
Then sending $m\to-\si_0$ in (\ref{PiB+BPi}), we obtain
(\ref{6.43}).

Now, for the case $A^*=-A$, by sending $m\to0$, we obtain
(\ref{6.44}). \endpf

\ms

Although the conditions stated in Theorem 7.5 and Corollary 7.6
still look lengthy, they are practically checkable. To illustrate
this, let us look at an interesting situation covered.

\ms

\bf Corollary 7.7. \sl Let {\rm(H0)} hold and
$(b,g,h)\in\cG_2\cap\cG_3$. Let
\bel{hy>0}h_y(y)+h_y(y)^*\ges0,\qq\forall y\in X.\ee
Let
\bel{7.48}\left\{\2n\ba{ll}
\ns\ds B_{11}(t,y,\psi)=B_{22}(t,y,\psi)^*,\\
\ns\ds B_{12}(t,y,\psi)+B_{12}^*(t,y,\psi)\les0,\\
\ns\ds
B_{21}(t,y,\psi)+B_{21}(t,y,\psi)^*\ges\d.\ea\right.\qq(t,y,\psi)\in[0,T]\times
X\times X,\ee
for some $\d>0$, and
\bel{7.49}\bar p_0(t)B_{22}+B_{22}^*\bar p_0(t)\les0,\qq\forall
t\in[0,T],\ee
where
\bel{}\bar p_0(t)=\left\{\2n\ba{ll}
\ns\ds\bar p_0e^{2(A+\si_0)t}+\bar q_0t\eta\big((A+\si_0)t\big),\qq
\hb{if }A^*=A,\\
\ns\ds[\bar p_0+\bar q_0t]I,\qq\qq\qq\qq\qq\q\hb{if
}A^*=-A,\ea\right.\ee
with $p_1,q_0,\bar p_0,\bar q_0>0$ and $\eta(\cd)$ is defined by
{\rm(\ref{eta})}. Then the FBEE generated by $(b,g,h)$ is
well-posed.

\ms

\it Proof. \rm First of all, by (\ref{hy>0}), and the boundedness of
$h_y(\cd)$ (since $(b,g,h)\in\cG_2\cap\cG_3$), we can find $p_1>0$
large so that (\ref{6.41}) holds, and $\g>0$ is allowed to be
arbitrarily large. Also, by letting $\th=0$, we see that
(\ref{6.42}) holds, as long as $q_0,\bar q_0>0$. Next, we define
\bel{}p_1(t)=\left\{\2n\ba{ll}
\ns\ds
p_1e^{2(A+\si_0)(T-t)}+q_0(T-t)\eta\big((A+\si_0)(T-t)\big),\q
\hb{if }A^*=A,\\
\ns\ds[p_1+q_0(T-t)]I,\qq\qq\qq\qq\qq\qq\q\hb{if }
A^*=-A.\ea\right.\ee
Then according to Corollary 7.6, the FBEE is well-posed if the
following holds:
$$\ba{ll}
\ns\ds\e\ges\begin{pmatrix}\sc p_1(t)&\sc\g I\\
\sc\g I&\sc-\bar p_0(t)\end{pmatrix}\begin{pmatrix}\sc B_{11}&\sc
B_{12}
\\ \sc-B_{21}&\sc-B_{22}
\end{pmatrix}+\begin{pmatrix}\sc B_{11}^*&\sc-B_{21}^*\\ \sc B_{12}^*&\sc-B_{22}^*\end{pmatrix}
\begin{pmatrix}\sc p_1(t)&\sc\g I\\
\sc\g I&\sc-\bar p_0(t)\end{pmatrix}\\
\ns\ds=\begin{pmatrix}\sc p_1(t)B_{11}-\g B_{21}&\sc p_1(t)B_{12}-\g
B_{22}\\ \sc\g B_{11}+\bar p_0(t)B_{21}&\sc\g B_{12}+\bar
p_0(t)B_{22}\end{pmatrix}+\begin{pmatrix}\sc B_{11}^*p_1(t)-\g
B_{21}^*&\sc\g
B_{11}^*+B_{21}^*\bar p_0(t)\\
\sc B_{12}^*p_1(t)-\g B_{22}^*&\sc\g B_{12}^*+B_{22}^*\bar p_0(t)\end{pmatrix}\\
\ns\ds=\begin{pmatrix}\sc
p_1(t)B_{11}+B_{11}^*p_1(t)-\g[B_{21}+B_{21}^*]&\sc
p_1(t)B_{12}+B_{21}^*\bar p_0(t)\\
\sc B_{12}^*p_1(t)+\bar p_0(t)B_{21}&\sc\bar
p_0(t)B_{22}+B_{22}^*\bar
p_0(t)+\g[B_{12}+B_{12}^*]\end{pmatrix}.\ea$$
This is equivalent to the following:
$$\begin{pmatrix}\sc\e+\g[B_{21}+B_{21}^*]-
[p_1(t)B_{11}+B_{11}^*p_1(t)]&\sc-[p_1(t)B_{12}+B_{21}^*\bar
p_0(t)]\\
\sc-[B_{12}^*p_1(t)+\bar
p_0(t)B_{21}]&\sc\e-\g[B_{12}+B_{12}^*]-[\bar
p_0(t)B_{22}+B_{22}^*\bar p_0(t)]\end{pmatrix}\ges0,$$
which is implied by
$$\g[B_{21}+B_{21}^*]-
[p_1(t)B_{11}+B_{11}^*p_1(t)]>0,$$
and
$$\ba{ll}
\ns\ds\e-\g[B_{12}+B_{12}^*]-[\bar p_0(t)B_{22}+B_{22}^*\bar
p_0(t)]\\
\ns\ds-[B_{12}^*p_1(t)+\bar p_0(t)B_{21}]\(\g[B_{21}+B_{21}^*]-
[p_1(t)B_{11}+B_{11}^*p_1(t)]\)^{-1}[p_1(t)B_{12}+B_{21}^*\bar
p_0(t)]\ges0.\ea$$
Note that
$$\ba{ll}
\ns\ds\Big\|\(\g[B_{21}+B_{21}^*]-
[p_1(t)B_{11}+B_{11}^*p_1(t)]\)^{-1}\Big\|\\
\ns\ds\les{1\over\g}\|(B_{21}+B_{21}^*)^{-1}\|\Big\|\(I-{1\over\g}(B_{21}+B_{21}^*)^{-1}
[p_1(t)B_{11}+B_{11}^*p_1(t)]\)^{-1}\Big\|\\
\ns\ds\les{1\over\g}\|(B_{21}+B_{21}^*)^{-1}\|{1\over1-{1\over\g}\|(B_{21}+B_{21}^*)^{-1}
[p_1(t)B_{11}+B_{11}^*p_1(t)]\|}\\
\ns\ds={\|(B_{21}+B_{21}^*)^{-1}\|\over\g-\|(B_{21}+B_{21}^*)^{-1}
[p_1(t)B_{11}+B_{11}^*p_1(t)]\|}.\ea$$
Therefore, it suffices to have (note (\ref{7.48})--(\ref{7.49}))
$$\ba{ll}
\ns\ds\e\ges{\|B_{12}^*p_1(t)+\bar
p_0(t)B_{21}\|^2\|(B_{21}+B_{21}^*)^{-1}\|\over\g-\|(B_{21}+B_{21}^*)^{-1}
[p_1(t)B_{11}+B_{11}^*p_1(t)]\|},\ea$$
which can be achieved by letting $\g>0$ sufficiently large. Then our
conclusion follows. \endpf

\ms

Let us look at some more cases.

\ms

\bf Corollary 7.8. \sl Let {\rm(H0)} hold. Suppose
$(b,g,h)\in\cG_2\cap\cG_3$ such that
\bel{>0}h_y(y)+h_y(y)^*\ges0,\qq\forall y\in X,\ee
and
\bel{<-d}\begin{pmatrix}0&I\\ I&0\end{pmatrix}\dbB(t,y,\psi)
+\dbB(t,y,\psi)^*\begin{pmatrix}0&I\\
I&0\end{pmatrix}\les-\d,\q\forall(t,y,\psi)\in[0,T]\times X\times
X,\ee
for some $\d>0$. Then the corresponding FBEE is well-posed.

\ms

\it Proof. \rm First, by letting $p_1,\bar p_0,q_0,\bar q_0>0$, with
$q_0,\bar q_0>0$ suitably large, we will have (\ref{6.42}). Then
letting $p_1>0$ large, we will have (\ref{6.41}). Next, by noting
$$0\les\eta(\k)\les1,\qq\forall\k\le0,$$
we see that under the condition $(b,g,h)\in\cG_4$, either $A^*=A$ or
$A^*=-A$, we always have the boundedness of all the terms involved
in the left-hand sides of (\ref{6.43}) and (\ref{6.44}),
respectively. Hence, under condition (\ref{<-d}), we can find $\g>0$
large enough so that (\ref{6.43}) and (\ref{6.44}) holds,
respectively. Due to the condition (\ref{>0}), by letting $\g>0$
large, (\ref{6.41}) will not be affected. Then Corollary 7.6 applies
to get the well-posedness of the corresponding FBEE. \endpf

\ms

Note that (\ref{<-d}) is equivalent to the following:
$$\ba{ll}
\ns\ds\begin{pmatrix}-[g_y(t,y,\psi)+g_y(t,y,\psi)^*]&
b_y(t,y,\psi)^*-g_\psi(t,y,\psi)\\
b_y(t,y,\psi)-g_\psi(t,y,\psi)^*&
b_\psi(t,y,\psi)+b_\psi(t,y,\psi)^*\end{pmatrix}\les-\d,\qq\forall(t,y,\psi)\in[0,T]\times
X\times X.\ea$$
This is further equivalent to the uniform monotonicity of the
following map
$$\begin{pmatrix}\sc y\\ \sc\psi\end{pmatrix}\mapsto\begin{pmatrix}\sc
g(t,y,\psi)\\\sc -b(t,y,\psi)
\end{pmatrix},$$
in the sense that for some $\d>0$,
$$\ba{ll}
\ns\ds\lan\begin{pmatrix}\sc g(t,y,\psi)-g(t,\bar y,\bar\psi)\\ \sc
-b(t,y,\psi)+b(t,\bar y,\bar\psi)\end{pmatrix},\begin{pmatrix}\sc y-\bar y\\
\sc\psi-\bar\psi\end{pmatrix}\ran\ges\d\big(\|y-\bar
y\|^2+\|\psi-\bar\psi\|^2),\qq\forall t\in[0,T],~y,\bar
y,\psi,\bar\psi\in X.\ea$$

\ms

It is possible to cook up many other cases from Theorem 7.5 and/or
Corollary 7.6, for which the corresponding FBEEs are well-posed. Let
us list some of them here.

\ms

\bf Corollary 7.9. \sl Let {\rm(H0)} hold and
$(b,g,h)\in\cG_2\cap\cG_3$. Then the corresponding FBEE is
well-posed if one of the following holds:

\ms

{\rm(i)} For some $\d,\e>0$,
$$I+h_y(y)+h_y(y)^*\ges\d,\qq\forall y\in X.$$
In the case $A^*=A$, for all $(t,y,\psi)\in[0,T]\times X\times X$,
$$\begin{pmatrix}e^{2(A+\si_0)(T-t)}& I\\
I&0\end{pmatrix}\dbB(t,y,\psi)+\dbB(t,y,\psi)^*\begin{pmatrix}
e^{2(A+\si_0)(T-t)}&I\\ I&0\end{pmatrix}\les0,$$
and in the case $A^*=-A$, for all $(t,y,\psi)\in[t,T]\times X\times
X$,
$$\begin{pmatrix}I&I\\
I&0\end{pmatrix}\dbB(t,y,\psi)+\dbB(t,y,\psi)^*\begin{pmatrix}I& I\\
I&0\end{pmatrix}\les0.$$

{\rm(ii)} For some $\d,\e>0$,
$$I-h_y(y)^*h_y(y)\ges\d,\qq\forall y\in X.$$
In the case $A^*=A$, for all $(t,y,\psi)\in[0,T]\times X\times X$,
$$\ba{ll}
\ns\ds\begin{pmatrix}e^{2(A+\si_0)(T-t)}&0\\ 0&
e^{2(A+\si_0)t}\end{pmatrix}\dbB(t,y,\psi)+\dbB(t,y,\psi)^*\begin{pmatrix}e^{2(A+\si_0)(T-t)}&0\\
0&e^{2(A+\si_0)t}\end{pmatrix}\les0,\ea$$
and in the case $A^*=-A$, for all $(t,y,\psi)\in[0,T]\times X\times
X$,
\bel{}\dbB(t,y,\psi)+\dbB(t,y,\psi)^*\le0.\ee

\ms

{\rm(iii)} For some $\d,\e>0$,
\bel{}I+h_y(y)+h_y(y)^*-h_y(y)^*h_y(y)\ges\d,\qq\forall y\in X.\ee
In the case $A^*=A$, for all $(t,y,\psi)\in[0,T]\times X\times X$,
\bel{}\ba{ll}
\ns\ds\begin{pmatrix}e^{2(A+\si_0)(T-t)}&I\\ I&
e^{2(A+\si_0)t}\end{pmatrix}\dbB(t,y,\psi)
+\dbB(t,y,\psi)^*\begin{pmatrix}e^{2(A+\si_0)(T-t)}&I\\
I&e^{2(A+\si_0)t}\end{pmatrix}\les0,\ea\ee
and in the case $A^*=-A$, for all $(t,y,\psi)\in[0,T]\times X\times
X$,
\bel{}\begin{pmatrix}I&I\\
I&I\end{pmatrix}\dbB(t,y,\psi)+\dbB(t,y,\psi)^*\begin{pmatrix}I& I\\
I&I\end{pmatrix}\les0.\ee

\ms

\it Proof. \rm (i) We take $p_1=\g>0$ large enough and take $\bar
p_0,q_0,\bar q_0>0$ small enough, $\th=0$. Then we may apply
Corollary 7.6 to get our claim.

\ms

(ii) and (iii) can be proved similarly. \endpf

\ms

Inspired by the above result, it is easy for us to prove many other
results of similar nature. We prefer not to get into exhausting
details.

\section{More General Cases}

In this section, we will briefly consider some more general cases.

\ms

First of all, we consider the case $(b,g,h)\in\cG_2$, i.e., the
generator $(b,g,h)$ only satisfies (H2), and might not be Fr\'echet
differentiable in $(y,\psi)$. Such a situation happens in many
optimal control problems. To study such a case, let us recall some
results from \cite{Pales-Zeidan 2007}.

\ms

Let $f:X\to X$ be Lipschitz continuous and $\bar y\in X$. For any
linear subspace $L\subseteq X$, we define {\it
$L$-G\^ateaux-Jacobian $D_Lf(\bar y)\in\cL(L;X)$} by the following
(if the limit exists):
$$D_Lf(\bar y)(x)=f'(\bar y;x)=\lim_{t\to0}{f(\bar y+tx)-f(\bar
y)\over t},\qq\forall x\in X.$$
The set of all points $\bar y\in X$ for which $D_Lf(\bar y)$ exist
is denoted by $\O_L(f)$. Next, we let
$$\pa_Lf(\bar y)=\bigcap_{\d>0}\coh\Big\{D_Lf(y)\bigm|y\in\O_L(f),\|y-\bar y\|\les\d\Big\}$$
and define the {\it generalized Jacobian} of $f(\cd)$ at $\bar y$ by
the following:
$$\pa f(\bar y)=\Big\{\Psi\in\cL(X;X)\bigm|\Psi\big|_L\in\pa_Lf(\bar
y),~\forall L\hb{ subspace of }X\Big\}.$$
For any $y,z\in X$, define $y\otimes z:\cL(X)\to\dbR$ by
$$(y\otimes z)(\Psi)=\lan\Psi(y),z\ran,\qq\forall\Psi\in\cL(X).$$
Then $y\otimes z\in\cL(\cL(X);\dbR)$. Let
$$X\otimes X=\span\{y\otimes z\bigm|y,z\in X\}\subseteq\cL(X)^*.$$
The weak topology induced by $X\otimes X$ on $\cL(X)$ is called the
{\it weak$^*$-operator-topology}, denoted by $\b(X)$. The following
can be found in \cite{Pales-Zeidan 2007}.

\ms

\bf Proposition 8.1. \sl If $f:X\to X$ is Lipschitz near $\bar y$,
then $\pa f(\bar y)$ is non-empty, bounded, and $\b(X)$-compact.

\ms

\rm

More interestingly, we have the following {\it mean-value theorem}
(see \cite{Pales-Zeidan 2007}, Theorem 4.4).

\ms

\bf Proposition 8.2. \sl Let $f:X\to X$ be locally Lipschitz. Then
for any $y,\bar y\in X$,
$$f(y)-f(\bar y)\in\[\coh\(\bigcup_{\l\in[0,1]}\pa f(\bar
y+\l(y-\bar y))\)\](y-\bar y).$$

\ms

\rm

With the above preparation, we now consider the case that
$(b,g,h)\in\cG_2$. Naturally, we need only to define
$$\dbB(s,y,\psi)=\begin{pmatrix}\sc b_y&\sc b_\psi\\
\sc-g_y&\sc-g_\psi\end{pmatrix},$$
with
$$\ba{ll}
\ns\ds b_y\in\pa_yb(t,y,\psi),\q b_\psi\in\pa_\psi b(t,y,\psi),\q
g_y\in\pa_y g(t,y,\psi),\q g_\psi\in\pa_\psi g(t,y,\psi).\ea$$
Then, all the results from previous sections for
$(b,g,h)\in\cG_2\cap\cG_3$ can be carried over properly to the case
$(b,g,h)\in\cG_2$.

\ms

Next, we consider $(b,g,h)\in\cG_3$, i.e., the generator $(b,g,h)$
may be not globally Lipschitz with respect to $y$ and/or $\psi$, or
equivalently, the Fr\'echet derivative of
$(y,\psi)\mapsto(b(t,y,\psi),g(t,y,\psi),h(y))$ might be not
bounded. In such cases, a priori uniform boundedness of the mild
solution $(y(\cd),\psi(\cd))$ could play an essential role. Let us
indicate one such a case. To this end, we introduce the following

\ms

{\bf(H3)$'$} In addition to (H3), there is a non-decreasing function
$f:[0,\infty)\to[0,\infty)$ such that
\bel{}\ba{ll}
\ns\ds\|b_y(t,y,\psi)\|+\|b_\psi(t,y,\psi)\|+\|g_y(t,y,\psi)\|
+\|g_\psi(t,y,\psi)\|+\|h_y(y)\|\les
f(\|y\|+\|\psi\|),\\
\ns\ds\qq\qq\qq\qq\qq\qq\qq\qq\qq\qq\forall(t,y,\psi)\in[0,T] \times
X\times X,\ea\ee
Moreover,
\bel{}\left\{\2n\ba{ll}
\ns\ds\lan b(t,y,\psi),y\ran\les L(1+\|y\|^2),\\
\ns\ds\lan g(t,y,\psi),\psi\ran\les
L(1+\|\psi\|^2),\ea\right.\qq\forall(t,y,\psi)\in[0,T]\times X\times
X.\ee

\ms

Under (H3)$'$, if $(y^\rho_\l(\cd),\psi^\rho_\l(\cd))$ is a solution
to (\ref{rho2}), then
$$\ba{ll}
\ns\ds\|y^\rho_\l(t)\|^2=\|x\|^2+2\int_0^t\lan y^\rho_\l(s),A_\l
y^\rho_\l(s)+\rho b(s,y^\rho_\l(s),\psi^\rho_\l(s))+b_0(s)
\ran ds\\
\ns\ds\qq\qq\les\|x\|^2+2L\int_0^t\(1+\|y^\rho_\l(s)\|^2+\|y^\rho_\l(s)\|\,
\|b_0(s)\|\)ds.\ea$$
Then by Gronwall's inequality, we have
$$\|y_\l^\rho(\cd)\|_\infty\les K\(1+\|x\|+\int_0^T\|b_0(r)\|dr\).$$
Similarly,
$$\ba{ll}
\ns\ds\|\psi^\rho_\l(t)\|^2=\|\psi^\rho_\l(T)\|^2
-2\int_t^T\lan\psi^\rho_\l(s),
-A^*_\l\psi^\rho_\l(s)-\rho g(s,y^\rho_\l(s),\psi^\rho_\l(s))-g_0(s)\ran ds\\
\ns\ds\qq\qq\les\|\psi^\rho_\l(T)\|^2+2L\int_t^T\(1+\|\psi^\rho_\l(s)\|^2
+\|\psi^\rho_\l(s)\|\,\|g_0(s)\|\)ds.\ea$$
Hence, it follows from Gronwall's inequality that
$$\ba{ll}
\ns\ds\|\psi^\rho_\l(\cd)\|_\infty\les
K\(1+\|\psi^\rho_\l(T)\|+\int_0^T
\|g_0(s)\|ds\)\\
\ns\ds\qq\qq\les
K\(1+\|h(0)\|+f(\|y^\rho_\l(T)\|)\|y^\rho_\l(T)+\|h_0\|
+\int_0^T\|g_0(s)\|ds\)\les K.\ea$$
Consequently, the relevant proofs will go through as if (H4) is assumed.

\ms

For concrete PDEs, there are some other ways to obtain uniform
boundedness of the (weak) solutions to the system. We will see some
of such below.

\ms

\section{Several Illustrative Examples}

In this section, we look at several examples.

\ms

\bf Example 9.1. (Linear-Convex Optimal Control Problem) \rm
Consider an optimal control problem with a linear state equation:
$$\left\{\2n\ba{ll}
\ns\ds\dot y(t)=Ay(t)+Bu(t),\\
\ns\ds y(0)=x,\ea\right.$$
and with the cost functional
$$J(x;u(\cd))=\int_0^T\(Q(y(t))+{1\over2}\lan Ru(t),u(t)\ran\)ds+G(y(T)),$$
where $y\mapsto Q(y)$ and $y\mapsto G(y)$ are $C^2$ and convex.
Then Pontryagin minimum principle leads to the optimality system:
\bel{FBEE(LC)}\left\{\2n\ba{ll}
\ns\ds\dot y(t)=Ay(t)-BR^{-1}B^*\psi(t),\\
\ns\ds\dot\psi(t)=-A^*\psi(t)-Q_y(y(t)),\\
\ns\ds y(0)=x,\qq\psi(T)=G_y(y(T)).\ea\right.\ee
In this case, we have
$$\ba{ll}
\ns\ds b(t,y,\psi)=-BR^{-1}B^*\psi,\q g(t,y,\psi)=Q_y(y),\q
h(y)=G_y(y).\ea$$
Thus,
$$\ba{ll}
\ns\ds b_y(t,y,\psi)=0,\qq b_\psi(t,y,\psi)=-BR^{-1}B^*,\\
\ns\ds g_y(t,y,\psi)=Q_{yy}(y),\qq g_\psi(t,y,\psi)=0,\qq
h_y(y)=G_{yy}(y).\ea$$
Then
$$\dbB(t,y,\psi)=\begin{pmatrix}\sc0&\sc-BR^{-1}B^*\\
\sc-Q_{yy}(y)&\sc0\end{pmatrix},$$
Hence, under conditions
$$R\ges\d,\qq M\ges G_{yy}(y)\ges0,\q M\ges Q_{yy}(y)\ges\d,\q\forall y\in X,$$
for some $M,\d>0$, all the conditions of Corollary 7.7 hold, and the
FBEE (\ref{FBEE(LC)}) admits a unique mild solution. A further
special case is the following:
$$Q(y)={1\over2}\lan Qy,y\ran,\q G(y)={1\over2}\lan Gy,y\ran,$$
for some $Q,G\in\dbS^+(X)$. In this case, the FBEE can be written as
\bel{FBEE(LQ)}\left\{\2n\ba{ll}
\ns\ds\dot y(t)=Ay(t)-BR^{-1}B^*\psi(t),\\
\ns\ds\dot\psi(t)=-A^*\psi(t)-Qy(t),\\
\ns\ds y(0)=x,\qq\psi(T)=Gy(T).\ea\right.\ee
Hence, according to the above, when
$$R\ges\d,\qq G\ges0,\qq Q\ges\d,$$
for some $\d>0$, the FBEE (\ref{FBEE(LQ)}) admits a unique mild solution.

\ms

\bf Example 9.2. (AQ Problem) \rm For the simplicity of
presentation, we let $S(\cd)=0$, and assume that all the involved
functions are time-independent. Then the optimality system reads
\bel{FBEE-AQ}\left\{\2n\ba{ll}
\ns\ds\dot
y(t)=Ay(t)+F(y(t))-BR^{-1}B^*\psi(t),\\
[1mm]
\ns\ds\dot\psi(t)=-A^*\psi(t)-Q_y(y(t))-F_y(y(t))^*\psi(t),\\
[1mm]
\ns\ds y(0)=x,\qq\psi(T)=G_y(y(T)),\ea\right.\ee
with $A^*=A$ or $A^*=-A$. Thus,
\bel{}\left\{\2n\ba{ll}
\ns\ds b(t,y,\psi)=F(y)-BR^{-1}B^*\psi,\\ [1mm]
\ns\ds g(t,y,\psi)=Q_y(y)+F_y(y)^*\psi,\\ [1mm]
\ns\ds h(y)=G_y(y).\ea\right.\ee
Let $\{\xi_n\}_{n\ge1}$ be an orthonormal basis of $X$, under which
we may let
$$F(y)=\sum_{n=1}^\infty\lan F(y),\xi_n\ran\xi_n\equiv\sum_{n=1}^\infty
f^n(y)\xi_n.$$
Then
$$F_y(y)z=\lim_{\d\to0}{F(y+\d
z)-F(y)\over\d}=\sum_{n=1}^\infty\lan f_y^n(y),z\ran\xi_n\equiv
\sum_{n=1}^\infty[\xi_n\otimes f_y^n(y)]z.$$
Thus,
$$F_y(y)=\sum_{n=1}^\infty\xi_n\otimes f_y^n(y),$$
and
$$F_y(y)^*\psi=\sum_{n=1}^\infty[f^n_y(y)\otimes\xi_n]\psi=\sum_{n=1}^\infty
f^n_y(y)\lan\xi_n,\psi\ran.$$
Hence,
$$[F_y(y)^*\psi]_y=\sum_{n=1}^\infty f^n_{yy}(y)\lan\xi_n,\psi\ran\in\dbS(X),
\qq\forall y\in X.$$
Consequently,
$$\left\{\2n\ba{ll}
\ns\ds b_y(s,y,\psi)=F_y(y),\q b_\psi(s,y,\psi)=-BR^{-1}B^*,\\
\ns\ds g_y(s,y,\psi)=Q_{yy}(y)+[F_y(y)^*\psi]_y,\q g_\psi(s,y,\psi)=F_y(y)^*=b_y(s,y,\psi)^*,\\
\ns\ds h_y(y)=G_{yy}(y).\ea\right.$$
Then
$$\dbB(s,y,\psi)=\begin{pmatrix}\sc F_y(y)&\sc-BR^{-1}B^*\\
\sc-Q_{yy}(y)-[F_y(y)^*\psi]_y&\sc-F_y(y)^*
\end{pmatrix}\equiv\begin{pmatrix}\sc B_{11}&\sc B_{12}\\ \sc-B_{21}&\sc-B_{22}\end{pmatrix}$$
From this, we can calculate
$$\ba{ll}
\ns\ds\begin{pmatrix}0&I\\
I&0\end{pmatrix}\dbB(s,y,\psi)+\dbB(s,y,\psi)^*
\begin{pmatrix}0&I\\ I&0\end{pmatrix}\\
\ns\ds=\begin{pmatrix}\sc-[g_y(s,y,\psi)+g_y(s,y,\psi)^*]&\sc b_y(s,y,\psi)^*-g_\psi(s,y,\psi)\\
\sc b_y(s,y,\psi)-g_\psi(s,y,\psi)^*&\sc
b_\psi(s,y,\psi)+b_\psi(s,y,\psi)^*\end{pmatrix}\\
\ns\ds=-2\begin{pmatrix}\sc Q_{yy}(y)+[F_y(y)^*\psi]_y&\sc0\\
\sc0&\sc BR^{-1}B^*\end{pmatrix}\ea$$
Next, we note that if $\psi(\cd)$ is a mild solution to the backward
evolution equation in (\ref{FBEE-AQ}), we have
$$\ba{ll}
\ns\ds\|\psi(t)\|^2=\|\psi(T)\|^2+2\int_t^T\lan[A+F_y(y(s))]\psi(s)+Q_y(y(s)),\psi(s)\ran
ds\\
\ns\ds\les\|G_y(y(T)))\|^2+2\int_t^T\|Q_y(y(s))\|\,\|\psi(s)\|ds\\
\ns\ds\les\|G_y(\cd)\|_\infty^2+2\int_t^T\|Q_y(\cd)\|_\infty\|\psi(s)\|ds\equiv\f(t).\ea$$
Then
$$\dot\f(t)=-2\|Q_y(\cd)\|_\infty\|\psi(t)\|\ges-2\|Q_y(\cd)\|_\infty\sqrt{\f(t)},$$
which leads to
$$\(\sqrt{\f(t)}\)'\ges-\|Q_y(\cd)\|_\infty.$$
Then
$$\ba{ll}
\ns\ds\|\psi(t)\|=\sqrt{\f(t)}=\sqrt{\f(T)}-\int_t^T\(\sqrt{\f(s)}\)'ds\les\|G_y(\cd)\|_\infty+\|Q_y(\cd)\|_\infty
T,\qq\forall t\in[0,T].\ea$$
Hence,
$$\ba{ll}
\ns\ds\|F_y(y)^*\psi\|=\|F_y(y)\|\,\|\psi\|\les\|F_y(y)\|\(\|G_y(\cd)\|_\infty
+\|Q_y(\cd)\|_\infty T\),\\
\ns\ds\qq\qq\qq\qq\qq\forall\,\|\psi\|\les\(\|G_y(\cd)\|_\infty
+\|Q_y(\cd)\|_\infty T\).\ea$$
Consequently, if we assume
$$\left\{\2n\ba{ll}
\ns\ds G_{yy}(y)\ges0,\q\forall y\in X,\qq BR^{-1}B^*\ges\d,\\
\ns\ds Q_{yy}(y)\ges\|F_y(y)\|\(\|G_y(\cd)\|_\infty
+\|Q_y(\cd)\|_\infty T\)+\d,\qq\forall y\in X,\ea\right.$$
for some $\d>0$, then (\ref{FBEE-AQ}) admits a unique mild solution,
by Corollary 7.8.

\ms

\bf Example 9.3. (Optimal Control of a Parabolic PDE). \rm We now
consider an optimal control problem for a parabolic equation. Such a
problem was studied in \cite{Stojanovic 1989}. The controlled state
equation reads:
\bel{}\left\{\2n\ba{ll}
\ns\ds y_t=\D y-(\l+u)y+f,\qq\hb{in }(0,T)\times\O,\\
\ns\ds y\big|_{\pa\O}=0,\\
\ns\ds y(0,x)=y_0(x),\qq x\in\O,\ea\right.\ee
where $y(t,x)$ is the state and $u(t,x)$ is the control, and $\O\subseteq
\dbR^n$ is a bounded domain with smooth boundary $\pa\O$. The cost
functional is the following:
\bel{}J(u(\cd))={1\over2}\int_0^T\int_\O\(L|y-y_d|^2+Nu^2\)dxdt+{1\over2}\int_\O
M|y(T,x)-z(x)|^2dx.\ee
We assume that
$$\left\{\2n\ba{ll}
\ns\ds f(t,x)\ges0,\q y_d(t,x)\les0,\qq(t,x)\in(0,T)\times\O,\\
[2mm]
\ns\ds y_0(x)\ges0,\q z(x)\les0,\qq x\in\O.\ea\right.$$
According to \cite{Stojanovic 1989}, optimal control exists and the
optimality system reads:
\bel{}\left\{\2n\ba{ll}
\ns\ds y_t=\D y-\l y-{1\over N}\psi y^2+f,\qq\hb{in }(0,T)\times\O,\\
\ns\ds\psi_t=-\D\psi+\l\psi+{1\over N}y\psi^2-L(y-y_d),\qq\hb{in }(0,T)\times\O,\\
\ns\ds y\big|_{\pa\O}=\psi\big|_{\pa\O}=0,\\ [2mm]
\ns\ds y(0,x)=y_0(x),\q\psi(T,x)=M\big(y(T,x)-z(x)\big),\qq
x\in\O,\ea\right.\ee
Then we have
$$\left\{\2n\ba{ll}
\ns\ds b(s,y,\psi)=-\l y-{1\over N}\psi y^2+f,\\ [2mm]
\ns\ds g(s,y,\psi)=-\l\psi-{1\over N}y\psi^2+L(y-y_d).\ea\right.$$
Hence,
$$\left\{\2n\ba{ll}
\ns\ds b_y=-\l-{2\over N}y\psi,\qq b_\psi=-{1\over N}y^2,\\
\ns\ds g_y=L-{1\over N}\psi^2,\qq g_\psi=-\l-{2\over
N}y\psi=b_y.\ea\right.$$
Then
$$\dbB(t,y,\psi)=\begin{pmatrix}\sc-\l-{2\over N}
y\psi&\sc-{1\over N}y^2\\ \sc-L+{1\over N}\psi^2&\sc\l+{2\over
N}y\psi\end{pmatrix}.$$
Thus,
$$\ba{ll}
\ns\ds\begin{pmatrix}\sc0&\sc I\\ \sc I&0\end{pmatrix}\dbB(t,y,\psi)
+\dbB(t,y,\psi)\begin{pmatrix}\sc0&\sc I\\ \sc
I&0\end{pmatrix}=2\begin{pmatrix}\sc-L+{1\over N}\psi^2&
\sc0\\
\sc0&\sc-{1\over N}y^2\end{pmatrix}\les0,\ea$$
provided $\psi$ is bounded (which was shown in \cite{Stojanovic
1989}) and $N$ is large enough.

\ms

\section{Concluding Remarks}

We have discussed the well-posedness of FBEEs which is mainly
motivated by the optimality systems of optimal control problems for
infinite dimensional evolution equations. We have presented some
basic results from two approaches: the decoupling method and the
method of continuity. It is seen that the theory is far from mature
and many challenging questions are left open. Here is a partial list
of these:

\ms

$\bullet$ In the direction of decoupling method, it is widely open
that how one can construct decoupling field, through solving a PDE
in Hilbert space.

\ms

$\bullet$ In the direction of method of continuity, more careful
analysis is need to make the stated condition easier to use.

\ms

$\bullet$ More general generators $A$ other than $A^*=A$ and
$A^*=-A$. Also, taking into account of PDEs, the generator $(b,g,h)$
might be unbounded (involving differential operators).

\ms

\bibliographystyle{amsplain}

\section{Appendix.} \rm

\it Proof of Proposition 4.2. \rm First, we have
$$\ba{ll}
\ns\ds\|\h y(s)\|^2\1n=\1n\|\h x\|^2\2n+\2n\int_t^s\2n\(2\lan\h
y(r),A_\l\h y(r)\1n+\1n\rho\wt b_y(r)\h y(r)\1n+\1n\rho\wt
b_\psi(s)\h\psi(s)
\1n+\1n\d b(r)\1n+\1n\d b_0(r)\ran\)dr\\
\ns\ds\les\1n\|\h
x\|^2\2n+\3n\int_t^s\2n\[\1n\lan\(A_\l+A_\l^*+\rho\big(\wt
b_y(r)\1n+\1n\wt b_y(r)^*\big)\)\h y(r),\h y(r)\ran\1n+2\|\h
y(r)\|\|\rho\wt b_\psi(r)\h\psi(r)\1n
+\1n\rho\d b(r)\1n+\1n\d b_0(r)\|\]dr\\
\ns\ds\les\|\h x\|^2\2n+2\2n\int_t^s\2n\[\rho L_{by}(r)\|\h
y(r)\|^2+\|\h y(r)\|\|\rho\wt b_\psi(r)\h\psi(r)\1n +\1n\rho\d
b(r)\1n+\1n\d b_0(r)\|\]dr\equiv\f(s)^2.\ea$$
Then
$$\ba{ll}
\ns\ds\f(s)\dot\f(s)=\rho L_{by}(s)\|\h y(s)\|^2+\|\h
y(s)\|\|\rho\wt b_\psi(s)\h\psi(s)\1n+\1n\rho\d b(s)\1n+\1n\d
b_0(s)\|\\
\ns\ds\qq\qq\les\rho L_{by}(s)\f(s)^2+\f(s)\|\rho\wt
b_\psi(s)\h\psi(s)\1n+\1n\rho\d b(s)\1n+\1n\d b_0(s)\|\\
\ns\ds\qq\qq\equiv a_2(s)\f(s)^2+a_1(s)\f(s).\ea$$
Consequently,
$$\dot\f(s)\les a_2(s)\f(s)+a_1(s).$$
Hence, by Gronwall's inequality,
$$\ba{ll}
\ns\ds\|\h y(s)\|\les\f(s)\le e^{\int_t^sa_2(\t)d\t}\|\h
x\|+\int_t^se^{\int_r^sa_2(\t)d\t}a_1(r)dr\\
\ns\ds=e^{\rho\int_t^sL_{by}(\t)d\t}\|\h
x\|+\int_t^se^{\rho\int_r^sL_{by}(\t)d\t}\|\rho\wt
b_\psi(r)\h\psi(r)\1n+\1n\rho\d b(r)\1n+\1n\d
b_0(r)\|dr\\
\ns\ds\les\rho\int_t^se^{\rho\int_r^sL_{by}(\t)d\t}\|\wt
b_\psi(r)\h\psi(r)\|dr+K\[\|\h x\|+\int_t^s\(\|\d b(r)\|+\|\d
b_0(r)\|\)dr\].\ea$$
This proves (\ref{|hy|}). We now prove (\ref{|hpsi|}). One has
$$\ba{ll}
\ns\ds\|\h\psi(T)\|^2-\|\h\psi(s)\|^2\\
\ns\ds=\2n\int_s^T2\lan\h\psi(r),-A_\l^*\h\psi(r)-\rho\wt g_y(r)\h
y(r)\1n-\1n\rho\wt g_\psi(r)\h\psi(r)-\1n\rho\d
g(r)\1n-\1n\d g_0(r)\1n\ran dr\\
\ns\ds\ges-2\int_s^T\3n\(\lan
\rho[\wt g_\psi(s)+\wt g_\psi(s)^*]\h\psi(s),\h\psi(s)\ran
+\lan\h\psi(s),\rho\wt g_y(s)\h y(s)+\rho\d g(s)+\d g_0(s)\ran\)ds\\
\ns\ds\ges2\int_t^T\Big\{-\rho
L_{g\psi}(r)\|\h\psi(r)\|^2-\|\h\psi(r)\|\,\|\rho\wt g_y(r)\h
y(r)+\rho\d g(r)+\d g_0(r)\|\Big\}ds,\ea$$
which leads to
$$\ba{ll}
\ns\ds\|\h\psi(s)\|^2\les\|\h\psi(T)\|^2+2\int_s^T\Big\{\rho
L_{g\psi}(r)\|\h\psi(r)\|^2+\|\h\psi(r)\|\,\|\rho\wt g_y(r)\h
y(r)+\rho\d g(r)+\d g_0(r)\|\Big\}dr\equiv\f(s)^2.\ea$$
Then
$$\ba{ll}
\ns\ds\f(s)\dot\f(s)=-\rho
L_{g\psi}(r)\|\h\psi(s)\|^2+\|\h\psi(s)\|\,\|\rho\wt
g_y(s)\h y(s)+\rho\d g(s)+\d g_0(s)\|\\
\ns\ds\qq\qq\ges-\rho L_{g\psi}(r)\f(s)^2+\|\rho\wt g_y(s)\h
y(s)+\rho\d g(s)+\d
g_0(s)\|\f(s)\equiv-a_2(s)\f(s)^2-a_1(s)\f(s).\ea$$
Thus,
$$\dot\f(s)+a_2(s)\f(s)\ges-a_1(s).$$
$$\(e^{-\int_s^Ta_2(\t)d\t}\f(s)\)'\ges-a_1(s)e^{-\int_s^Ta_2(\t)d\t}$$
$$\f(T)-e^{-\int_s^Ta_2(\t)d\t}\f(s)\ges-\int_s^Ta_1(r)e^{-\int_r^Ta_2(\t)d\t}dr.$$
Hence,
$$\ba{ll}
\ns\ds\|\h\psi(s)\|\les\f(s)\les
e^{\int_s^Ta_2(\t)d\t}\f(T)+\int_s^Ta_1(r)e^{\int_s^ra_2(\t)d\t}dr\\
\ns\ds\les e^{\rho\int_s^TL_{g\psi}(\t)d\t}\|\wt h_y\h y(T)+\rho\d
h+\d h_0\|+\int_s^Te^{\rho\int_s^rL_{g\psi}(\t)d\t}\|\rho\wt
g_y(r)\h y(r)+\rho\d g(r)+\d g_0(r)\|dr\\
\ns\ds\les\rho\[e^{\rho\int_s^TL_{g\psi}(\t)d\t}\|\wt h_y\h
y(T)\|+\int_s^Te^{\rho\int_s^rL_{g\psi}(\t)d\t}\|\wt g_y(r)\h
y(r)\|dr\]\\
\ns\ds\qq+K\[\|\d h\|+\|\d h_0\|+\int_s^T\(\|\d g(r)\|+\|\d
g_0(r)\|\)dr\].\ea$$
This completes the proof. \endpf

\bs

Note that if we let $(y^\rho_0(\cd),\psi^\rho_0(\cd))$ be the mild
solution of the following:
\bel{rho0}\left\{\2n\ba{ll}
\ns\ds\3n\ba{ll}\dot y_0^\rho(s)=Ay_0^\rho(s)+\rho b(s,0,0)+b_0(s),\\
\ns\ds\dot\psi_0^\rho(s)=-A^*\psi_0^\rho(s)-\rho g(s,0,0)-g_0(s),\ea\qq s\in[t,T],\\
\ns\ds y_0^\rho(t)=x,\qq\psi^\rho_0(T)=\rho h(0)+h_0.\ea\right.\ee
Then
$$y^\rho_0(s)=e^{A(s-t)}x+\int_t^se^{A(s-r)}[\rho
b(r,0,0)+b_0(r)]dr.$$
Hence,
$$\|y^\rho_0(\cd)\|_\infty\les\|x\|+\int_t^T\|\rho
b(r,0,0)+b_0(r)\|dr.$$
Also,
$$\psi^\rho_0(s)=e^{A^*(T-s)}\psi^\rho_0(T)+\int_s^Te^{A^*(r-s)}[\rho g(r,0,0)+g_0(r)]dr.$$
Hence,
$$\|\psi^\rho_0(\cd)\|_\infty\les\|\rho h(0)+h_0\|+\int_t^T\|\rho
g(r,0,0)+g_0(r)\|dr.$$
Now, taking
$$\ba{ll}
\ns\ds\bar b(s,y,\psi)=b(s,0,0),\q\bar b_0(s)=b_0(s),\\
\ns\ds\bar g(s,y,\psi)=g(s,0,0),\q\bar g_0(s)=g_0(s),\\
\ns\ds\bar h(y)=h(0),\qq\bar h_0=h_0.\ea$$
Then
\bel{}\left\{\2n\ba{ll}
\ns\ds\d b(s)=b(s,0,0)-b(s,y_\l^\rho(s),\psi_\l^\rho(s)),\q\d
b_0(s)=\bar
b_0(s)-b_0(s)=0,\\
\ns\ds\d g(s)=g(s,0,0)-g(s,y_\l^\rho(s),\psi_\l^\rho(s)),\q\d
g_0(s)=\bar
g_0(s)-g_0(s)=0,\\
\ns\ds\d h=h(0)-h(y_\l^\rho(T)),\q\d h_0=\bar h_0-h_0=0,\q\h
x=x-\bar x=0.\ea\right.\ee
\bel{|hy|*}\ba{ll}
\ns\ds\|\h
y(\cd)\|_\infty\les\rho\int_t^Te^{\rho\int_r^TL_{by}(\t)d\t}\|\wt
b_\psi(r)\h\psi(r)\|dr.\ea\ee
and
\bel{|hpsi|*}\ba{ll}
\ns\ds\|\h\psi(\cd)\|_\infty\les\rho\[e^{\rho\int_t^TL_{g\psi}(\t)d\t}\|\wt
h_y\h y(T)\|+\int_t^T\3n e^{\rho\int_t^rL_{g\psi}(\t)d\t}\|\wt
g_y(r)\h y(r)\|dr\].\ea\ee

\end{document}